\documentclass[preprint,10pt]{elsarticle}

\usepackage[margin=0.8in]{geometry}

\usepackage[usenames,dvipsnames]{xcolor}

\let\today\relax
\makeatletter
\def\ps@pprintTitle{%
    \let\@oddhead\@empty
    \let\@evenhead\@empty
    \def\@oddfoot{\footnotesize\itshape
         {} \hfill\today}%
    \let\@evenfoot\@oddfoot
    }
\makeatother

\usepackage{amssymb,amsfonts,amsmath}
\usepackage{mathtools} % shawn
\usepackage{color}
\usepackage{graphicx,float}

\usepackage{booktabs}
\usepackage{dcolumn,cancel}
\newcolumntype{d}[1]{D{.}{.}{#1}}  % define "d" column type
\usepackage[normalem]{ulem}
\usepackage{accents}

\newtheorem{remark}{Remark}
\newtheorem{claim}{Claim}

\newcommand{\eps}{\epsilon}
\newcommand{\bv}{\boldsymbol{v}}

\newcommand{\be}{\boldsymbol{e}}

\newcommand{\dd}{\,\textrm{d}}
\newcommand{\bR}{\boldsymbol{R}}

\newcommand{\SSS}{S}

\newcommand{\bk}{{\boldsymbol{k}}}

\newcommand{\bi}{{\boldsymbol{i}}}

\newcommand{\bmu}{{\boldsymbol{\mu}}}

\newcommand{\bchi}{{\boldsymbol{\chi}}}

\newcommand{\bx}{\boldsymbol{x}}
\newcommand{\bX}{\boldsymbol{X}}

\newcommand{\bg}{\boldsymbol{g}}
\newcommand{\bh}{\boldsymbol{h}}
\newcommand{\br}{\boldsymbol{r}}
\newcommand{\bc}{\boldsymbol{c}}
\newcommand{\bth}{\boldsymbol{\theta}}
\newcommand{\btau}{\boldsymbol{\tau}}

\newcommand{\bnu}{\boldsymbol{\nu}}

\newcommand{\hh}{\hspace*{0.7pt}}
\newcommand{\nes}{\hspace*{-0.2pt}}

\newcommand{\chsig}{\check{\sigma}}
\newcommand{\hasig}{\hat{\sigma}}
\newcommand{\sip}{\!\cdot\!}

\newcommand{\tz}{\mbox{\tiny{(0)}}}
\newcommand{\tmm}{\mbox{\tiny{(m)}}}
\newcommand{\so}{\mbox{\tiny{(1)}}}
\newcommand{\st}{\mbox{\tiny{(2)}}}
\newcommand{\sh}{\mbox{\tiny{(3)}}}

\newcommand{\tn}{\mbox{\tiny{N}}}
\newcommand{\td}{\mbox{\tiny{D}}}

\newcommand{\tnn}{\mbox{\tiny{(N)}}}
\newcommand{\tdd}{\mbox{\tiny{(D)}}}
\newcommand{\tpp}{\prime}
\newcommand{\bpp}{\mbox{\tiny{[p]}}}

\newcommand{\bu}{\boldsymbol{u}}

\begin{document}
\begin{frontmatter}
%
%%% Title, authors and addresses
%
%%% use the tnoteref command within \title for footnotes;
%%% use the tnotetext command for theassociated footnote;
%%% use the fnref command within \author or \address for footnotes;
%%% use the fntext command for theassociated footnote;
%%% use the corref command within \author for corresponding author footnotes;
%%% use the cortext command for theassociated footnote; % \mathcal{K}_\eps = \{\hat\bk\in\mathbb{R}^d: \eps\hat\bk\in\mathbb{C}, ~~ \omega_n(\bk_s+\eps\hat\bk)=\omega ~\Leftrightarrow~ \bth^{\tz}\sip (i\hat{\bk})+\chsig \rho^{\tz}\check{\omega}^2=0\}.
% \]

%%% use the ead command for the email address,
%%% and the form \ead[url] for the home page:
%%% \title{Title\tnoteref{label1}}
%%% \tnotetext[label1]{}
%%% \author{Name\corref{cor1}\fnref{label2}}
%%% \ead{email address}
%%% \ead[url]{home page}
%%% \fntext[label2]{}
%%% \cortext[cor1]{}
%%% \address{Address\fnref{label3}}
%%% \fntext[label3]{}
%
\title{On the spectral asymptotics of waves in periodic media with Dirichlet or Neumann exclusions }

\author[1]{Othman Oudghiri-Idrissi} \author[1]{Bojan B. Guzina \corref{BG}} \author[2]{Shixu Meng}
%\renewcommand{\thefootnote}{\fnsymbol{footnote}}
%\footnotetext[1]{Department of Civil, Environmental, \& Geo- Engineering, University of Minnesota, Twin Cities,
% Minneapolis, MN 55405, USA. {\tt guzin001@umn.edu} and {\tt oudgh001@umn.edu}}
%\footnotetext[2]{Department of Mathematics, University of Michigan,
%  Ann Arbor, MI 48109, USA.  {\tt shixumen@umich.edu}}

%% use optional labels to link authors explicitly to addresses:
%% \author[label1,label2]{}
\address[1]{Department of Civil, Environmental \& Geo- Engineering, University of Minnesota, Twin Cities, Minneapolis, MN 55405, USA. {\tt oudgh001@umn.edu} and {\tt guzin001@umn.edu} }
\address[2]{Institute of Applied Mathematics, Academy of Mathematics and Systems Science, Chinese Academy of Sciences, Beijing, 100190, China. {\tt shixumeng@amss.ac.cn }}

\cortext[BG]{Corresponding author}

\begin{abstract}
We consider homogenization of the scalar wave equation in periodic media at finite wavenumbers and frequencies, with the focus on continua characterized by: (a) arbitrary Bravais lattice in~$\mathbb{R}^d$, $d\!\geqslant\!2$, and (b) exclusions i.e.~``voids'' that are subject to homogenous (Neumann or Dirichlet) boundary conditions. Making use of the Bloch wave expansion, we pursue this goal via asymptotic ansatz featuring the ``spectral distance'' from a given wavenumber-eigenfrequency pair (situated \emph{anywhere} within the first Brillouin zone) as the perturbation parameter. We then introduce the effective wave motion via projection(s) of the scalar wavefield onto the Bloch eigenfunction(s) for the unit cell of periodicity, evaluated at the origin of a spectral neighborhood. For generality, we account for the presence of the source term in the wave equation and  we consider -- at a given wavenumber -- generic cases of isolated, repeated, and nearby eigenvalues. In this way we obtain a palette of effective models, featuring both wave- and Dirac-type behaviors, whose applicability is controlled by the local band structure and eigenfunction basis. In all spectral regimes, we pursue the homogenized description up to at least first order of expansion, featuring asymptotic corrections of the homogenized Bloch-wave operator and the homogenized source term. Inherently, such framework provides a convenient platform for the synthesis of a wide range of intriguing wave phenomena, including negative refraction and  topologically protected states in metamaterials and phononic crystals. The proposed homogenization framework is illustrated by approximating asymptotically the dispersion relationships for (i) Kagome lattice featuring hexagonal Neumann exclusions, and (ii) ``pinned'' square lattice with circular Dirichlet exclusions. We complete the numerical portrayal of analytical developments by studying the response of a Kagome lattice due to a dipole-like source term acting near the edge of a band gap.
\end{abstract}

%%Research highlights
%\begin{highlights}
%\item Research highlight 1
%\item Research highlight 2
%\end{highlights}

\begin{keyword}
Waves in periodic media \sep dynamic homogenization \sep finite wavenumber \sep finite frequency \sep phononic crystals \sep acoustic metamaterials
\end{keyword}

\end{frontmatter}
\bibliographystyle{elsarticle-num} 
\bibliography{<your bibdatabase>}
\bibliographystyle{plain}

%\maketitle

%--------------------------------------------------------------------------------------------------------------%
\section{Introduction} \label{Introduction}
%--------------------------------------------------------------------------------------------------------------%

\noindent Dynamic homogenization of periodic media such as composites, phononic crystals, and metamaterials serves to (i)~deepen insight into the underpinning physical phenomena such as wave directivity, stop bands, and negative refraction \cite{ACG14,LZMZYCS2000,W16}, (ii)~reduce the burden of multi-scale numerical simulations, and (iii)~aid the ``microstructural'' design catering for applications such us cloaking \cite{Capo2009}, vibration control \cite{Bava2013}, logic circuits \cite{MSG19}, or seismic shielding~\cite{AABCEG17}. To survey the lay of the land in terms of homogenization strategies, it is generally convenient to distinguish between competing frequency and wavelength regimes. 

In the low-wavenumber, low-frequency (LW-LF) regime, one assumes the separation in scale between some finite wavelength and the vanishing lengthscale of medium periodicity, which then provides a perturbation parameter for the two-scale homogenization method \cite{ABDW08, CF01, PBL78, WG15}  and the Bloch-wave expansion (BWE) approach \cite{ABV16,SS91,W78} towards establishing a macroscopic i.e.~effective description of the medium. In this regime, a non-asymptotic effective medium model introduced by Willis \cite{JRW83,MW07,W11}, that links in coupled form the momentum and stress to particle velocity and strain, has also gained notable traction in the literature~\cite[e.g.][]{MG17,W16}. 

In principle, the two- (or more generally multiple-) scale homogenization framework applies equally to the low-wavenumber, high-frequency (LW-HF) regime where the lengtscale ``$\ell$'' of medium periodicity vanishes in the limit, while the dominant wavelength is kept at~$O(1)$ irrespective of the vibration frequency. In this case, eigenfrequencies of the unit cell problem grow as $O(\ell^{-1})$, which justifies the ``high frequency'' designation. In this vein, Allaire and Conca~\cite{Alla1998} introduced the Bloch wave homogenization method, which is essentially a combination of BWE and two-scale convergence analysis~\cite{Alla1992}. The subject of LW-HF homogenization via multiple-scale expansion has since been pursued in a number of studies, with applications to e.g. Maxwell equations~\cite{Sjob2005} and Navier equations~\cite{Auri2012}.

In the finite-wavenumber, finite-frequency (FW-FF) regime, the wavelength and the unit cell may be of the same order. For this class of problems, the ``high frequency homogenization'' (HFH) treatment introduced by Craster and co-workers \cite{CKP10} loosely exploits the small size of the unit cell relative to dimensions of the periodic domain to establish a two-scale analysis of the homogenous wave equation in periodic media. The method yields zeroth- i.e. leading-order effective description of the medium, in the vicinity of simple and multiple eigenfrequencies, that is shown to capture dynamic phenomena such as anisotropic wave motion \cite{ACG14,CAAMCSTECG15} and all-angle negative refraction \cite{Cra11}. In the follow-up studies, the HFH has extended to deal with zero-frequency stop-band media, i.e.~those with Dirichlet inclusions \cite{ACG13b,MAMGC16}, and to periodic media with Neumann inclusions \cite{ACG13a}. From the mathematical viewpoint, the existence of FW-FF effective differential operators was formally established by Birman and Suslina~\cite{Birm2004,Birm2006}, who considered the behavior of periodic systems near the edge of an internal band gap. On the engineering side, an effort was also made to extend the Willis' homogenization approach to finite frequencies and wavelengths~\cite{NHA16a} by introducing additional kinematic degrees of freedom; however the uniqueness of such effective model remains an open question. Recently, a systematic framework for homogenization of the scalar wave equation in the FW-FF regime was proposed in~\cite{GMO19}. This study, catering for rectangular lattices, makes use of PWE to secure a ``tight handle'' on the wavenumber, and defines the perturbation parameter as a distance  in the frequency-wavenumber space in order to  obtain effective medium description in the vicinity of simple, repeated, and and nearby eigenfrequencies. As a means to deal with finite (non-zero) wavenumbers, the authors in~\cite{GMO19} make use of the so-called multicell technique \cite{Doss2008,Gone2010}, which essentially restricts the applicability of their model to the apexes of the first Brillouin zone and its quadrants in~$\mathbb{R}^2$ (resp. octants in~$\mathbb{R}^3$). In essence, such restriction on the wavenumber ensures that the germane Bloch eigenfunctions have constant phase over the unit cell, which lends itself to an immediate proof that the effective medium properties are real-valued there. Unlike the HFH approach, this study also includes a homogenization treatment of the source term, which (as it turns out) is also subject to asymptotic corrections; see also~\cite{ABV16,MG17} in the context of LW-LF homogenization.  

In the present work, we generalize the FW-FF homogenization framework in~\cite{GMO19} to enable treatment of periodic media in $\mathbb{R}^d$, $d\!\geqslant\!2$ that are supported by generic Bravais lattice and may contain ``voids'' that are subject to homogenous Neumann and/or Dirichlet boundary conditions. The scalar wave equation under consideration bears relevance, for instance, to the description of anti-plane shear waves in two-dimensional composites, transverse electric (TE) or transverse magnetic (TM) fields in photonic crystals, and acoustic waves in three-dimensional phononic crystals. In a departure from the preceding work, our analysis further (a) permits expansion about an \emph{arbitrary} wavenumber-eigenfrequency pair within the first Brillouin zone; (b) allows for spatially-varying (as opposed to constant) Bloch-wave representation of the source term, and (c) pursue the homogenized description near simple, repeated, and nearby eigenvalues up to at least the first order of asymptotic correction. We illustrate the proposed homogenization framework through the study of (i) Kagome lattice featuring hexagonal Neumann exclusions, and (ii) ``pinned'' square lattice with circular Dirichlet exclusions. We complete the numerical portrayal by studying the response of a Kagome lattice due to a dipole-like source term acting near the edge of an internal band gap.

%--------------------------------------------------------------------------------------------------------------%
\section{Preliminaries} \label{prelim}
%--------------------------------------------------------------------------------------------------------------%
%--------------------------------------------------------------------------------------------------------------%
\subsection{Geometry} \label{geo}
\begin{figure}[h!] 
\centering{\includegraphics[width= 130 mm]{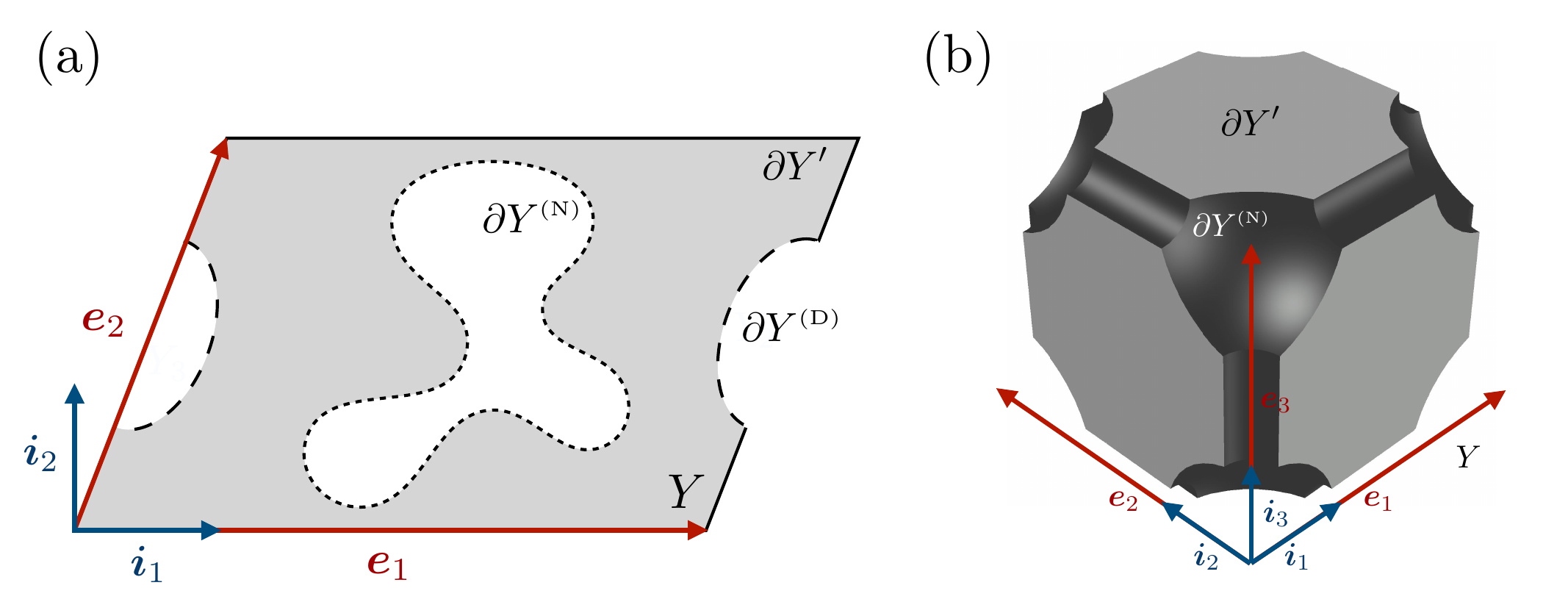}}\vspace*{-5mm}
\caption{Examples of the unit cell indicating the lattice basis vectors $\be_j$ as well as boundaries $\partial Y^{\tpp}$, $\partial Y^{\tnn}$ and $\partial Y^{\tdd}$ in: (a) $\mathbb{R}^2$, and (b) $\mathbb{R}^3$.} \label{fig1} \vspace*{-5mm}
\end{figure}

%--------------------------------------------------------------------------------------------------------------%
\noindent Consider an infinite periodic medium in $\mathbb{R}^{d}$ ($d \geqslant2$) affiliated with a Bravais lattice
\begin{equation}\label{Brav}
 \boldsymbol{R} \;=\; \big\{\sum_{j=1}^d n^j \be_j\!: ~ n^j\in\mathbb{Z}\big\},  
\end{equation}
featuring the basis $\be_j\in\mathbb{R}^d$, $j=\overline{1,d}$. Letting hereon $j\in\overline{1,d}$ implicitly unless stated otherwise, we denote by~$x^{j}$ the contravariant components of the position vector $\bx\in\mathbb{R}^d$ with reference to the lattice basis $\be_j$, and by $r^j\in\mathbb{Z}$ the contravariant coordinates of the lattice point $\br\in\boldsymbol{R}$. Next,  let
\[
Y_0 \,=\, \{\bx\!: 0 < x^j < 1\}
\]
denote the ``elemental parallelepiped'' of the lattice attached to the origin, and let $Y^{\tn}, Y^{\td}\subset Y_0$ denote a pair of disjointed open sets, each representing a union of 1-connected sets as illustrated in Fig.~\ref{fig1}. With such definitions, one may define the support of the periodic medium as  
 \begin{equation}\label{S}
\SSS \;=\; \mathbb{R}^d\setminus \bigcup_{\br\in\boldsymbol{R}} \big(\br + \overline{Y^{\tn}\cup Y^{\td}}\big), 
\end{equation}
 \noindent whose unit cell of periodicity is given by
 \begin{equation}\label{Y}
Y \;=\; Y_0\cap \SSS.
\end{equation}
Here it is useful to observe that~$Y$ is connected set thanks to the foregoing restrictions on~$Y^{\tn}$ and~$Y^{\td}$. We further define the domain boundaries $\partial Y^{\tnn}$, $\partial Y^{\tdd}$, $\partial Y^{\tpp}$, $\partial\SSS^{\tn}$ and $\partial\SSS^{\td}$ respectively as
\begin{eqnarray}\label{SNSD} 
&\partial Y^{\tnn} = \partial Y \cap \partial Y^{\tn}, \quad \partial Y^{\tdd} = \partial Y \cap \partial Y^{\td}, \quad \partial Y^{\tpp} = \partial Y \setminus (\partial Y^{\tnn}\cup \partial Y^{\tdd}),  \nonumber \\
&\partial\SSS^{\tn}=\bigcup_{\br\in\boldsymbol{R}} (\br+\partial Y^{\tnn}), \quad \partial\SSS^{\td}=\bigcup_{\br\in\boldsymbol{R}} (\br+\partial Y^{\tdd}),
\end{eqnarray}
such that $\partial Y = \partial Y^{\tpp}\cup\partial Y^{\tnn} \cup \partial Y^{\tdd}$ and $\partial\SSS = \partial\SSS^{\tn}\cup\partial\SSS^{\td}$. In physical terms, subtraction of~$\overline{Y^{\tn} \cup Y^{\td}}$ from~$Y$ in~\eqref{S} accounts for the ``holes'' featured by the periodic medium. Accordingly, the boundary conditions defined on $\partial\SSS^{\tn}$ and $\partial Y^{\tnn}$ (resp. $\partial\SSS^{\td}$ and $\partial Y^{\tdd}$) are considered only if $Y^{\tn}$ (resp. $Y^{\td}$) is a nonempty set. Examples of 2D and 3D unit cells geometries as defined above are illustrated in Fig. \ref{fig1}.  To facilitate the analysis, we will also make use of the short-hand notation
\begin{equation}\label{perbcs}
\partial{Y^{\tpp}_{jm}} = \{\bx\in\partial Y^{\tpp}\!:\; x^{j}=m \}, \quad m=\overline{0,1}.
\end{equation}

For further reference, we denote by $\be^{j}\in\mathbb{R}^d$ the covariant lattice basis that spans the reciprocal space $\mathbb{R}^d$ and satisfies $\be^j \sip \be_i = 2\pi\delta_{ij}$ ($i,j=\overline{1,d}$) where $\delta_{ij}$~is the Kronecker delta; by
\begin{equation}\label{recBrav}
\boldsymbol{R}^* \:=\; \big\{\sum_{j=1}^d n_j \be^j: ~ n_j\in\mathbb{Z}\big\} 
\end{equation}
the reciprocal Bravais lattice; by $k_{j}$ the covariant components of wave vector $\bk\in\mathbb{R}^d$ tied to the basis $\be^j$; by $r^*_j\in\mathbb{Z}$ the covariant coordinates of the lattice point $\br^*\in\boldsymbol{R}^*$; by $Y_0^*$ the reciprocal of~$Y_0$ defined by
 \begin{equation}\label{Y_0s} \notag
Y_0^* \,=\, \{\bk\!:0 < k_j < 1 \},
\end{equation}
and by 
\begin{equation}\label{Brillouin}
\mathcal{B} = \big\{\bk\in\mathbb{R}^d\!: \;\bk\cdot\boldsymbol{\kappa}\leqslant\tfrac{1}{2}\|\boldsymbol{\kappa}\|^2,~\boldsymbol{\kappa} \!=\! \sum_{j=1}^{d}n_j\be^j, \; n_j\!\in\{-1,0,1\}\big\}
\end{equation}
the first Brillouin zone of the lattice. We also denote by $|\mathcal{D}|$ the volume of a finite domain $\mathcal{D}\subset\mathbb{R}^d$ and by $v=(1-|Y||Y_0|^{-1})$ the \emph{porosity} of periodic medium $\SSS$. With such definitions, one may note that $|\mathcal{B}| = |Y_0^*|=(2\pi)^{d}|Y_0|^{-1}$.

%--------------------------------------------------------------------------------------------------------------%
\subsection{Function spaces} \label{funcSpac}
%--------------------------------------------------------------------------------------------------------------%

\noindent In what follows, we will deal with mappings of type $g:D\mapsto\mathbb{C}$ for some~$D\subseteq\mathbb{R}^d$, and their tensorial generalizations. To help set up the table for discussion, we first define the space 
\begin{eqnarray} \label{spaces1} 
 L^2(D) ~=~ \big\{ g\!:\! \int_D g\, \overline{g} \dd\bx <\infty \big\},  
\end{eqnarray}
and we introduce a generalized inner product over~$Y$, namely 
\begin{equation} \label{innerp}
(\bg,\bh) \;=\; |Y|^{-1} \int_{Y} \bg \!:\! \overline{\bh} \dd\bx,\qquad \bg \!\in (L^2(Y))^{d^m}, ~\bh \!\in (L ^2(Y))^{d^n}, \quad m,n\geqslant 0,
\end{equation}
where
\[
(L ^2(Y))^{d^m} \;=\; (L ^2(Y))^{d\times d\times\ldots\times d} \quad (m \text{ times}),
\]
and  ``:'' stands for the usual product, the inner product, and the $\min(m,n)$-tuple tensor contraction when $m=n=0$, $m=n=1$, and $\max(m,n)>1$, respectively. In situations when $g:S\mapsto\mathbb{C}$ is $Y$-periodic, we will also make use of the periodic function spaces 
\begin{eqnarray} \label{spaces} 
L^2_p(Y) &=& \{g: g|_{Y}\!\in\! L^2(Y),~ g(\bx\!+\!\br)=g(\bx) ~~\forall\hh \bx\!\in\!S,\, \br\!\in\!\bR\}, \notag \\ 
L^2_{p0}(Y) &=& \{ g \!\in\! L^2_p(Y): ~g|_{\partial Y^{\tdd}}=0 \}, \\
H^1_{p0}(Y) &=& \{ g \in L^2_{p0}(Y): \nabla  g\in(L ^2(Y))^d \}. \notag
\end{eqnarray}

%--------------------------------------------------------------------------------------------------------------%
\subsection{Wave equation, boundary conditions, and Bloch wave expansion} \label{we}
%--------------------------------------------------------------------------------------------------------------%
\noindent Consider the time-harmonic wave equation in~$\SSS$ at frequency~$\omega$, namely 
\begin{equation}\label{PDE}
-\omega^2\rho(\bx)\hh u - \nabla\sip \big(G(\bx)\nabla u\big) \;=\; f(\bx) \qquad \text{in~}\SSS, 
\end{equation}
where $f\!\in\! L^2(S).$ In what follows, we assume the coefficients~$0<G<\infty$ and~$0<\rho<\infty$ to be $Y$-periodic and bounded away from zero. In this case, we note that $L^2_{p0}(Y)$ equivalently defines the class of ``kinematically-compatible'' periodic functions satisfying~$(\rho g,g)<\infty$. To complete the formulation of the problem, we assume that~$u$ satisfies homogeneous Neumann and Dirichlet boundary conditions on $\partial\SSS^{\tn}$ and $\partial\SSS^{\td}$, respectively. In other words, we let 
\begin{eqnarray}\label{BC}
\bnu\cdot\!G\hh\nabla u &\!\!=\!\!& 0  \qquad \text{on~} \partial \SSS^{\tn}, \\
u &\!\!=\!\!& 0 \qquad \text{on~} \partial \SSS^{\td}, \label{BCD}
\end{eqnarray}
where $\bnu$ is the unit outward normal on $\partial\SSS^{\tn}$. For generality, we note that~\eqref{PDE} pertains to a wide class of physical processes including: (i) anti-plane shear waves (when $d=2$) in an elastic solid with~mass density $\rho$ and shear modulus $G$, (ii) transverse electric (TE) or transverse magnetic (TM) waves (when $d=2$) in a dielectric medium endowed with permittivity~$\varepsilon=G^{-1}$ and permeability $\mu=\rho^{-1}$, and (iii) acoustic i.e. pressure waves (when~$d=2,3$) in a fluid characterized by the mass density~$\varrho = G^{-1}$ and bulk modulus~$\kappa=\rho^{-1}$. \\

At this point, we can deploy the results in~\cite{OK64,PBL78} to demonstrate (see~\ref{BWEappa} for details) that any $g\in L^2(\SSS)$ permits the Bloch wave expansion (BWE) as 
\begin{equation}\label{BWE}
g(\bx)= |\mathcal{B}|^{-1} \int_{\bk_s+\mathcal{B}}\tilde{g}_{\bk}(\bx)e^{ i \bk\cdot\bx} \dd\bk,
\end{equation}
where $\bk_s\in\mathbb{R}^d$ is an arbitrary shift vector, and 
\begin{equation}\label{BWE2}
\tilde{g}_{\bk}(\bx)= \sum_{\br\in\boldsymbol{R}}g(\bx+\br)e^{- i \bk\cdot(\bx+\br)} 
\end{equation}
belongs to $L^2_p(Y)$. This motivates us to consider a relatively broad class of source terms given by 
\begin{equation}\label{forceexp}
f(\bx)= |\mathcal{C}|^{-1} \int_{\bk_s+\mathcal{C}}\tilde{f}_{\bk}(\bx)e^{ i \bk\cdot\bx} \dd\bk,\quad \bk_s\in\overline{\mathcal{B}},
\end{equation}
where $\mathcal{C}\subset\mathcal{B}$ and $\tilde{f}_{\bk}\!\in L^2_p(Y)$ for $\bk\in\bk_s\!+\mathcal{C}$. Since~$\mathcal{C}\subset\mathcal{B}$, it is clear that~\eqref{forceexp} is nothing but a restriction of~\eqref{BWE} which implicitly defines a subset of~$L^2(\mathbb{R}^d)$. The main motivation behind~\eqref{forceexp} is to \emph{spectrally localize}~$f$, and thus~$u$, to a neighborhood of some~$\bk_s\!\in\overline{\mathcal{B}}$, which then greatly facilitates the asymptotic treatment. For future reference, we note that~\eqref{forceexp} covers the special cases where: (i) $\tilde{f}_\bk(\bx) = \rho(\bx) \hat{f}(\bk,\bx)$ with $\hat{f}(\bk,\cdot)\in L^2_p(Y)$, and (ii) $\tilde{f}_\bk(\bx)=1$ (with~$Y^{\tn}=Y^{\td}=\emptyset$) as in the related asymptotic treatments~\cite{GMO19,MG17} that rely on the plane wave expansion.

\begin{claim} \label{BPWESR}
Assuming $f\in L^2(\SSS)$, the Fourier integral representation 
\begin{equation}\label{PWE}
f(\bx)= \int_{\mathbb{R}^d}\mathfrak{F}(\bk)e^{ i \bk\cdot\bx} \dd\bk,
\end{equation}
is referred to as the plane wave expansion (PWE) of $f$. The relationship between BWE~\eqref{BWE} and PWE~\eqref{PWE} is given by
\begin{equation}\label{RPBWE}
|\mathcal{B}|^{-1}\tilde{f}_{\bk}(\bx) ~= \sum_{\br^*\in\boldsymbol{R}^*}\mathfrak{F} (\br^*\!+\bk)e^{ i \br^*\!\cdot\bx}.
\end{equation}
Further, when $\mathfrak{F} (\bk)$ is compactly supported within $\bk_s+\mathcal{B}$, relationship \eqref{RPBWE} simplifies to 
\begin{equation}\label{RPBWES}
|\mathcal{B}|^{-1}\tilde{f}_{\bk}(\bx)= \mathfrak{F} (\bk), 
\end{equation}
see~\ref{RB2PWEappa} for proof.
\end{claim}

\begin{remark}\label{commut}
For $\phi\in L^2_p(Y)$ and $f\in L^2(S)$, function $g(\bx)=\phi(\bx)f(\bx)\in L^2(S)$ and its BWE is given by \eqref{BWE}, where $\tilde{g}_{\bk}(\bx)=\phi(\bx) \tilde{f}_\bk(\bx)$ and $\tilde{f}_\bk$ is given by~\eqref{BWE2}.  
\end{remark}

By the linearity of \eqref{PDE}, we can account for~\eqref{forceexp} by focusing our analysis on the reduced field equation  
\begin{equation}\label{RPDES}
-\omega^2\rho(\bx)\hh u - \nabla\sip \big(G(\bx)\nabla u\big) ~=~ \tilde{f}_{\bk}(\bx)~e^{ i \bk\cdot\bx} \quad\text{in~~}\SSS.
\end{equation}
Thanks to the periodicity of~$\rho$ and~$G$ and the fact that $\tilde{f}_{\bk}\in L^2_p(Y)$, \eqref{RPDES} admits a Bloch wave solution $\tilde{u}(\bx)e^{ i \bk\cdot\bx}$, where $\tilde{u}:=\tilde{u}_\bk$ is $Y$-periodic and solves
\begin{eqnarray}\label{RPDE}
-\omega^2\rho(\bx)\hh \tilde{u} - \nabla_{\!\bk}\sip \big(G(\bx)\nabla_{\!\bk} \tilde{u}\big) ~=~ \tilde{f}_{\bk}(\bx) \quad\text{in~~}Y,
\end{eqnarray}
subject to boundary conditions 
\begin{eqnarray}\label{BCPDE}
\tilde{u}|_{\partial{Y^{\tpp}_{j0}}} &=& \tilde{u}|_{\partial{Y^{\tpp}_{j1}}}, \notag\\
\bnu\cdot G\nabla_{\!\bk}\tilde{u}|_{\partial{Y^{\tpp}_{j0}}} & =& -\bnu\cdot G\nabla_{\!\bk}\tilde{u}|_{\partial{Y^{\tpp}_{j1}}}, \notag \\%\quad j\in\overline{1,d},  \\
\bnu\cdot G\nabla_{\!\bk}\tilde{u}|_{\partial Y^{\tnn}} & = & 0,\notag \\
\tilde{u}|_{\partial Y^{\tdd}} & = & 0,
\end{eqnarray}
where~$\partial{Y^{\tpp}_{jm}}$ ($m=\overline{0,1}$) are given by~\eqref{perbcs}, $\nabla_{\!\bk} \!=\! \nabla +  i \bk$, and $\bnu$ is the unit outward normal on $\partial Y$. For brevity of notation, the dependence of~$\tilde{u}$ on~$\bk$ and~$\omega$ in \eqref{RPDE} and thereon is assumed implicitly.

%--------------------------------------------------------------------------------------------------------------%
\subsection{Eigenvalue problem} \label{eigen_problem}
% -------------------------------------------------------------------------------------------------------------%
\noindent For $\bk\in\mathbb{R}^d$ and $u,v \in H^1_{p0}(Y)$, we have
\begin{eqnarray}\label{varf}
\big(\nabla_{\!\bk}\sip (G(\bx)\nabla_{\!\bk} u),v\big) & = & |Y|^{-1} \int_{Y} \nabla_{\!\bk}\sip (G(\bx)\nabla_{\!\bk} u(\bx)) \, \overline{v(\bx)} \dd\bx \notag \\
& = &  -|Y|^{-1} \int_{Y} G(\bx)\nabla_{\!\bk} u(\bx) \cdot \overline{\nabla_{\!\bk} v(\bx)} \dd\bx  \;+\; |Y|^{-1}  \int_{\partial Y} \bnu\cdot (G(\bx)\nabla_{\!\bk} u(\bx)) \overline{v(\bx)} \dd\bx \notag \\
&=& -(G(\bx)\nabla_{\!\bk} u,\nabla_{\!\bk} v)   
\end{eqnarray}
thanks to the divergence theorem and boundary conditions~\eqref{BCPDE}. As a result, we find from the variational formulation that the operator $(-\rho^{-1}\nabla_{\!\bk}\sip(G\nabla_{\!\bk}))^{-1}$ from $L^2_{p}(Y)$ to~$L ^2_p(Y)$, subject to the germane boundary conditions, is self-adjoint and compact. Accordingly, \eqref{RPDE}--\eqref{BCPDE} are affiliated with the eigensystem $\{\tilde{\lambda}_n(\bk)
\!\in\!\mathbb{R}, \, \tilde{\phi}_n(\bk)\!\in\!H^1_{p0}(Y)\}$, that satisfies  
\begin{equation}\label{EVP}
-\tilde{\lambda}_n\rho(\bx)\hh \tilde{\phi}_n - \nabla_{\!\bk}\sip \big(G(\bx)\nabla_{\!\bk}~\tilde{\phi}_n\big) ~=~0 \quad\text{in~~}Y,
\end{equation}
subject to the boundary conditions  
\begin{eqnarray}\label{BCEVP}
\tilde{\phi}_n|_{\partial{Y^{\tpp}_{j0}}} &=& \tilde{\phi}_n|_{\partial{Y^{\tpp}_{j1}}}, \notag\\
\bnu\cdot G\nabla_{\!\bk}\tilde{\phi}_n|_{\partial{Y^{\tpp}_{j0}}} & =& -\bnu\cdot G\nabla_{\!\bk}\tilde{\phi}_n|_{\partial{Y^{\tpp}_{j1}}},  \notag \\%\quad j\in\overline{1,d},  \\
\bnu\cdot G\nabla_{\!\bk}\tilde{\phi}_n|_{\partial Y^{\tnn}} & = & 0, \notag \\
\tilde{\phi}_n|_{\partial Y^{\tdd}} & = & 0.
\end{eqnarray}
Note that the sequence $\{\tilde{\phi}_n\}$ is complete in $H^1_{p0}(Y)$ and $\rho$-orthogonal. We normalize the eigenfunctions so that $\|\tilde{\phi}_n\|=1$, whereby  
\begin{equation}\label{evpo}
(\rho \tilde{\phi}_n,\tilde{\phi}_m) ~=~ (\rho \tilde{\phi}_n,\tilde{\phi}_n)~\delta_{nm}. 
\end{equation}
For given~$\bk\!\in\!\mathbb{R}^d$, periodic medium~$\SSS$ thus permits the propagation of ``free'' Bloch waves $\tilde{\phi}_n(\bk) e^{i(\bk\cdot\bx-\omega_nt)}$ at eigenfrequency $\omega_n(\bk)=(\tilde{\lambda}_n)^{\frac{1}{2}}$. The set of all wavenumber-eigenfrequency pairs $(\bk,\omega_n)$ defines the \textit{Bloch dispersion relationship} of the medium. The latter is periodic in the reciprocal space, and is described completely by the first Brillouin zone $\mathcal{B}$ of the lattice. By the completeness of~$\tilde{\phi}_n$ in~$H_{p0}^1(Y)$, the solution $\tilde{u}$ of \eqref{RPDE}--\eqref{BCPDE} can be expanded as 
\begin{equation} \label{u_decomp}
\tilde{u}(\bx) = \sum_{n=1}^{\infty} \alpha_n \tilde{\phi}_n(\bx). 
\end{equation}
Provided that~$\omega^2\neq\tilde{\lambda}_n$ $\forall\hh n\in\mathbb{Z}^+$, \eqref{u_decomp} yields 
\begin{equation}\label{uts}
\tilde{u}(\bx) = -\sum_{n=1}^{\infty} \frac{(\tilde{f}_{\bk},\tilde{\phi}_n)~\tilde{\phi}_n(\bx)}{ (\rho \tilde{\phi}_n,\tilde{\phi}_n) (\omega^2-\tilde{\lambda}_n) } ,
\end{equation}
thanks to~\eqref{RPDE} and~\eqref{varf}. Then, by the linearity of~\eqref{PDE}, the total solution is expressed as 
 \begin{equation} \label{total_solution}
u(\bx) = \frac{1}{|\mathcal{C}|}\int_{\bk_s+\mathcal{C}} \Big( -\sum_{n=1}^{\infty} \frac{(\tilde{f}_{\bk},\tilde{\phi}_n)~\tilde{\phi}_n(\bx) e^{ i \bk\cdot\bx}}{ (\rho \tilde{\phi}_n,\tilde{\phi}_n) (\omega^2-\tilde{\lambda}_n) } \Big) \dd \bk.
\end{equation}
\noindent For future reference, we note that the weight of the $nth$ Bloch eigenmode in \eqref{total_solution} is inversely proportional to the spectral distance $|\omega^2-\tilde{\lambda}_n|$. 

%--------------------------------------------------------------------------------------------------------------%
\subsection{Scaling}\label{scaling}
%--------------------------------------------------------------------------------------------------------------%

\noindent In what follows, we seek a homogenized description of~\eqref{RPDE}--\eqref{BCPDE} in a spectral neighborhood of the wavenumber-frequency pair
\[
\big(\bk_s,\omega_n(\bk_s)=(\tilde{\lambda}_n(\bk_s))^{1/2}\big)\in\overline{\mathcal{B}}\times\mathbb{R}, \qquad n\in\mathbb{Z}^+,
\]
and we assume all quantities to be \emph{a priori} normalized by some reference ``mass density'' $\rho_0$, ``shear modulus'' $G_0$ and lengthscale $\ell_0$. On making use of the short-hand notation $\tilde{\lambda}_n=\tilde{\lambda}_n(\bk_s)$ and~$\omega_n=\omega_n(\bk_s)$ hereon, we next introduce the perturbation parameter $\eps=o(1)$ defining the spectral neighborhood as 
\begin{eqnarray}
\bk=\bk_s+\eps\hat{\bk},\qquad \omega^2 = \tilde{\lambda}_n + \eps\hh\chsig \check{\omega}^2 + \eps^2\hasig \hat{\omega}^2, \qquad \chsig ,\hasig \in\{-1,0,1\}, \quad \chsig \hh \hasig  = 0, \quad |\chsig + \hasig | = 1. \label{scal} 
\end{eqnarray}

\begin{remark}\label{scaling:choice}
Through the design of~$\check\sigma$ and~$\hat\sigma$, frequency separation parameters $\check\omega$ and~$\hat\omega$ are meant to be used in the ``either or'' sense, depending on the driving frequency (when~$\tilde{f}_\bk\!\neq\! 0$) and the local geometry of germane dispersion surface (when~$\tilde{f}_\bk\!=\!0$).  Specifically when~$\tilde{f}_\bk\neq 0$ whereby~$\omega$ is given, we have 
\begin{equation}  \label{komega}
\omega^2-\tilde{\lambda}_n \;=\; \left\{\!\!\begin{array}{lcrr}
O(\eps) & \Rightarrow &~  |\check\sigma|=1, & \hat\sigma~=0 \\
O(\eps^2)   & \Rightarrow &~ \check\sigma~=0, & |\hat\sigma|=1 \end{array}\right..
\end{equation}
When~$\tilde{f}_\bk\!=\!0$, on the other hand, it will be for instance shown that for dispersion surfaces with locally parabolic (resp.~conical) sections, the frequency in those $\bk$-directions scales as $\omega_n^2(\bk)\!-\!\tilde{\lambda}_n \!= \eps^2\hat\omega^2$ (resp. $\omega_n^2(\bk)\!-\!\tilde{\lambda}_n \!= \eps\hh\breve\omega^2$). Since such information is not available beforehand, the idea is to substitute~\eqref{scal} ``as is'' into the field equation~\eqref{RPDE}, and then to identify the appropriate frequency scaling (by letting either~$\check\sigma=0$ or~$\hat\sigma=0$) depending on the local eigenfunction structure. In order to bring the analyses of both forced ($\tilde{f}_\bk\!\neq\!0$) and free ($\tilde{f}_\bk\!=\!0$) wave motion under a common umbrella, we will uniformly start from the agnostic scaling law~\eqref{scal} throughout the remainder of this work. 
\end{remark}

In the context of~\eqref{scal}, we are now in position to pursue the ansatz %~of \eqref{RPDE}
\begin{equation}\label{sol_scal}
\tilde{u}(\bx) \;=\; \eps^{-2}  \sum_{m=0}^\infty \, \eps^m \hh \tilde{u}_m(\bx),
\end{equation}
via the asymptotic expansion of~\eqref{RPDE}--\eqref{BCPDE}, see also~\cite{MG17,GMO19}. For completenes, we note that the presence of the factor~$\eps^{-2}$ in front of the series is motivated by~\eqref{uts} and the smallness of $|\omega^2-\tilde\lambda_n(\bk)|$ suggested by~\eqref{scal}. On inserting \eqref{scal}--\eqref{sol_scal} into \eqref{RPDE}--\eqref{BCPDE} and letting~$\tilde{f}_{\bk}=O(1)$, we obtain a cascade of field equations over $Y$, namely 
\begin{eqnarray}
O(\eps^{-2}): &\!\!\!\! -\tilde{\lambda}_n\rho\hh \tilde{u}_0 - \nabla_{\!\bk_s}\!\sip \big(G\hh\nabla_{\!\bk_s} \tilde{u}_0\big) = 0, \label{eps-2}\\
O(\eps^{-1}): &\!\!\!\! -\tilde{\lambda}_n\rho\hh \tilde{u}_1 - \nabla_{\!\bk_s}\!\sip \big(G(\nabla_{\!\bk_s} \tilde{u}_1+\tilde{u}_0\hh i \hat{\bk})\big) - G(\nabla_{\!\bk_s} \tilde{u}_0)\sip  i \hat{\bk} - \check\sigma\rho \hh \check{\omega}^2\tilde{u}_0 = 0,\label{eps-1}\\
O(1): &\!\!\!\! -\tilde{\lambda}_n\rho\hh \tilde{u}_{2} - \nabla_{\!\bk_s}\!\sip \big(G(\nabla_{\!\bk_s} \tilde{u}_{2}+\tilde{u}_{1}\hh i \hat{\bk})\big) - G(\nabla_{\!\bk_s} \tilde{u}_{1}+\tilde{u}_{0}\hh i \hat{\bk})\sip  i \hat{\bk} -\chsig \rho \hh \check{\omega}^2\tilde{u}_{1}-\hasig \rho \hh \hat{\omega}^2\tilde{u}_{0} = \tilde{f}_{\bk},  \label{eps0} \\
O(\eps^{m\geqslant 1}): &\!\!\!\! -\tilde{\lambda}_n\rho\hh \tilde{u}_{m+2} - \nabla_{\!\bk_s}\!\sip \big(G(\nabla_{\!\bk_s} \tilde{u}_{m+2}+\tilde{u}_{m+1}\hh i \hat{\bk})\big) - G(\nabla_{\!\bk_s} \tilde{u}_{m+1}+\tilde{u}_{m}\hh i \hat{\bk})\sip  i \hat{\bk} \label{epsp} \notag  
-\chsig \rho \hh \check{\omega}^2\tilde{u}_{m+1}-\hasig \rho \hh \hat{\omega}^2\tilde{u}_{m} = 0, \\  & 
\end{eqnarray}
\noindent along with the sequence of boundary conditions 
\begin{eqnarray}\label{BCPDES}
\tilde{u}_m|_{\partial{Y^{\tpp}_{j0}}} \!\!&=&\!\! \tilde{u}_m|_{\partial{Y^{\tpp}_{j1}}}, \notag\\
\bnu\sip G(\nabla_{\!\bk_{s}}\tilde{u}_m+\tilde{u}_{m-1}\hh i \hat{\bk})|_{\partial{Y^{\tpp}_{j0}}} \!\!&=&\!\! -\bnu\sip G(\nabla_{\!\bk_{s}}\tilde{u}_m+\tilde{u}_{m-1}\hh i \hat{\bk})|_{\partial{Y^{\tpp}_{j1}}} \notag \\
\bnu\sip G(\nabla_{\!\bk_{s}}\tilde{u}_m+\tilde{u}_{m-1}\hh i \hat{\bk})|_{\partial Y^{\tnn}} \!\!&=&\!\! 0,\notag \\
\tilde{u}_m|_{\partial Y^{\tdd}} \!\!&=&\!\! 0,\quad  m\geqslant 0,
\end{eqnarray}
\noindent where $\tilde{u}_{-1}\equiv 0$. In the sequel, we say that tensor $\boldsymbol{g}\!\in\! (\bar{H}^{1}_{p0}(Y) )^{d^q}$, $q \!\geqslant\! 1$ satisfies the \textit{``flux boundary conditions''} if 
\begin{eqnarray}\label{gbcs}
\bnu\sip \boldsymbol{g}|_{\partial{Y^{\tpp}_{j0}}} & =& -\bnu\sip \boldsymbol{g}|_{\partial{Y^{\tpp}_{j1}}},\notag \\ %\quad j\in\overline{1,d}, 
\bnu\sip \boldsymbol{g}|_{\partial Y^{\tnn}} & = & \boldsymbol{0}.
\end{eqnarray}

%--------------------------------------------------------------------------------------------------------------%
\subsection{Averaging operators and effective solution}\label{effective solution}
%--------------------------------------------------------------------------------------------------------------% 
\noindent Let $n_{\nes q}\!\in\mathbb{Z}^+$ ($q=\overline{1,Q}$) collect the ``nearby'' dispersion branches $\omega_{n_{\nes q}}(\bk)$ traversing the vicinity of~$(\bk_s,\tilde{\lambda}_n^{1/2})$, where we aim to pursue ansatz~\eqref{sol_scal}. With such setup in mind, we introduce the \textit{averaging} operators $\langle\cdot\rangle^{n_{\nes q}}$ and $\langle\cdot\rangle^{n_{\nes q}}_\rho$ and the ``zero mean'' Sobolev space $\bar{H}^{1}_{p0}(Y) $ as
\begin{eqnarray}\label{aveop}
\langle \tilde{g} \rangle^{n_{\nes q}} &=& (\tilde{g},\tilde{\phi}_{n_{\nes q}}), \label{aveop1} \\
\langle \tilde{g} \rangle^{n_{\nes q}}_\rho &=& (\rho \tilde{\phi}_{n_{\nes q}},\tilde{\phi}_{n_{\nes q}})^{-1}(\rho \tilde{g},\tilde{\phi}_{n_{\nes q}}), \label{aveop2} \\
\bar{H}^{1}_{p0}(Y) &=& \{\tilde{g}\in H^{1}_{p}(Y): \langle \tilde{g} \rangle^{n_{\nes q}}_\rho=0, ~q=\overline{1,Q}\}. \label{aveop3}
\end{eqnarray}
For completeness, we note that our definition~\eqref{aveop3} of the ``zero mean'' Sobolev space~$\bar{H}^{1}_{p0}(Y)$ is different from that in~\cite{GMO19} which postulates~$\langle \tilde{g} \rangle^{n_{\nes q}}=0$ in lieu of~$\langle \tilde{g} \rangle^{n_{\nes q}}_\rho=0$, and from that in \cite{CKP10} where the functions~$\{\tilde{g},\tilde{\varphi}_{\nes q} (q=\overline{1,Q})\}$ are assumed to be linearly independent. 
For~$\tilde{g}=\tilde{u}$, we will use the short-hand notation
\begin{equation}
u_{mq}(\eps\hat{\bk}) \;:=\; \langle \tilde{u}_m \rangle^{n_{\nes q}}_\rho, \qquad q=\overline{1,Q}. \label{umq}
\end{equation}
On the basis of~\eqref{sol_scal} and~\eqref{umq}, we can adapt the definition of \textit{effective solution}~\cite{GMO19} at wavenumber $\bk_s\!+\eps\hat{\bk}$ as 
\begin{equation}\label{effsolud}
\langle \tilde{u} \rangle^{n_{\nes q}}_{\!\rho} (\eps\hat{\bk}) \;=\; \sum _{m=0}^{\infty} \eps^{m-2}\hh u_{mq}, \qquad q=\overline{1,Q} 
\end{equation}
which then provides the basis for computing the (set of) effective solution(s) near $\bk_s$ in the physical space as
\begin{eqnarray}\label{effsolu}
\langle  u \rangle^{n_{\nes q}}_{\!\rho} (\bx) \;=\; |\mathcal{C}|^{-1} \int_{\mathcal{C}} \langle \tilde{u} \rangle^{n_{\nes q}}_{\!\rho}(\eps\hat{\bk})\, e^{ i (\bk_s\!+\eps\hat{\bk})\cdot\bx} \dd (\eps\hat{\bk}), \quad \bx\in\mathbb{R}^d.
\end{eqnarray}

\begin{remark}
In situations where $Q\!=\!1$ and~$n_1\!=\!n$ which corresponds to the case of an isolated branch, $\langle \cdot \rangle^{n}$, $\langle \cdot \rangle^{n}_\rho$, and~$u_{m1}$ will be conveniently denoted as $\langle\cdot\rangle$, $\langle\cdot\rangle_{\!\rho}$ and~$u_m$, respectively. In this case, \eqref{effsolud} and~\eqref{effsolu} reduce to 
\begin{equation}\label{effsoluds}
\langle\tilde{u}\rangle_{\!\rho}(\eps\hat{\bk}) \;=\; \sum_{m=0}^\infty\eps^{m-2} \hh u_m, \qquad
\langle  u \rangle_{\!\rho} (\bx) \;=\; |\mathcal{C}|^{-1} \int_{\mathcal{C}} \langle \tilde{u} \rangle_{\!\rho}(\eps\hat{\bk})\, e^{ i (\bk_s\!+\eps\hat{\bk})\cdot\bx} \dd (\eps\hat{\bk}).  
\end{equation}
\end{remark}

For future reference, we also define the symmetrization operators $\{\cdot\}$ and~$\{\cdot\}'$ on tensors $\btau\in\mathbb{C}^{d^n}$, $n\geqslant 2$ as
\begin{eqnarray}
\{\btau\}_{p_1,p_2,\dots,p_n} &=& \frac{1}{n!} \sum_{(q_1,q_2,\dots,q_n)\hh \in\hh \Pi_n} \btau_{q_1,q_2,\dots,q_n}, \label{symop1} \\
\{\btau\}_{p_1,p_2,\dots,p_n}' &=& \frac{1}{(n-1)!} \sum_{(q_2,\dots,q_n)\hh \in\hh \Pi_{n-1}} \btau_{p_1,q_2,\dots,q_n}, \quad p_1,p_2,\dots,p_n \in \overline{1,d}, \label{symop2}
\end{eqnarray}
respectively, where $\Pi_n$  is the set of all permutations of set $\overline{1,n}$.

\begin{remark}
In the context of~\eqref{scal}, we first recall the multicell homogenization i.e. ``folding'' technique \cite{ACG13b,Doss2008,GMO19,Gone2010} which enables evaluation of the effective properties of a periodic medium at rational wavenumbers $\bk_s\!=\! \sum_j q_j\be^j\in\overline{\mathcal{B}}$, $q_j\!\in\!\mathbb{Q}$. A leading-order expansion about an arbitrary (rational or irrational) wavenumber~$\bk_s\!\in\overline{\mathcal{B}}$ was implicitly considered in \cite{MAMGC16}, via multiple scales approach, near simple and repeated eigenfrequencies. With reference to~\eqref{umq}, we specifically pursue explicit (first- or second-order) effective descriptions governing~$u_{mq}$ ($m=\overline{0,2}$, $q=\overline{1,Q}$) at an arbitrary wavenumber $\bk_s\!\in\overline{\mathcal{B}}$ near simple, multiple, and nearby eigenfrequencies.
\end{remark}

\begin{remark}\label{warning!}
When $f_\bk\!=\!0$ identically, the applicability of any effective model for given perturbation vector~$\hat\bk$ also implies its validity for~$\alpha\hat\bk$, $\alpha \!\leqslant\! O(1)$ thanks to the arbitrariness of~$\eps=o(1)$ in~\eqref{eps-2}--\eqref{epsp}. When~$f_\bk\neq 0$, on the other hand, this implication holds as long as the point~$(\bk_s\!+\!\alpha\eps\hat\bk,\omega)$ does not lie on the germane dispersion branch, i.e.~as long as $\omega_{n_q}(\bk_s\!+\alpha\eps\hat\bk)\neq \omega$. To provide a focus for the analysis, we hereon (i) identify the wavenumber perturbations by their direction~$\hat\bk/\|\hat\bk\|$, and (ii) for ~$f_\bk\neq 0$ we restrict our consideration to~$\hat\bk\in\mathcal{K}_\eps$, where
\begin{equation}\label{Kset}
\mathcal{K}_\eps = \{\hat\bk\in\mathbb{R}^d: \eps\hat\bk\in\mathcal{C}, ~\omega_{n_q}^{\tmm}(\bk_s+\eps\hat\bk) \neq \omega\}, \qquad q=\overline{1,Q}, 
\end{equation}
where~$\omega_{n_q}^{\tmm}(\bk)$ denotes the~$m$th order approximation of~$\omega_{n_q}(\bk)$ affiliated with~$u_{mq}$ in~\eqref{umq}.
\end{remark}

% --------------------------------------------------------------------------------------------------------------%
\section{Simple eigenvalue} \label{simple_eig}
%--------------------------------------------------------------------------------------------------------------%
%--------------------------------------------------------------------------------------------------------------%
\subsection{Leading-order approximation}\label{leador}
% -----------------------------------------------------------------------------------------------------------------------------------------------------------------------------------%

\noindent With reference to the eigenvalue problem \eqref{EVP}-\eqref{BCEVP}, the solution of \eqref{eps-2} in the vicinity of a simple eigenfrequency $\omega_n$ is expressed as 
\begin{equation}\label{ut0exp}
\tilde{u}_0(\bx) \;=\;  u_{0} \hh \tilde{\phi}_n(\bx), \quad u_{0} \in\mathbb{C},
\end{equation}
where~$u_0=\langle\tilde{u}_0\rangle_{\!\rho}$ as stated before. Then, by inserting~\eqref{ut0exp} into~\eqref{eps-1} and integrating~$(\eqref{eps-1}, \tilde{\phi}_n)$ by parts via the boundary conditions~\eqref{BCPDES} with $m=1$, we obtain the averaged $O(\eps^{-1})$ statement as
\begin{equation}\label{solvacon-1}
-(\bth^{\tz}\sip (i\hat{\bk}) + \chsig \check{\omega}^2 \rho^{\tz}) u_{0}~= 0,
\end{equation}
where 
\begin{equation}\label{teta&rho}
\bth^{\tz} =  \langle G\nabla_{\!\bk_{s}} \tilde{\phi}_n \rangle -  \overline{ \langle G\nabla_{\!\bk_{s}}\tilde{\phi}_n\rangle} \in  i \mathbb{R}^d \quad~ \text{and} ~\quad \rho^{\tz} = \langle \rho\tilde{\phi}_n\rangle \in \mathbb{R}^{+}.
\end{equation}
On substituting~\eqref{ut0exp} in~\eqref{eps-1}, one finds by the linearity of the problem that
\begin{equation}\label{ut1exp}
\tilde{u}_1(\bx) = u_0~\bchi^{\so}(\bx)\sip (i\hat{\bk})+ u_1 \tilde{\phi}_n(\bx), \quad u_1\in\mathbb{C},
\end{equation}
where~$u_1=\langle\tilde{u}_1\rangle_{\!\rho}$ and $\bchi^{\so}\in(\bar{H}^{1}_{p0}(Y))^d$ uniquely solves the unit cell problem
\begin{equation}\label{chi1}
\tilde{\lambda}_n\hh\rho\hh\bchi^{\so} + \nabla_{\!\bk_s}\!\sip \big(G(\nabla_{\!\bk_s}\hh\bchi^{\so}\!+\tilde{\phi}_n \boldsymbol{I})\big) + G\nabla_{\!\bk_s} \tilde{\phi}_n - \frac{\rho}{\rho^{\tz}}\tilde{\phi}_n\bth^{\tz}~=~\boldsymbol{0}, 
\end{equation} 
subject to the boundary conditions~\eqref{gbcs} with $\bg\!=\!G(\nabla_{\!\bk_s}\hh\bchi^{\so}\!+\tilde{\phi}_n \boldsymbol{I})$ and $\boldsymbol{I}$ denoting the second-order identity tensor.

We next consider the $O(1)$ field equation \eqref{eps0}. On recalling~\eqref{ut0exp} and~\eqref{ut1exp}, we can integrate~$(\eqref{eps0}, \tilde{\phi}_n)$ by parts aided by the boundary conditions~\eqref{BCPDES} with $m=2$ to obtain the averaged $O(1)$ statement
\begin{equation}\label{solvacon0}
-\big( \bmu^{\tz} : (i\hat{\bk})^2 + \hasig \hat{\omega}^2 \rho^{\tz}\big )\hh u_0-(\bth^{\tz}\sip (i\hat{\bk})+\chsig \check{\omega}^2\rho^{\tz})\hh u_1 \;=\; \langle \tilde{f}_{\bk} \rangle, 
\end{equation}
where
\[
(i\hat\bk)^m  ~=~ (i\hat\bk) \otimes (i\hat\bk) \ldots \otimes (i\hat\bk) \quad (m \text{ times}),
\]
and
\begin{equation}\label{mu0}
\bmu^{\tz} = \langle G\{\nabla_{\!\bk_s}~\bchi^{\so}\!+ \tilde{\phi}_n \boldsymbol{I}\} \rangle - \big \{\big( G \bchi^{\so}\otimes \overline{\nabla_{\!\bk_s}\tilde{\phi}_n},1 \big )\big \}.
\end{equation}

\begin{claim}\label{mu0real}
For any~$\bk_s\!\in\overline{\mathcal{B}}$, effective tensor~$\bmu^{\tz}$ is real-valued, i.e. $\bmu^{\tz}\in\mathbb{R}^{d\times d}$. See~\ref{Effcoef} for proof.
%It is further definite positive at $\bk_s=0$ for $n=1$ if $\omega_1=0$.
\end{claim}

\begin{claim} \label{specase1}
For wavenumbers $\bk_s = \frac{1}{2}(\sum_j n_j\hh\be^j)$, $n_j\in\{-1,0,1\}$, which include the origin and apexes of the first Brillouin zone $\mathcal{B}$, Bloch eigenfunction $\tilde{\phi}_n(\bx) e^{ i \bk_s\cdot\bx}$ is real-valued up to a constant multiplier $e^{i\varphi_0}$. As a result, in such cases we find that $\bth^{\tz}=\boldsymbol{0}$. See~\ref{specases} for proof.
\end{claim}

Claim~\ref{specase1} motivates us to consider separately the situations when $\bth^{\tz}\!\neq\boldsymbol{0}$ and $\bth^{\tz}\!=\boldsymbol{0}$, which we address next.

%--------------------------------------------------------------------------------------------------------------%
\subsubsection{Effective model for non-trivial $\bth^{\tz}$}\label{simpleadcas1}
% -----------------------------------------------------------------------------------------------------------------------------------------------------------------------------------%

\noindent As can be seen from the foregoing analysis, presence of the source term in the~$O(1)$ statement~\eqref{solvacon0} requires that its~$O(\eps^{-1})$ predecessor \eqref{solvacon-1} be satisfied \emph{identically}. When~$\tilde{f}_{\bk}\neq 0$ and $\omega^2\!-\omega_n^2=O(\eps)$ whereby~$|\check\sigma|=1$ due to~\eqref{komega}, we must have $u_0=0$ in~\eqref{solvacon-1} thanks to Remark~\ref{warning!} which guarantees that the multiplier $(\bth^{\tz}\sip (i\hat{\bk})+\chsig \rho^{\tz}\check{\omega}^2)$ is non-trivial. As a result, \eqref{solvacon0} yields the leading-order effective equation 
\begin{equation}\label{leadorder11}
-(\bth^{\tz}\sip (i\hat{\bk})+\chsig \check{\omega}^2 \rho^{\tz})\hh u_1 ~=~ \langle \tilde{f}_{\bk}\rangle.
\end{equation}
A similar treatment can be pursued for the situation when $\omega^2\!-\omega_n^2=O(\eps^2)$, in which case~$|\hat\sigma|=1$. This case is not addressed for the reasons of brevity.

In the absence of the source term $\tilde{f}_{\bk}$, on the other hand, the existence of a non-trivial wavefield solving~\eqref{solvacon-1} and~\eqref{solvacon0} independently requires that $|\check\sigma|=1$. In this case, \eqref{leadorder11} with~$\langle \tilde{f}_{\bk}\rangle=0$ furnishes the leading-order asymptotic approximation of the dispersion relationship and group velocity near $(\bk_s,\omega_n\!>\!0)$ as 
\begin{equation}\label{drasymp11}
\omega_n^2(\bk) \;=\; \omega_n^2 - \frac{1}{\rho^{\tz}}\hh i \bth^{\tz}\sip(\eps \hat{\bk}), \qquad
\boldsymbol{c}_g = \frac{\dd\omega_n(\bk)}{\dd\bk} \;=\; \frac{-1}{2\omega_n\hh\rho^{\tz}}\hh i \bth^{\tz} 
\end{equation}
respectively, where $\omega_n$ (without an argument) refers to~$\omega_n(\bk_s)$ as stated earlier. Geometrically, \eqref{drasymp11} describes the $n$th dispersion (hyper-) surface locally as a (hyper-) plane, where~$\boldsymbol{c}_g$ signifies its ``steepest slope''.

%--------------------------------------------------------------------------------------------------------------%
\subsubsection{Effective model for trivial $\bth^{\tz}$ }\label{simpleadcas2}
% -----------------------------------------------------------------------------------------------------------------------------------------------------------------------------------%

\noindent When $\tilde{f}_\bk\neq 0$ and $\omega^2\!-\omega_n^2=O(\eps^2)$, we have that~$|\hat\sigma|=1$ thanks to~\eqref{komega}. In this case the $O(\eps^{-1})$ statement~\eqref{solvacon-1} is satisfied identically, while its~$O(1)$ companion~\eqref{solvacon0} produces the effective equation 
\begin{equation}\label{leadorder12}
-\big( \bmu^{\tz} \!:\! (i\hat{\bk})^2 + \hasig \rho^{\tz} \hat{\omega}^2 \big)\hh u_0~=~\langle \tilde{f}_{\bk} \rangle,
\end{equation}
With reference to Claim \ref{specase1}, \eqref{leadorder12} in particular describes the response of a periodic medium near the \emph{origin} and \emph{apexes} of the first Brillouin zone. The nature of such response depends on (i) the sign definiteness of $\bmu^{\tz}$, and (ii) the sign of $\omega^2 -\omega_n^2$. For example, when $\bmu^{\tz}$ is sign-definite oppositely to the sign of $\omega^2 -\omega_n^2$, the effective medium is ``dissipative'' in that~$\omega$ resides inside a band gap~\cite{GMO19} terminating at~$\omega_n$.

When $\tilde{f}_\bk=0$, on the other hand, from~\eqref{solvacon-1} and~\eqref{solvacon0} we find that a non-trivial solution is possible only if~$\check\sigma=0$, i.e.~$|\hat\sigma|=1$. In this case~\eqref{solvacon-1} is again satisfied identically, while \eqref{solvacon0} provides the leading-order approximation of dispersion relationship and group velocity near $(\bk_s,\omega_n\!>\!0)$ as 
\begin{equation}\label{drasymp12}
\omega_n^2(\bk) \;=\; \omega_n^2 + \frac{1}{\rho^{\tz}}~\bmu^{\tz}\!:\!(\eps \hat{\bk})^2, \qquad \boldsymbol{c}_g(\bk)\;=\;  \frac{1}{\omega_n\hh\rho^{\tz}}~\bmu^{\tz}\cdot(\eps \hat{\bk}).
\end{equation}

%--------------------------------------------------------------------------------------------------------------%
\subsection{Second-order correctors}\label{secondor1}
% -----------------------------------------------------------------------------------------------------------------------------------------------------------------------------------%

\noindent With~\eqref{ut0exp} and~\eqref{ut1exp} at hand, one can make use of the averaged statements~\eqref{solvacon-1} and~\eqref{solvacon0} to solve the~$O(1)$ field equation~\eqref{eps0} in terms of~$\tilde{u}_2$ as 
\begin{equation}\label{ut2exp}
\tilde{u}_2(\bx) = u_0\hh\bchi^{\st}(\bx):(i\hat{\bk})^2+ u_1\hh\bchi^{\so}(\bx)\cdot (i\hat{\bk}) +u_2\hh\tilde{\phi}_n (\bx)+ \eta^{\tz}(\bx),
\end{equation}
\noindent where $\bchi^{\st}\!\in\big(\bar{H}^{1}_{p0}(Y)\big)^{d\times d}$ and $\eta^{\tz}\!\in\bar{H}^{1}_{p0}(Y)$ uniquely solve the respective cell problems  
\begin{equation}\label{chi2}
\tilde{\lambda}_n\hh\rho\hh\bchi^{\st} + \nabla_{\!\bk_s}\!\sip \big(G(\nabla_{\!\bk_s}\hh\bchi^{\st}+\{\boldsymbol{I}\otimes\bchi^{\so}\}')\big) + G(\nabla_{\!\bk_s}\hh\bchi^{\so}\!+\tilde{\phi}_n \boldsymbol{I})~=~\frac{\rho}{\rho^{\tz}}\{\bth^{\tz}\otimes\bchi^{\so}\} + \frac{\rho}{\rho^{\tz}}\tilde{\phi}_n\hh\bmu^{\tz},
\end{equation}
\begin{equation}\label{eta0}
-\tilde{\lambda}_n\rho\hh\eta^{\tz} - \nabla_{\!\bk_s}\!\sip \big(G\nabla_{\!\bk_s} \eta^{\tz}\big)~=~\tilde{f}_{\bk}-\frac{\rho}{\rho^{\tz}}\langle \tilde{f}_{\bk}\rangle\hh\tilde{\phi}_n,
\end{equation}
\noindent with $G(\nabla_{\!\bk_s}\hh\bchi^{\st}+ \{\boldsymbol{I}\otimes\bchi^{\so}\}')$ and $G\nabla_{\!\bk_s}\eta^{\tz}$ each satisfying the flux boundary conditions \eqref{gbcs}.

\begin{remark}\label{eta0_lin}
Cell function $\eta^{\tz}$ depends implicitly on $\eps\hat{\bk}$ via $\tilde{f}_{\bk}=\tilde{f}_{\bk}(\bx)$. In situations when $\tilde{f}_{\bk}(\bx) = F(\eps\hat{\bk})\hh\phi(\bx)$, the solution of~\eqref{eta0} is given by $\eta^{\tz}(\bx)=F(\eps\hat{\bk})\hh\zeta^{\tz}(\bx)$, where $\zeta^{\tz}$ uniquely solves 
\begin{equation}\label{zeta0}
-\tilde{\lambda}_n\rho\hh\zeta^{\tz} - \nabla_{\!\bk_s}\!\sip \big(G\nabla_{\!\bk_s} \zeta^{\tz}\big)~=~\phi-\frac{\rho}{\rho^{\tz}}\langle\phi\rangle\hh\tilde{\phi}_n,
\end{equation}
\noindent with $G\nabla_{\!\bk_s}\zeta^{\tz}$ satisfying the flux boundary conditions \eqref{gbcs}.
\end{remark}

\begin{claim}\label{eta_chi1_id}
The following identity holds: 
\begin{equation}\label{eta0_chi1_id}
(G\eta^{\tz},\nabla_{\!\bk_s}\hh\tilde{\phi}_n)-\langle G \nabla_{\!\bk_s}\hh\eta^{\tz} \rangle \;=\; (\tilde{f}_{\bk},\bchi^{\so}). 
\end{equation}
\noindent See \ref{identities} for proof.
\end{claim}
We next consider the $O(\eps)$ field equation \eqref{epsp} with $m\!=\!1$. On substituting \eqref{ut0exp}, \eqref{ut1exp} and \eqref{ut2exp} into~\eqref{epsp}, integrating~$(\eqref{epsp}, \tilde{\phi}_n)$ by parts via the boundary conditions~\eqref{BCPDES} with $m=3$, and exploiting Claim~\ref{eta_chi1_id}, we obtain the averaged $O(\eps)$ statement 
\begin{eqnarray}\label{solvacon1}
-\bth^{\so}:(i\hat{\bk})^3 u_0 -\big( \bmu^{\tz} : (i\hat{\bk})^2 + \hasig \rho^{\tz} \hat{\omega}^2 \big )~u_1-(\bth^{\tz}\sip (i\hat{\bk})+\chsig \check{\omega}^2\rho^{\tz})~u_2~=~M_1(\hat{\bk}),
\end{eqnarray}
where
\begin{eqnarray}
\bth^{\so} &\!\!\!=\!\!\!& \langle G\{\nabla_{\!\bk_s}\bchi^{\st}+ \boldsymbol{I}\otimes\bchi^{\so}\} \rangle- \big \{\big( G \bchi^{\st}\otimes\overline{\nabla_{\!\bk_s}\tilde{\phi}_n},1 \big )\big \}, \\
M_1(\hat{\bk}) &\!\!\!=\!\!\!&  - (\tilde{f}_{\bk},\bchi^{\so})\sip (i\hat{\bk}).  \label{force1}
\end{eqnarray}

\begin{claim}\label{theta1r}
For any~$\bk_s\!\in\overline{\mathcal{B}}$, effective tensor $\bth^{\so}$ is imaginary-valued, namely $\bth^{\so}\!\in i\mathbb{R}^{d\times d \times d}$. See \ref{Effcoef} for proof.
\end{claim}

Proceeding with the analysis, we make use of the solutions~\eqref{ut0exp}, \eqref{ut1exp} and~\eqref{ut2exp} in conjunction with ~averaged statements \eqref{solvacon-1}, \eqref{solvacon0} and~\eqref{solvacon1} to solve the~$O(\eps)$ field equation~\eqref{epsp} with $m=1$ in terms of~$\tilde{u}_3$  as
\begin{equation}\label{ut3exp}
\tilde{u}_3(\bx) = u_0\hh\bchi^{\sh}(\bx):(i\hat{\bk})^3+ u_1\hh\bchi^{\st}(\bx) : (i\hat{\bk})^2 +u_2\hh\bchi^{\so}(\bx)\cdot (i\hat{\bk}) + u_3\hh\tilde{\phi}_n(\bx) + \boldsymbol{\eta}^{\so}(\bx)\sip (i\hat{\bk})+\chsig \check{\omega}^2\eta^{\st}(\bx),
\end{equation}
where $\bchi^{\sh}\!\in\big(\bar{H}^{1}_{p0}(Y)\big)^{d \times d \times d}$, $\boldsymbol{\eta}^{\so}\!\in\big(\bar{H}^{1}_{p0}(Y)\big)^{d}$ and~$\eta^{\st}\!\in\bar{H}^{1}_{p0}(Y)$ uniquely solve the respective cell problems 
\begin{eqnarray}
&\hspace*{-50mm}\tilde{\lambda}_n\hh\rho\hh\bchi^{\sh} + \nabla_{\!\bk_s}\!\sip \big(G(\nabla_{\!\bk_s}\hh\bchi^{\sh}+\{\boldsymbol{I}\otimes\bchi^{\st}\}' )\big) + G\{ \nabla_{\!\bk_s}\hh\bchi^{\st}+\boldsymbol{I}\otimes\bchi^{\so}\} &\notag \\
&\hspace*{50mm} = \frac{\rho}{\rho^{\tz}} \big \{\bth^{\tz}\otimes\bchi^{\st} \big \} + \frac{\rho}{\rho^{\tz}}\hh\big \{\bmu^{\tz}\otimes \bchi^{\so} \big \} + \frac{\rho}{\rho^{\tz}}\tilde{\phi}_n\hh\bth^{\so}, & \label{chi3} \\
 \notag \\
&\hspace*{0mm}\tilde{\lambda}_n\rho\hh\boldsymbol{\eta}^{\so} + \nabla_{\!\bk_s}\!\sip \big(G(\nabla_{\!\bk_s} \boldsymbol{\eta}^{\so}+\eta^{\tz}\boldsymbol{I} )\big) +G\nabla_{\!\bk_s}\eta^{\tz} +\frac{\rho}{\rho^{\tz}}\tilde{\phi}_n(\tilde{f}_{\bk},\bchi^{\so})~=~ \frac{\rho}{\rho^{\tz}}\langle\tilde{f}_{\bk}\rangle\bchi^{\so},& \label{eta_1}\\
\notag  \\
&-\tilde{\lambda}_n \rho \eta^{\st} - \nabla_{\!\bk_s}\sip(G\nabla_{\!\bk_s}\hh\eta^{\st})~=~\rho\eta^{\tz},& \label{eta_2}
\end{eqnarray}
with $G(\nabla_{\!\bk_s}\hh\bchi^{\sh}+\tilde{\phi}_n~\{\boldsymbol{I}\otimes\bchi^{\st}\}')$, $G(\nabla_{\!\bk_s} \boldsymbol{\eta}^{\so}+\eta^{\tz}\boldsymbol{I} )$ and~$G\nabla_{\!\bk_s}\eta^{\st}$ each satisfying the flux boundary conditions~\eqref{gbcs}. On recalling~\eqref{eta0_chi1_id}, one may note that both~$\boldsymbol{\eta}^{\so}$ and~$\eta^{\st}$ are $\tilde{f}_\bk$-dependent.

\begin{claim}\label{eta2_idc}
The following identity holds:
\begin{equation}\label{eta2_ide}
(G\eta^{\st},\nabla_{\!\bk_s}\hh\tilde{\phi}_n)-\langle G \nabla_{\!\bk_s}\hh\eta^{\st} \rangle \;=\; (\rho\eta^{\tz},\bchi^{\so}).
\end{equation}
See \ref{identities} for proof.
\end{claim}

\begin{remark}\label{eta1_lin}
\noindent  With reference to Remark \ref{eta0_lin}, we find assuming $\tilde{f}_{\bk}(\bx) = F(\eps\hat{\bk})\hh\phi(\bx)$ that the respective solutions of~\eqref{eta_1} and~\eqref{eta_2} can be computed as  $\boldsymbol{\eta}^{\so}(\bx)=F(\eps\hat{\bk})\hh\boldsymbol{\zeta}^{\so}(\bx)$ and $\eta^{\st}(\bx)=F(\eps\hat{\bk})\hh\zeta^{\st}(\bx)$ , where $\boldsymbol{\zeta}^{\so}\!\in\big(\bar{H}^{1}_{p0}(Y)\big)^{d}$ and $\zeta^{\st}\!\in\bar{H}^{1}_{p0}(Y)$ are independent of~$\tilde{f}_\bk$ and uniquely solve 
\begin{eqnarray}\label{zeta0}
&\tilde{\lambda}_n\rho\hh\boldsymbol{\zeta}^{\so} + \nabla_{\!\bk_s}\!\sip \big(G(\nabla_{\!\bk_s} \boldsymbol{\zeta}^{\so}+\zeta^{\tz}\boldsymbol{I} )\big) +G\nabla_{\!\bk_s}\zeta^{\tz} +\frac{\rho}{\rho^{\tz}}\tilde{\phi}_n(\phi,\bchi^{\so}) \,=\, \frac{\rho}{\rho^{\tz}}\langle\phi\rangle\bchi^{\so},& \label{zeta_1} ~~~\\
\notag  \\
&-\tilde{\lambda}_n \rho \zeta^{\st} - \nabla_{\!\bk_s}\cdot(G\nabla_{\!\bk_s}\hh\zeta^{\st})~=~\rho\zeta^{\tz},& \label{eta_2}
\end{eqnarray}
\noindent with $G(\nabla_{\!\bk_s} \boldsymbol{\zeta}^{\so}+\zeta^{\tz}\boldsymbol{I} )$ and~$G\nabla_{\!\bk_s}\zeta^{\st}$ each being subject to the flux boundary conditions \eqref{gbcs}.\\
\end{remark}

In order to ``average'' the  $O(\eps^2)$ field equation~\eqref{epsp} with $m=2$, we insert the solutions~\eqref{ut0exp}, \eqref{ut1exp}, \eqref{ut2exp} and \eqref{ut3exp} into \eqref{epsp} and integrate~$(\eqref{epsp}, \tilde{\phi}_n)$ by parts using boundary conditions~\eqref{BCPDES} with $m=4$. In this way, we obtain  
\begin{eqnarray}\label{solvacon2}
-\bmu^{\st}\!:(i\hat{\bk})^4 u_0-\bth^{\so}:(i\hat{\bk})^3 u_1 -\big( \bmu^{\tz} : (i\hat{\bk})^2 + \hasig \rho^{\tz} \hat{\omega}^2 \big )~u_2 -(\bth^{\tz}\sip (i\hat{\bk})+\chsig \check{\omega}^2\rho^{\tz})~u_3~=~M_2(\hat{\bk},\hasig \check{\omega}^2)
\end{eqnarray}
\noindent where 
\begin{eqnarray}
\bmu^{\st} &\!\!\!=\!\!\!& \langle G\{\nabla_{\!\bk_s}\bchi^{\sh}+ \boldsymbol{I}\!\otimes\bchi^{\st}\} \rangle- \big \{\big( G \bchi^{\sh}\!\otimes \overline{\nabla_{\!\bk_s}\tilde{\phi}_n},1 \big )\big \}, \label{mu2} \\*[1mm]
M_2(\hat{\bk},\hasig \check{\omega}^2) &= &\hasig \check{\omega}^2 \big( \langle G\nabla_{\!\bk_s} \eta^{\st} \rangle - (G\eta^{\st},\nabla_{\!\bk_s}\tilde{\phi}_n) \big ) \cdot (i\hat{\bk}) \notag \\
&& +\; \big(\langle G\{\nabla_{\!\bk_s} \boldsymbol{\eta}^{\so}+\eta^{\tz}\boldsymbol{I} \}\rangle - (G \{\boldsymbol{\eta}^{\so} \!\otimes \overline{\nabla_{\!\bk_s}\tilde{\phi}_n}\},1) \big ) \!:\! (i\hat{\bk})^2, \label{force2} 
 \end{eqnarray}

\begin{claim}\label{mu2r}
For any~$\bk_s\!\in\overline{\mathcal{B}}$, effective tensor $\bmu^{\st}$ is real-valued, i.e. $\bmu^{\st}\in\mathbb{R}^{d\times d \times d \times d}$. See \ref{Effcoef} for proof.
\end{claim}

\begin{claim}\label{eta1_id}
We have the following identity 
\begin{equation}
\langle \{G\nabla_{\!\bk_s} \boldsymbol{\eta}^{\tiny{\mbox{(1)}}}+\eta^{\tz}\boldsymbol{I} \}\rangle -(G \{\boldsymbol{\eta}^{\tiny{\mbox{(1)}}} \!\otimes \overline{\nabla_{\!\bk_s}\tilde{\phi}_n}\},1) =(\tilde{f}_{\bk},\bchi^{\st}) + \frac{\langle\tilde{f}_{\bk}\rangle}{\rho^{\tz}}\{(\rho\bchi^{\tiny{\mbox{(1)}}}\!\otimes \overline{\bchi^{\tiny{\mbox{(1)}}}},1)\} +\frac{1}{\rho^{\tz}} \{\bth^{\tz}\otimes(\rho\eta^{\tz},\bchi^{\so})\},
\end{equation}
see \ref{identities} for proof. 
\end{claim}

With the effective equations \eqref{solvacon1} and \eqref{solvacon2} featuring $u_2$ and $u_3$ in place, we next evaluate the second-order counterparts of~\eqref{leadorder11} and~\eqref{leadorder12} depending on the triviality of $\bth^{\tz}$.

%--------------------------------------------------------------------------------------------------------------%
\subsubsection{Effective model for non-trivial $\bth^{\tz}$}\label{simpleadcas12}
%--------------------------------------------------------------------------------------------------------------%

\noindent When~$\tilde{f}_\bk\neq 0$ and~$\omega^2\!-\omega_n^2=O(\eps)$, we have~$|\check\sigma|=1$ due to Remark~\ref{scaling:choice}. In this case, we evaluate $\eqref{leadorder11}+\eps\eqref{solvacon1}+\eps^2\eqref{solvacon2}$ to obtain the second-order effective equation 
\begin{eqnarray}\label{3oee}
-\big(\bth^{\so}\!:(i\eps\hat{\bk})^3 \hh+\hh \bmu^{\tz} \!: (i\eps\hat{\bk})^2 \hh+\hh \bth^{\tz}\sip (i\eps\hat{\bk}) + \eps\hh\chsig\check{\omega}^2 \rho^{\tz} \big)\hh \langle \tilde{u} \rangle_{\!\rho} \;\stackrel{\eps^3}{=}\;  M_2'(\eps\hat{\bk}, \eps\chsig\check{\omega}^2),
\end{eqnarray}
where ``$\stackrel{\eps^3}{=}$'' denotes equality with an $O(\eps^3)$ residual, and the effective source term is given by
\begin{eqnarray}\label{force3}
M_2'(\eps\hat{\bk}, \eps\chsig\check{\omega}^2) &\!\!\!=\!\!\!& \langle \tilde{f}_{\bk}\rangle(\eps\hat\bk) - (\tilde{f}_{\bk},\bchi^{\so})\sip (i\eps\hat{\bk}) \notag + \eps\hh\chsig\check{\omega}^2 \big( \langle G\nabla_{\!\bk_s} \eta^{\st} \rangle - (G\eta^{\st},\nabla_{\!\bk_s}\tilde{\phi}_n) \big ) \cdot (i\eps\hat{\bk}) \notag\\
&&  +\; \big(\langle G\{\nabla_{\!\bk_s} \boldsymbol{\eta}^{\so}+\eta^{\tz}\boldsymbol{I} \}\rangle - (G \{\boldsymbol{\eta}^{\so} \!\otimes \overline{\nabla_{\!\bk_s}\tilde{\phi}_n}\},1) \big ) :(i\eps\hat{\bk})^2. \notag
\end{eqnarray}

In the absence of the source term, we find from~\eqref{leadorder11}, \eqref{leadorder12}, \eqref{solvacon1} and~\eqref{solvacon2} that a non-trivial effective solution in terms of~$u_m=\langle\tilde{u}_m\rangle_{\!\rho}$ ($m=\overline{0,3}$) is possible only if~$|\check\sigma|=1$. In this case, \eqref{3oee} with $\tilde{f}_{\bk}\!=\!0\,$ i.e.~$M_2'=0$ can be shown to describe a cubic approximation of the dispersion relationship near $(\bk_s,\omega_n\!>\!0)$ as
\begin{eqnarray}\label{drasymp21}
  \omega_n^2(\bk) \;=\; \omega_n^2 - \frac{1}{\rho^{\tz}}\hh i \bth^{\tz}\cdot(\eps \hat{\bk}) + \frac{1}{\rho^{\tz}}~\bmu^{\tz}:(\eps \hat{\bk})^2 +\frac{1}{\rho^{\tz}}\hh i \bth^{\so}:(\eps \hat{\bk})^3.
\end{eqnarray}

%--------------------------------------------------------------------------------------------------------------%
\subsubsection{Effective model for trivial $\bth^{\tz}$ }\label{simpleadcas22}
%--------------------------------------------------------------------------------------------------------------%

\noindent Assuming~$\tilde{f}_\bk\neq 0$ and $\omega^2\!-\omega_n^2=O(\eps^2)$ whereby~$|\hat\sigma|\!=\!1$, we obtain the second-order effective equation by evaluating $\eqref{leadorder12}+\eps\eqref{solvacon1}+\eps^2\eqref{solvacon2}$, namely  
\begin{eqnarray}\label{2oee}
-\big(\bmu^{\st}\!:(i\eps\hat{\bk})^4 \hh+\hh \bth^{\so}\!:(i\eps\hat{\bk})^3 \hh+\hh \bmu^{\tz} \!:(i\eps\hat{\bk})^2 + \eps^2\hasig\hat{\omega}^2 \rho^{\tz}\big)\hh \langle \tilde{u} \rangle_{\!\rho} ~\stackrel{\eps^3}{=}~ M_2''(\eps\hat{\bk}),
\end{eqnarray}
where 
\begin{eqnarray}\label{forceg}
M_2''(\eps\hat{\bk}) &\!\!\!=\!\!\!& \langle \tilde{f}_{\bk}\rangle(\eps\hat\bk) -(\tilde{f}_{\bk},\bchi^{\so})\sip (i\eps\hat{\bk})+\big(\langle G\{\nabla_{\!\bk_s} \boldsymbol{\eta}^{\so}+\eta^{\tz}\boldsymbol{I} \}\rangle - (G \{\boldsymbol{\eta}^{\so} \!\otimes \overline{\nabla_{\!\bk_s}\tilde{\phi}_n}\},1) \big ) :(i\eps\hat{\bk})^2.
\end{eqnarray}

\begin{claim} \label{specase2}
For  $\bk_s = \frac{1}{2}(\sum_j n_j\hh\be^j)$, $n_j\in\{-1,0,1\}$, which includes the origin and apexes of~$\overline{\mathcal{B}}$, Bloch functions $\tilde{\phi}_n(\bx) e^{ i \bk_s\cdot\bx}$, $\bchi^{\so}(\bx) e^{ i \bk_s\cdot\bx}$, $\bchi^{\st}(\bx) e^{ i \bk_s\cdot\bx}$ and~$\bchi^{\sh}(\bx) e^{ i \bk_s\cdot\bx}$ are real-valued up to a common factor $e^{i\varphi_0}$. Consequently, $\bth^{\so}\!=\!\boldsymbol{0}$ which in particular motivates our pursuit of the second-order approximation. See \ref{specases} for proof.
\end{claim}

When $\tilde{f}_{\bk}=0$, on the other hand, we deduce from~\eqref{leadorder11}, \eqref{leadorder12}, \eqref{solvacon1} and~\eqref{solvacon2} that~$|\hasig |=1$ in order to have a non-trivial solution. As a result, one finds that \eqref{2oee} with~$M_4=0$ furnishes a quartic approximation of the dispersion relationship near $(\bk_s,\omega_n\!>\!0)$ as 
\begin{eqnarray}\label{drasymp22}
\omega_n^2(\bk) \;=\; \omega_n^2 + \frac{1}{\rho^{\tz}}~\bmu^{\tz} \!:\!(\eps \hat{\bk})^2 + \frac{1}{\rho^{\tz}}\hh i \bth^{\so} \!:\!(\eps \hat{\bk})^3 -\frac{1}{\rho^{\tz}}~\bmu^{\st}\!:\!(\eps \hat{\bk})^4.
\end{eqnarray}

\begin{remark}
In the special case where $\tilde{f}_{\bk}(\bk,\bx)=1$ and $\bk_s = \frac{1}{2}(\sum_jn_j\hh\be^j)$, $n_j\in\{-1,0,1\}$, effective equation \eqref{2oee} reduces, thanks to Claim~\ref{specase1}, Claim~\ref{eta1_id} and Claim~\ref{specase2}, to 
\begin{eqnarray}\label{2oees}
-\big(\bmu^{\st}\!:(i\eps\hat{\bk})^4 \hh+\hh \bmu^{\tz}\!: (i\eps\hat{\bk})^2 + \eps^2\hasig\hat{\omega}^2 \rho^{\tz} \big) \langle \tilde{u} \rangle_{\!\rho} ~\stackrel{\eps^3}{=}~ M_2'''(\eps\hat{\bk}),
\end{eqnarray}
where 
\begin{eqnarray}\label{force3s}
M_2'''(\eps\hat{\bk}) &\!\!\!=\!\!\!& \langle 1\rangle  - \overline{(\bchi^{\so},1)}\sip (i\eps\hat{\bk})+ \big(\langle G\{\nabla_{\!\bk_s} \boldsymbol{\eta}^{\so}+\eta^{\tz}\boldsymbol{I} \}\rangle -(G \{\boldsymbol{\eta}^{\so} \!\otimes \overline{\nabla_{\!\bk_s}\tilde{\phi}_n}\},1) \big ) :(i\eps\hat{\bk})^2.
\end{eqnarray}
We note that~\eqref{2oees}--\eqref{force3s} can be reduced to the FF-FW effective model \cite{GMO19} upon: (i) accounting for the ``unfolding'' of~$Y$ according to the multicell homogenization approach \cite{ACG13b,Doss2008} and (ii) using~$\langle\tilde{u}\rangle$ instead of~$\langle \tilde{u} \rangle_{\!\rho}$ to define the ``mean'' motion. 
\end{remark}
%--------------------------------------------------------------------------------------------------------------%
\section{Repeated eigenvalues} \label{repeated_eig}
%--------------------------------------------------------------------------------------------------------------%
\noindent Let $\omega_n$ be an eigenfrequency of multiplicity $Q\!>\!1$, and let $n_q$ ($q\!=\!\overline{1,Q}$) be the indexes of associated eigenfunctions. 
 
 \begin{remark}\label{pqinQ}
In what follows, we assume $p,q,s\in\overline{1,Q}$, unless stated otherwise. Further, we will use the short-hand notation~$\sum_q$ for~$\sum_{q=1}^Q$. 
 \end{remark}
 
%--------------------------------------------------------------------------------------------------------------%
\subsection{Leading-order approximation}\label{rep_lead}
%--------------------------------------------------------------------------------------------------------------%
\noindent With reference to the eigenvalue problem \eqref{EVP}-\eqref{BCEVP}, the solution of \eqref{eps-2} in the vicinity of a repeated eigenfrequency $\omega_n$ can be decomposed as 
\begin{equation}\label{ut0expme}
\tilde{u}_0(\bx) = \sum_q u_{0q}~\tilde{\phi}_{n_{\nes q}}(\bx), \quad u_{0q} \in\mathbb{C},
\end{equation}
consistent with the definition~\eqref{umq} of $u_{0q}$. Then, by inserting \eqref{ut0expme}~in \eqref{eps-1} and integrating~$(\eqref{eps-1}, \tilde{\phi}_{n_p})$ by parts via the boundary conditions~\eqref{BCPDES} with $m=1$, we obtain the averaged $O(\eps^{-1}$) statement
\begin{equation}\label{solvacon-1me}
\sum_q \big(\bth_{pq}^{\tz}\!\cdot\hh i \hat{\bk}+\chsig \check{\omega}^2 \rho^{\tz}_{p}\delta_{pq} \big)\hh u_{0q} \;=\; 0, \qquad p=\overline{1,Q}
\end{equation}
where 
\begin{equation}\label{teta&rhome}
\bth_{pq}^{\tz} =  \langle G\nabla_{\!\bk_{s}} \tilde{\phi}_{n_{\nes q}} \rangle^{n_p} -  \overline{ \langle G\nabla_{\!\bk_{s}}\tilde{\phi}_{n_p}\rangle^{n_{\nes q}} }\quad \text{ and } \quad \rho^{\tz}_p = \langle \rho\tilde{\phi}_{n_p}\rangle^{n_p}.
\end{equation}
For convenience, system of equations~\eqref{solvacon-1me} can be expressed in the matrix form as 
\begin{equation} \label{solconmat1}
(\boldsymbol{A}^{\tz}(\hat{\bk}) + \chsig \check{\omega}^2\boldsymbol{D})\hh \bu_0 \;=\; \boldsymbol{0},
\end{equation}
where 
\begin{equation}\label{ADU0}
A^{\tz}_{pq}(\hat\bk) = \bth_{pq}^{\tz}\cdot\hh i \hat{\bk}, \qquad D_{pq}=\rho^{\tz}_p\delta_{pq}.
\end{equation}

\begin{remark}\label{Aherm}
Effective vectors $\bth_{qq}^{\tz}$ are imaginary-valued,  i.e.~$\bth_{qq}^{\tz}\in i\mathbb{R}^d$. Coefficient matrix~$\boldsymbol{D}\!\in\mathbb{R}^{Q\times Q}$ is diagonal, and $\boldsymbol{A}^{\tz}\!\in\mathbb{C}^{Q\times Q}$ is Hermitian.
\end{remark}

On the basis of \eqref{ut0expme}--\eqref{solvacon-1me}, we can solve the~$O(\eps^{-1})$ field equation~\eqref{eps-1} in terms of~$\tilde{u}_1$ as 
\begin{equation}\label{ut1expme}
\tilde{u}_1(\bx) =\sum_q \big(u_{0q}~\bchi^{\so}_q(\bx)\cdot (i\hat{\bk}) + u_{1q}~\tilde{\phi}_{n_{\nes q}}(\bx)\big), \quad u_{1q} \in\mathbb{C},
\end{equation}
where $\bchi_q^{\so}\in\big(\bar{H}^{1}_{p0}(Y)\big)^d$ solves uniquely the unit cell problem
\begin{equation}\label{chi1q}
\tilde{\lambda}_n\hh\rho\hh\bchi^{\so}_q + \nabla_{\!\bk_s}\!\sip \big(G(\nabla_{\!\bk_s}\hh\bchi^{\so}_q+\tilde{\phi}_{n_{\nes q}} \boldsymbol{I})\big) + G\nabla_{\!\bk_s} \tilde{\phi}_{n_{\nes q}} - \sum_{s} \frac{\rho}{\rho^{\tz}_s}\tilde{\phi}_{n_s}\bth^{\tz}_{sq}~=~\boldsymbol{0},
\end{equation} 
subject to the flux boundary conditions \eqref{gbcs} in terms of $G(\nabla_{\!\bk_s}\hh\bchi^{\so}_q+\tilde{\phi}_{n_{\nes q}} \boldsymbol{I})$.

We next consider the $O(1)$ field equation~\eqref{eps0}. On recalling~\eqref{ut0expme} and~\eqref{ut1expme} and integrating~$(\eqref{eps0}, \tilde{\phi}_{n_p})$ by parts via the boundary conditions~\eqref{BCPDES} with $m=2$, we obtain the averaged $O(1)$ statement
\begin{equation}\label{solvacon0me}
- \sum_q \Big(\big( \bmu^{\tz}_{pq} : (i\hat{\bk})^2 + \hasig \hat{\omega}^2 \rho^{\tz}_p\delta_{pq} \big)\hh u_{0q}+(\bth^{\tz}_{pq}\cdot (i\hat{\bk})+\chsig \check{\omega}^2\rho^{\tz}_p\delta_{pq})\hh u_{1q} \Big) \;=\; \langle \tilde{f}_{\bk}\rangle^{n_p},
\end{equation}
where
\begin{equation}\label{mu0pq}
\bmu^{\tz}_{pq} = \langle G\{\nabla_{\!\bk_s}\bchi^{\so}_q+ \tilde{\phi}_{n_{\nes q}}\boldsymbol{I}\} \rangle^{n_p} - \big \{\big( G \bchi^{\so}_q\!\otimes \overline{\nabla_{\!\bk_s}\tilde{\phi}_{n_p}},1 \big)\big\}. 
\end{equation}
We rewrite~\eqref{solvacon0me} in the matrix form as
\begin{equation} \label{solconmat2}
-(\boldsymbol{B}^{\tz}(\hat\bk) + \hasig \hat{\omega}^2 \boldsymbol{D})\hh \bu_0 \hh-\hh (\boldsymbol{A}^{\tz}(\hat\bk) + \chsig \check{\omega}^2 \boldsymbol{D})\hh \bu_1 \;=\; \mathbf{f}_0
\end{equation}
where 
\begin{equation}\label{BDU1}
B^{\tz}_{pq}(\hat\bk) = \bmu_{pq}^{\tz}:(i\hat{\bk})^2, \qquad \text{f}_{0p} = \langle {\tilde{f}_{\bk}} \rangle^{n_p}.% \quad r,q\in\overline{1,Q}. 
\end{equation}

\begin{claim}\label{Bherm}
Matrix $\boldsymbol{B}^{\tz}\!\in\mathbb{C}^{Q\times Q}$ is Hermitian. See~\ref{Effcoef} for proof.
\end{claim}

%--------------------------------------------------------------------------------------------------------------%
\subsubsection{Eigenfunction basis }\label{EFB}
%--------------------------------------------------------------------------------------------------------------%
\noindent For a fixed direction $\hat{\bk}/\|\hat\bk\|$, let  $\boldsymbol{P}\!\in\!\mathbb{C}^{Q\times Q}$ denote the matrix of eigenvectors associated with the generalized eigenvalue problem
\begin{equation}\label{gep1}
\boldsymbol{A}^{\tz}\boldsymbol{v} \;=\; \tau \boldsymbol{D} \hh \boldsymbol{v}. 
\end{equation}
In this setting, we conveniently introduce the ``recombined'' eigenfunctions $\tilde{\psi}_q$ as
\begin{equation}\label{psi0}
\tilde{\psi}_q = \sum_s P_{sq} \,\tilde{\phi}_{n_s}, \quad q\in\overline{1,Q}.
\end{equation}
Then, by taking the eigenfunctions $\{\tilde{\psi}_q\}$ as the projection basis in~\eqref{ut0expme} and~\eqref{teta&rhome} instead of~$\{\tilde{\phi}_{n_q}\}$, we find that
\begin{equation} \label{Adiag}
\boldsymbol{A}^{\tz}(\hat\bk) = \text{diag}(0,\dots,0,\tau_{N_0+1}, \tau_{N_0+2},\dots, \tau_{Q}) 
\end{equation}
where $\tau_q=\bth^{\tz}_{qq} \cdot  i \hat{\bk}$ and~$0\leqslant N_0\leqslant Q$ is the number of trivial diagonal entries of $\boldsymbol{A}^{\tz}$. In this setting, we also define the sub-matrices $\bar{\boldsymbol{B}}^{\tz}(\hat\bk)\in\mathbb{C}^{N_0\times N_0}$ and~$\bar{\boldsymbol{D}}\in\mathbb{R}^{N_0\times N_0}$ such that 
\[
\bar{B}^{\tz}_{pq} \;=\; B^{\tz}_{pq} \quad \text{and}\quad  \bar{D}_{pq} \;=\; D_{pq}, \qquad p,q\in\overline{1,N_0}.
\]

\begin{claim}\label{Aimag}
For $\bk_s=\frac{1}{2}(\sum_j n_j\be^j)$, $n_j\in\{-1,0,1\}$ which include the origin and apexes of the first Brillouin zone, Bloch eigenfunctions $\tilde{\phi}_{n_p}e^{ i \bk_s\cdot\bx}$ have constant phase and can be taken as real-valued. In this case, vector $\bth_{pq}^{\tz}\in\mathbb{R}^d$ and $\bth_{pp}^{\tz}=\boldsymbol{0}$, whereby  $\boldsymbol{A}^{\tz}$ is imaginary-valued and skew-symmetric. Consequently, the nonzero eigenvalues of $\boldsymbol{A}^{\tz}$ consist of pairs $\{\tau,-\tau\}$ whose respective eigenvectors are complex conjugates of each other. See~\ref{specases} for proof.
\end{claim}

When~$N_0\!>\!0$, we denote by $\bar{\boldsymbol{P}}$ the matrix of eigenvectors of the generalized eigenvalue problem 
\begin{equation}\label{gep2}
\bar{\boldsymbol{B}}^{\tz} \boldsymbol{v} \;=\; \tau \bar{\boldsymbol{D}} \hh \boldsymbol{v},
\end{equation} 
and we define the eigenfunction basis $\tilde{\psi}'_q$ as
\begin{equation}\label{Psi}
\tilde{\psi}'_q \;=\; \left\{\begin{array}{rl}
\textstyle{\sum_{s=1}^{N_0}} \bar{P}_{sq} \, \tilde{\psi}_{s}, & q \leqslant N_0 \\*[1mm]
\tilde{\psi}_q, &  q>N_0
\end{array}\right. .
\end{equation}

\begin{remark}\label{psiortho}
Eigenfunctions $\tilde{\psi}'_q$ ($q=\overline{1,Q}$) are $\rho$-orthogonal in the sense of~\eqref{evpo}, see~\ref{Effcoef} for proof. For simplicity of discussion, we hereon relabel $\tilde{\psi}'_q$ as~$\tilde{\phi}_{n_q}$. In this setting, we have
\begin{equation} \label{Bdiag}
\bar{\boldsymbol{B}}^{\tz} \;=\; \text{diag}(\bmu^{\tz}_{11}\!:\! (i\hat{\bk})^2,\dots,\bmu^{\tz}_{N_0N_0} \!:\! (i\hat{\bk})^2).
\end{equation}
\end{remark} 

%--------------------------------------------------------------------------------------------------------------%
\subsubsection{Additional scaling}\label{scal2}
%--------------------------------------------------------------------------------------------------------------%

\noindent Depending on the perturbation direction, certain non-zero diagonal entries of~$\boldsymbol{A}^{\tz}(\hat\bk)$ in \eqref{Adiag} can become vanishingly small, namely $\tau_{q}=o(1)$ for some~$q$. In the context of Section~\ref{leador}, for instance, this  situation would correspond to directions~$\hat\bk/\|\hat\bk\|$ for which~$\bth^{\tz}\sip (i\hat{\bk})=o(1)$. To account for such situations, we decompose~$\boldsymbol{A}^{\tz}$ as
\begin{eqnarray}\label{addscal}
\boldsymbol{A}^{\tz}(\hat\bk) &\!\!\!=\!\!\!& \text{diag} (0,\dots,0,\underbrace{\tau_{N_0+1},\dots,\tau_{N}}_{O(\eps)},\underbrace{\tau_{N+1},\dots \tau_{Q}}_{O(1)}) ~=~  \boldsymbol{\dot{A}}^{\tz}\!(\hat\bk) + \eps\hh\boldsymbol{\ddot{A}}^{\tz}\!(\hat\bk) \label{addscalb}\\
\boldsymbol{\dot{A}}^{\tz}\!(\hat\bk)  &\!\!\!=\!\!\!& \text{diag} (0,\dots,0,\underbrace{\tau_{N+1},\dots \tau_{Q}}_{O(1)}) ~\in~\mathbb{R}^{Q\times Q}, \\
\boldsymbol{\ddot{A}}^{\tz}\!(\hat\bk)  &\!\!\!=\!\!\!& \text{diag} (0,\dots,0, \underbrace{\eps^{-1}\tau_{N_0+1},\dots,\eps^{-1}\tau_{N}}_{O(1)},0,\dots,0),  ~\in~\mathbb{R}^{Q\times Q} \label{addscalc}
\end{eqnarray}
and we \emph{carry over} thus incurred $O(\eps)$ residual in~\eqref{solconmat1} to~\eqref{solconmat2} as
\begin{eqnarray}
(\boldsymbol{\dot{A}}^{\tz}\!(\hat\bk) + \chsig \check{\omega}^2\boldsymbol{D}) \hh \bu_0 &\!\!\!=\!\!\!&\boldsymbol{0},  \label{solconmat11} \\
-(\boldsymbol{B}^{\tz}(\hat\bk) + \boldsymbol{\ddot{A}}^{\tz}\! (\hat\bk) + \hasig \hat{\omega}^2 \boldsymbol{D}) \hh \bu_0 \hh-\hh (\boldsymbol{\dot{A}}^{\tz}\! (\hat\bk) +\chsig \check{\omega}^2 \boldsymbol{D}) \hh \bu_1 &\!\!\!=\!\!\!& \mathbf{f}_0.  \label{solconmat21}
\end{eqnarray}

On the basis of~\eqref{solconmat11}--\eqref{solconmat21}, we next pursue a family of the first-order effective field equations (in prescribed direction $\hat{\bk}/\|\hat\bk\|$) as controlled by: (i) proximity of the driving frequency $\omega^2$ to~$\tilde{\lambda}_n$ (see Remark~\ref{scaling:choice}), and~(ii) the nature of $\boldsymbol{A}^{\tz}(\hat\bk)$ according to~\eqref{addscal}--\eqref{addscalc}.

%--------------------------------------------------------------------------------------------------------------%
\subsubsection{Effective solution for full-rank $\boldsymbol{A}^{\tz}$ when $\boldsymbol{\ddot{A}}^{\tz}\!\!=\boldsymbol{0}$}\label{Afullrank}
%--------------------------------------------------------------------------------------------------------------%

\noindent We first consider the case where $\text{rank}(\boldsymbol{A}^{\tz}(\hat{\bk}))=\!Q$ and $\boldsymbol{\ddot{A}}^{\tz}\!\!=\boldsymbol{0}$. With reference to~\eqref{Adiag}, this specifically implies that $\tau_q=\bth^{\tz}_{qq} \cdot  i \hat{\bk}=O(1)$, $q=\overline{1,Q}$. Letting further $\tilde{f}_\bk\neq\boldsymbol{0}$ and $\omega^2\!-\omega_n^2=O(\eps)$ so that $|\check\sigma|=1$ by Remark~\ref{scaling:choice},  we find from the~$O(\eps^{-1})$ statement~\eqref{solconmat11} that~$\bu_0=\boldsymbol{0}$ thanks to Remark~\ref{warning!}. From~\eqref{solconmat21}, we then obtain the leading-order model 
\begin{equation} \label{leadsolfr}
- (\boldsymbol{A}^{\tz}(\hat\bk) +\chsig \check{\omega}^2 \boldsymbol{D}) \bu_1 ~=~ \mathbf{f}_0.
\end{equation}

In the absence of the driving source $\tilde{f}_{\bk}$, the existence of a non-trivial solution to~\eqref{solconmat11}--\eqref{solconmat21} also requires that~$|\check\sigma|=1$. As a result, \eqref{solconmat11} constitutes a generalized eigenvalue problem (GEP) whose eigenvalues
\begin{equation}\label{drasympta1}
\omega_{n_{\nes q}}^2(\bk) \;=\; \omega_n^2 - \frac{1}{\rho^{\tz}_q}\hh i \bth^{\tz}_{qq}\cdot(\eps \hat{\bk}),
\end{equation}
describe the leading-order, \emph{linear} dispersion relationship in direction $\hat{\bk}$.

%--------------------------------------------------------------------------------------------------------------%
\subsubsection{Effective solution for near-trivial $\boldsymbol{A}^{\tz}$}\label{Antrivial}
%--------------------------------------------------------------------------------------------------------------%

\noindent When $\boldsymbol{A}^{\tz}(\hat\bk) = \eps\hh \boldsymbol{\ddot{A}}^{\tz}$ i.e. $\boldsymbol{\dot{A}}^{\tz}\!=\boldsymbol{0}$, we consider the situation where $\tilde{f}_\bk\neq0$ and $\omega^2\!-\omega_n^2=O(\eps^2)$ so that $|\hat{\sigma}|=1$. In this case \eqref{solconmat11} is satisfied identically, and we find from \eqref{solconmat21} that the leading-order solution $\bu_0$ solves
\begin{equation}\label{leadsol21}
-(\boldsymbol{B}^{\tz}(\hat\bk) + \boldsymbol{\ddot{A}}^{\tz}\! (\hat\bk) +\hasig \hat{\omega}^2 \boldsymbol{D}) \bu_0 ~=~ \mathbf{f}_0. 
\end{equation}
In the degenerate case when $\boldsymbol{A}^{\tz}= \eps\hh \boldsymbol{\ddot{A}}^{\tz}\!=\boldsymbol{0}$, \eqref{leadsol21} becomes 
\begin{equation}\label{leadsol22}
-(\boldsymbol{B}^{\tz}(\hat\bk) +\hasig \hat{\omega}^2 \boldsymbol{D}) \bu_0~=~\mathbf{f}_0. 
\end{equation}
In this case we conveniently let~$P_{sq}=\delta_{sq}$ in~\eqref{psi0}, and we have~$N_0=Q$ whereby $\boldsymbol{B}^{\tz}=\boldsymbol{\bar{B}}^{\tz}$ becomes diagonal according to~\eqref{Bdiag}. 

When $\tilde{f}_\bk=0$, the existence of a non-trivial solution requires that~$\check\sigma=0$ i.e.~$|\hat\sigma|=1$. In this case, the leading-order approximation of the dispersion relationships $\omega_{n_q}(\bk)$, $q\!=\!\overline{1,Q}$ is obtained by solving the eigenvalue problem $(\boldsymbol{B}^{\tz}+ \boldsymbol{\ddot{A}}^{\tz})\bv = \tau \boldsymbol{D} \bv$. When $\boldsymbol{\ddot{A}}^{\tz}\!$ vanishes, the solution is given explicitly by  
\begin{equation}\label{drasympta2}
\omega_{n_{\nes q}}^2(\bk) \;=\; \omega_n^2 + \frac{1}{\rho^{\tz}_q}~\bmu^{\tz}_{qq}:(\eps \hat{\bk})^2,
\end{equation}
thanks to the fact that~$\boldsymbol{B}^{\tz}$ is diagonal in this case.
 
%--------------------------------------------------------------------------------------------------------------%
\subsubsection{Effective solution for partial rank $\boldsymbol{A}^{\tz}$ }\label{Apartialrank}
%--------------------------------------------------------------------------------------------------------------%

\noindent We next assume that $\boldsymbol{A}^{\tz}$ has a partial rank, i.e.~$0\!<\!N_0\!<\!Q$. Letting~$\tilde{f}_\bk\neq 0$ and $\omega^2-\omega_n^2=O(\eps^2)$, we have~$|\hat\sigma|=1$. Thanks to the fact that~$\boldsymbol{{A}}^{\tz}$ is diagonal due to~\eqref{Adiag}, the last~$Q-N$ components of~$\bu_0$ must vanish by enforcing~\eqref{solconmat11} to the leading order. By virtue of this result and~\eqref{solconmat21}, we find that  
\begin{eqnarray} 
-\sum_{q=1}^N (B^{\tz}_{pq}(\hat\bk) + \ddot{A}^{\tz}_{pq}(\hat\bk) + \hasig \hat{\omega}^2 D_{pq})\hh u_{0q}~=&  \text{f}_{0p}, &  p\in\overline{1,N}, \label{leadsol31} \\*[-1mm]
u_{0p}~=&0, &  p\in\overline{N+1,Q}. \label{leadsol31xxx}
\end{eqnarray}
When $\tilde{f}_\bk=0$, we enable a non-trivial solution to~\eqref{solconmat11}--~\eqref{solconmat21} in terms of~$\bu_0$ by taking~$|\hat\sigma|=1$. In this case, \eqref{leadsol31} with~$\text{f}_{0p}=0$ constitute a GEP yielding the leading-order approximation the first~$N$ dispersion branches~$\omega_{n_q}(\bk)$, $q\!=\!\overline{1,N}$.

Letting~$\tilde{f}_\bk\neq 0$ and $\omega^2-\omega_n^2=O(\eps)$, on the other hand, we have~$|\check\sigma|=1$ whereby $\bu_0=\boldsymbol{0}$ thanks to~\eqref{solconmat11}. From~\eqref{solconmat21}, we accordingly find that~$\bu_1$ solves
\begin{equation} \label{leadsolfro}
-(\boldsymbol{\dot{A}}^{\tz}\!(\hat\bk) +\chsig \check{\omega}^2 \boldsymbol{D}) \bu_1 ~=~ \mathbf{f}_0, 
\end{equation}
to the leading order (specifically, we discard the~$O(\eps)$ residual in~\eqref{solconmat21} by superseding~$\boldsymbol{A}^{\tz}$ with~$\boldsymbol{\dot{A}}^{\tz}$). Assuming~$\tilde{f}_\bk=0$, we are now left with exposing the leading-order behavior the last~$Q\!-\!N$ dispersion branches~$\omega_{n_q}(\bk)$, $q\!=\!\overline{N\!+\!1,Q}$. In this case we must set~$|\check\sigma|=1$ because all dispersion branches permitting the~$\hat\sigma$-description are already given by~\eqref{leadsol31} with~$\text{f}_{0p}=0$. This yields the sought approximation via~\eqref{leadsolfro} with~$\mathbf{f}_0=\boldsymbol{0}$ as 
\begin{equation}\label{drasympta1o}
\omega_{n_{\nes q}}^2(\bk) \;=\; \omega_n^2 - \frac{1}{\rho^{\tz}_q}\hh i \bth^{\tz}_{qq}\cdot(\eps \hat{\bk}), \qquad q=\overline{N\!+\!1,Q}.
\end{equation}

%--------------------------------------------------------------------------------------------------------------%
\subsection{First-order correctors}\label{secondor}
%--------------------------------------------------------------------------------------------------------------%

\noindent With the aid of the averaged~$O(\eps^{-1})$ statement~\eqref{solvacon0me}, one may solve the~$O(1)$ field equation~\eqref{eps0} as
\begin{equation}\label{ut2expme}
\tilde{u}_2(\bx) = \sum_q (u_{0q}~\bchi^{\st}_q(\bx)\!:\!(i\hat{\bk})^2 + u_{1q}~\bchi^{\so}_q(\bx)\sip (i\hat{\bk}) + u_{2q}~\tilde{\phi}_{n_{\nes q}}(\bx)) + \eta^{\tz}(\bx), \qquad u_{2q} \in\mathbb{C},
\end{equation}
where $\bchi^{\st}_q\!\in\big(\bar{H}^{1}_{p0}(Y)\big)^{d\times d}$ and~$\eta^{\tz}\!\in\bar{H}^{1}_{p0}(Y)$ solve uniquely the respective equations 
\begin{eqnarray}\label{chi2q}
\tilde{\lambda}_n\hh\rho\hh\bchi^{\st}_q + \nabla_{\!\bk_s}\!\sip \big(G(\nabla_{\!\bk_s}\hh\bchi^{\st}_q+\tilde{\phi}_{n_{\nes q}}\{\boldsymbol{I}\!\otimes\bchi^{\so}_q\}')\big) + \{G(\nabla_{\!\bk_s}\hh\bchi^{\so}_q+\tilde{\phi}_{n_{\nes q}} \boldsymbol{I})\}
= \sum_s \Big( \frac{\rho}{\rho^{\tz}_{s}}\{\bth^{\tz}_{sq}\!\otimes\bchi^{\so}_s\} + \frac{\rho}{\rho^{\tz}_s}\tilde{\phi}_{n_s}\hh\bmu^{\tz}_{sq} \Big),
\end{eqnarray}
\begin{equation}\label{eta02}
-\tilde{\lambda}_n\rho\hh\eta^{\tz} - \nabla_{\!\bk_s}\!\sip \big(G\nabla_{\!\bk_s} \eta^{\tz}\big)~=~\tilde{f}_{\bk}-\sum_s \frac{\rho}{\rho^{\tz}_s}\langle \tilde{f}_{\bk}\rangle^{n_s} \hh\tilde{\phi}_{n_s}.
\end{equation}
with $G(\nabla_{\!\bk_s}\hh\bchi^{\st}_q+\{\boldsymbol{I}\!\otimes\bchi^{\so}_q\}')$ and $G\nabla_{\!\bk_s}\eta^{\tz}$ satisfying the flux boundary conditions \eqref{gbcs}.\\

\begin{claim}\label{eta_chi1_idmm}
We have the following identity 
\begin{equation}\label{eta0_chi1_idmm}
(G\eta^{\tz},\nabla_{\!\bk_s}\hh\tilde{\phi}_{n_{\nes q}})-\langle G \nabla_{\!\bk_s}\hh\eta^{\tz} \rangle^{n_{\nes q}} \;=\;  (\tilde{f}_{\bk},\bchi^{\so}_q),
\end{equation}
See Appendix section \ref{identities} for proof.
\end{claim}

\begin{remark}\label{eta02_lin}
When $\tilde{f}_{\bk}(\bx) = F(\eps\hat{\bk})\hh\phi(\bx)$, the solution of~\eqref{eta02} is given by $\eta^{\tz}(\bx)=F(\eps\hat{\bk})\zeta^{\tz}(\bx)$, where $\zeta^{\tz}\!\in\bar{H}^{1}_{p0}(Y)$ uniquely solves
\begin{equation}\label{zeta0}
-\tilde{\lambda}_n\rho\hh\zeta^{\tz} - \nabla_{\!\bk_s}\!\sip \big(G\nabla_{\!\bk_s} \zeta^{\tz}\big)~=~\phi-\sum_s \frac{\rho}{\rho^{\tz}_s}\langle \phi\rangle^{n_s} \hh\tilde{\phi}_{n_s}, 
\end{equation}
subject to the flux boundary conditions \eqref{gbcs} in terms of~$G\nabla_{\!\bk_s}\zeta^{\tz}$.
\end{remark}

On substituting~\eqref{ut0expme} and~\eqref{ut1expme} into~\eqref{epsp} with~$m\!=\!1$, integrating~$(\eqref{epsp}, \tilde{\phi}_n)$ by parts via boundary conditions~\eqref{BCPDES} with $m=3$, and exploiting the result of Claim \ref{eta_chi1_idmm}, we obtain the averaged~$O(\eps)$ statement
\begin{equation}\label{solvacon1me}
- \sum_q \Big( \bth^{\so}_{pq}:(i\hat{\bk})^3 u_{0q} + \big( \bmu^{\tz}_{pq} : (i\hat{\bk})^2 + \hasig \hat{\omega}^2 \rho^{\tz}_p\delta_{pq} \big)~u_{1q}+(\bth^{\tz}_{pq}\cdot (i\hat{\bk})+\chsig \check{\omega}^2\rho^{\tz}_p\delta_{pq})~u_{2q} \Big)= \text{f}_{1p},
\end{equation}
\noindent where
\begin{eqnarray}\label{theta1pq}
\bth^{\so}_{pq} \;=\; \langle G\{\nabla_{\!\bk_s}\bchi^{\st}_q+ \boldsymbol{I}\!\otimes\bchi^{\so}_q\} \rangle^{n_p} - \big \{\big( G \bchi^{\st}_q\!\otimes \overline{\nabla_{\!\bk_s}\tilde{\phi}_{n_p}},1 \big)\big \}, \qquad \text{f}_{1p} \;=\;  -(\tilde{f}_{\bk},\bchi^{\so}_p)\cdot  (i \hat{\bk}). 
\end{eqnarray}
We can conveniently rewrite~\eqref{solvacon1me} in the matrix form as
\begin{equation} \label{solconmat3}
- \boldsymbol{A}^{\so}(\hat\bk) \hh\hh \bu_0-(\boldsymbol{B}^{\tz}(\hat\bk) +\hasig \hat{\omega}^2 \boldsymbol{D}) \hh \bu_1 - (\boldsymbol{A}^{\tz}(\hat\bk) +\chsig \check{\omega}^2 \boldsymbol{D}) \hh \bu_2 \;=\; \mathbf{f}_1, 
\end{equation}
\noindent where 
\begin{equation}\label{A1DU2}
A^{\so}_{pq} = \bth_{pq}^{\so}:(i\hat{\bk})^3.
\end{equation}

\begin{remark}
At this point, we recall the decomposition of~$\boldsymbol{A}^{\tz}$ according to~\eqref{addscal}--\eqref{addscalc}, and the fact that~\eqref{solconmat21} accordingly incurs an~$O(\eps)$ residual, manifest in the term $-\eps\boldsymbol{\ddot{A}}^{\tz}\!\bu_1$, that carries over to the next order of asymptotic approximation.
\end{remark}

On accounting in~\eqref{solconmat3} for the~$O(\eps)$ residual stemming from~\eqref{solconmat21}, we obtain the averaged statement 
\begin{equation}\label{solconmat31}
- \boldsymbol{A}^{\so}(\hat\bk)\hh\hh \bu_0-(\boldsymbol{B}^{\tz}(\hat\bk) +\boldsymbol{\ddot{A}}^{\tz}\!(\hat\bk)+\hasig \hat{\omega}^2 \boldsymbol{D})\hh \bu_1 - (\boldsymbol{\dot{A}}^{\tz}\! (\hat\bk) +\chsig \check{\omega}^2 \boldsymbol{D}) \hh \bu_2 ~=~ \mathbf{f}_1.
\end{equation}
which allows us to compute the first-order correction~$\bu_1$ (resp.~$\bu_2$) of the leading-order model~$\bu_0$ (resp.~$\bu_1$). With reference to sections~\S\ref{Afullrank}--\S\ref{Apartialrank}, we pursue such task for three canonical situations driven by the nature of~$\boldsymbol{A}^{\tz}$. Before proceeding, we conveniently denote by~$\langle\tilde{\bu}\rangle_{\!\rho}$ the effective solution vector collecting the left-hand sides in~\eqref{effsolud}, which gives   
\begin{equation}\label{solvec1}
\langle\tilde{\bu}\rangle_{\!\rho}(\eps\bk) = \sum_{m=0}^\infty \eps^{m-2} \hh \bu_m. 
\end{equation}
For brevity, we focus our attention on the effective equations only, noting that the respective approximations of the dispersion relationship can be uniformly obtained by: (i) setting the source term in the effective equation to zero, and (ii) solving te resulting GEP.

%--------------------------------------------------------------------------------------------------------------%
\subsubsection{Effective solution for full-rank $\boldsymbol{A}^{\tz}$ when $\boldsymbol{\ddot{A}}^{\tz}\!\!=\boldsymbol{0}$}\label{Afullrank1}
%--------------------------------------------------------------------------------------------------------------%

\noindent Assuming~$\text{rank}(\boldsymbol{A}^{\tz})=Q$ and~$\boldsymbol{A}^{\tz}=\boldsymbol{\dot{A}}^{\tz}$, we let~$\tilde{f}_\bk\neq 0$ and~$\omega^2\!-\omega_n^2=O(\eps)$. In this case~$|\check\sigma|=1$, $\bu_0=\boldsymbol{0}$, and by~\eqref{solconmat31} the first-order order corrector~$\bu_2$ solves 
\begin{equation}\label{secordcor1}
-\boldsymbol{B}^{\tz}(\hat\bk)\hh \bu_1 \hh-\hh (\boldsymbol{A}^{\tz}(\hat\bk) +\chsig \check{\omega}^2 \boldsymbol{D}) \bu_2 \;=\; \mathbf{f}_1,
\end{equation}
where~$\bu_1$ is given by~\eqref{leadsolfr}. Thanks to~\eqref{solvec1}, we can now evaluate~$\eqref{leadsolfr}+\eps\eqref{secordcor1}$ in order to obtain the first-order effective model
\begin{equation}\label{bojan1}
- (\boldsymbol{B}^{\tz}(\eps\hat\bk) + \boldsymbol{A}^{\tz}(\eps\hat\bk) + \eps\chsig\check{\omega}^2 \boldsymbol{D})\hh \langle\tilde{\bu}\rangle_{\!\rho} \;\stackrel{\eps^2}{=}\; \mathbf{f}_0 + \eps\hh\mathbf{f}_1,
\end{equation}
where the components of~$\mathbf{f}_0$ and~$\mathbf{f}_1$  are given respectively in~\eqref{BDU1} and~\eqref{theta1pq}. For completeness, one may note that~\eqref{bojan1} carries the same structure as its simple-eigenvalue counterpart~\eqref{3oee} when truncated to the first order.

%--------------------------------------------------------------------------------------------------------------%
\subsubsection{Effective solution for near-trivial $\boldsymbol{A}^{\tz}$ }\label{Atrivial1}
%--------------------------------------------------------------------------------------------------------------%

\noindent When $\boldsymbol{A}^{\tz}=\eps\hh\boldsymbol{\ddot{A}}^{\tz}$, $\tilde{f}_\bk\neq 0$ and $\omega^2\!-\omega_n^2=O(\eps^2)$ i.e.~$|\hat\sigma|=1$, the leading-order model is given by~$\bu_0$ solving~\eqref{leadsol21}. In this case we discard the second-order correction~$\bu_2$ in~\eqref{solconmat31}, which then yields 
\begin{equation}\label{secordcor2}
-\boldsymbol{A}^{\so}(\hat\bk)\hh \bu_0 \hh-\hh (\boldsymbol{B}^{\tz}(\hat\bk) +\boldsymbol{\ddot{A}}^{\tz}\! (\hat\bk) +  \hasig \hat{\omega}^2 \boldsymbol{D}) \bu_1~=~\mathbf{f}_1.
\end{equation}
From~$\eqref{leadsol21}+\eps\eqref{secordcor2}$, we obtain the first-order effective model 
\begin{equation}\label{bojan2}
-(\boldsymbol{A}^{\so}(\eps\hat\bk) + \boldsymbol{B}^{\tz}(\eps\hat\bk) + \boldsymbol{\ddot{A}}^{\tz}\! (\eps\hat\bk) + \eps^2\hasig \hat{\omega}^2 \boldsymbol{D}) \hh \langle\tilde{\bu}\rangle_{\!\rho} \;\stackrel{\eps^2}{=}\; \mathbf{f}_0 + \eps\hh \mathbf{f}_1.
\end{equation}
When $\boldsymbol{\ddot{A}}^{\tz}\!=\boldsymbol{0}$, \eqref{bojan2} carries the same structure as its simple-eigenvalue counterpart~\eqref{2oee} after truncation to the first order.

%--------------------------------------------------------------------------------------------------------------%
\subsubsection{Effective solution for partial rank $\boldsymbol{A}^{\tz}$ }\label{Apartialrank1}
%--------------------------------------------------------------------------------------------------------------%
\noindent When $\text{rank}(\boldsymbol{A}^{\tz})<Q$, $\tilde{f}_\bk\neq 0$ and $\omega^2\!-\omega_n^2=O(\eps)$, the leading-order solution~$\bu_1$ satisfies~\eqref{leadsolfro}, while the first-order corrector~$\bu_2$ solves~\eqref{secordcor1}. On the other hand, when $\omega^2-\lambda_n=O(\eps^2)$, the leading-order model~$\bu_0$ is given by~\eqref{leadsol31}--\eqref{leadsol31xxx}, while its corrector~$\bu_1$ satisfies
\begin{eqnarray}\label{secorccor2}
-\sum_{q=1}^{N}(B^{\tz}_{pq}(\hat\bk) +\ddot{A}^{\tz}_{pq}(\hat\bk)+\hasig \hat{\omega}^2 D_{pq})\hh u_{1q}~= &\text{f}_{1p}+ \sum_{q=1}^{N} A^{\so}_{pq}(\hat\bk)\hh u_{0q}, & p\in\overline{1,N},  \\
-\dot{A}^{\tz}_{pp}(\hat\bk)\hh u_{1p}~= & \text{f}_{0p}+ \sum_{q=1}^{N}B^{\tz}_{pq}(\hat\bk)\hh u_{0q}, & p\in\overline{N+1,Q}. \label{secorccor2x}
\end{eqnarray}
In this case, we obtain the first-order model as~$\langle\tilde{\bu}\rangle_{\!\rho} \stackrel{1}{=} \eps^{-2}\bu_0+\eps^{-1}\bu_1$, where the scaled summands are directly computable from~\eqref{leadsol31}--\eqref{leadsol31xxx} and~\eqref{secorccor2}--\eqref{secorccor2x} upon replacing~$\hat\bk$ and~$ \hat\sigma\hat\omega$ respectively by~$\eps\hh\hat\bk$ and~$ \eps^2\hat\sigma\hat\omega$.

%--------------------------------------------------------------------------------------------------------------%
\section{Cluster of nearby eigenvalues} \label{cluster_eig}
%--------------------------------------------------------------------------------------------------------------%

\begin{figure}[h!] 
\centering{\includegraphics[width=80mm]{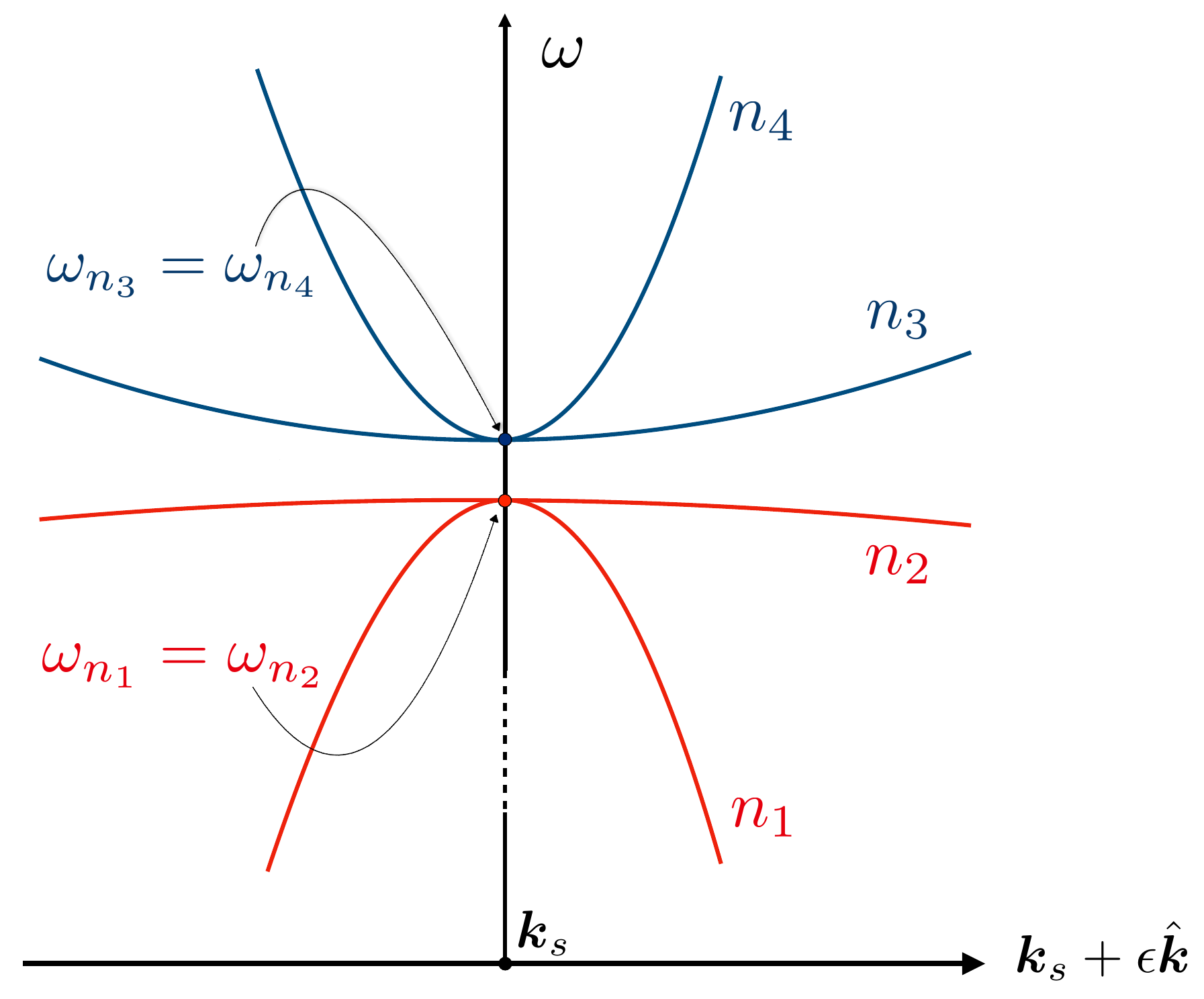}}\vspace*{-1mm}
\caption{Example cluster of nearby dispersion branches.} \label{fig2}
\end{figure}

\noindent We conclude the general analysis by letting the driving frequency be near a cluster of nearby eigenfrequencies $\{\omega_{n_q}\}$, $q=\overline{1,Q}$ depicted in~Fig. \ref{fig2}. This situation was originally considered in~\cite{GMO19} in an effort to handle the ``short asymptotic range''  exhibited by single- and repeated-eigenfrequency models within~$(\bk,\omega)$ regions characterized by closely spaced disperson curves. 
Our goal is to extend analysis in~\cite{GMO19} by: (i) permitting expansion about an arbitrary point~$(\bk_s,\omega)$, $\bk_s\in\overline{\mathcal{B}}$ and (ii) exposing the first-order correction of the leading-order model.\\

We let~$\bar{Q}$ be the number of distinct eigenvalues within set~$\{\omega_{n_q}\}$, and we denote by~$(\bk_s,\omega_{n_0})$ for some $n_0\in \{n_q\}$ the origin of asymptotic expansion in~\eqref{scal}.  In this setting, we conveniently redeploy the scaling parameter~$\eps=o(1)$ to quantify the ``smallness'' of distances between the neighboring eigenvalues by letting 
\begin{equation}\label{eig_scal}
\lambda_{n_{\nes q}} \:=\; \lambda_{n_0} - \eps \hh \gamma_q, \qquad q=\overline{1,Q}.
\end{equation}

\begin{remark}
Note that the logic behind such use of~$\eps$ is consistent with previous developments. Specifically, in sections~\S\ref{simple_eig} and \S \ref{repeated_eig}, we defined the size of the ``asymptotic box'' surrounding~$(\bk_s,\omega)$ as either~$O(\eps)^d\times O(\eps)$ or~$O(\eps)^d\times O(\eps^2)$ depending on (a) the driving frequency when~$\tilde{f}_\bk\neq 0$, and (b) the flatness of the $n$ dispersion branch for~$\tilde{f}_\bk= 0$. In the present case, by~\eqref{eig_scal} we ensure that such  ``asymptotic box'' captures all (relevant) dispersion surfaces in the cluster.
\end{remark}

With~\eqref{eig_scal} in place, we consider the local eigenfunction basis $\{\tilde{\phi}_{n_q}(\bk)\!\in\!H^1_{p0}(Y)\}$ that satisfies  
\begin{equation}\label{EVPclus}
-(\tilde{\lambda}_{n0}-\eps\hh\gamma_q) \rho(\bx) \hh \tilde{\phi}_{n_q} - \nabla_{\!\bk}\sip \big(G(\bx)\nabla_{\!\bk}~\tilde{\phi}_{n_q}\big) ~=~0 \quad\text{in~~}Y, \qquad q\!=\!\overline{1,Q}
\end{equation}
together with boundary conditions~\eqref{BCEVP}. As can be seen from~\eqref{EVPclus}, the current problem can be described as an ``almost repeated'' eigenvalue case, which allows us to take advantage of the foregoing developments.

With the insight on~$\tilde{u}$ solving~\eqref{RPDE} gained in Section~\ref{simple_eig} and Section~\ref{repeated_eig}, we skip intermediate steps and proceed by specifying ansatz~\eqref{sol_scal} up to~$m=2$ as
\begin{eqnarray}\label{utildeexpc}
\tilde{u}(\bx)  &\!\!\! \stackrel{\eps}{=}\!\!\!& \sum_q \Big( \eps^{-2} u_{0q}~\tilde{\phi}_{n_{\nes q}}(\bx) + \eps^{-1} u_{1q}~\tilde{\phi}_{n_{\nes q}}(\bx) + \eps^{-1}    u_{0q}~\bchi^{\so}_{q}(\bx)\cdot (i\hat{\bk} ) + u_{0q}~\bchi^{\st}_q(\bx):(i\hat{\bk})^2 \notag \\*[-2mm]
&& \qquad~~ +\, u_{1q}~\bchi^{\so}_q(\bx)\cdot (i\hat{\bk}) + u_{2q}~\tilde{\phi}_{n_{\nes q}}(\bx) \Big) \;+\; \sum_{q=1}^{\bar{Q}} \eta^{\tz}_q(\bx), \qquad u_{0q}, u_{1q}, u_{2q} \in \mathbb{C},
\end{eqnarray}

\noindent where $\eta^{\tz}_q\!\in\!\bar{H}^{1}_{p0}(Y)$, $\bchi^{\so}_{q}\!\in\!\big(\bar{H}^{1}_{p0}(Y)\big)^d$ and $\bchi^{\st}_{q}\!\in\! \big(\bar{H}^{1}_{p0}(Y)\big)^{d\times d}$ uniquely solve the respective equations
\begin{eqnarray}\label{chi1c}
&&\hspace*{-8mm}\tilde{\lambda}_{n_{\nes q}}\hh\rho\hh\bchi^{\so}_q + \nabla_{\!\bk_s}\!\sip \big(G(\nabla_{\!\bk_s}\hh\bchi^{\so}_q+\tilde{\phi}_{n_{\nes q}} \boldsymbol{I})\big) + G\nabla_{\!\bk_s} \tilde{\phi}_{n_{\nes q}} ~=~\sum_s \frac{\rho}{\rho^{\tz}_s}\tilde{\phi}_{n_s}\bth^{\tz}_{sq},\label{chi1c} \\
&&\hspace*{-8mm}\tilde{\lambda}_{n_{\nes q}}\hh\rho\hh\bchi^{\st}_q + \nabla_{\!\bk_s}\!\sip \big(G(\nabla_{\!\bk_s}\hh\bchi^{\st}_q+\tilde{\phi}_{n_{\nes q}}~\{\boldsymbol{I}\!\otimes\bchi^{\so}_q\}')\big) + \{G(\nabla_{\!\bk_s}\hh\bchi^{\so}_q+\tilde{\phi}_{n_{\nes q}} \boldsymbol{I})\} ~=~\sum_s \Big( \frac{\rho}{\rho^{\tz}_{s}}\{\bth^{\tz}_{sq}\!\otimes\bchi^{\so}_s\} + \frac{\rho}{\rho^{\tz}_s}\tilde{\phi}_{n_s}\hh\bmu^{\tz}_{sq} \Big), \notag \\ \label{chi2c}\\
&&\hspace*{-8mm}-\tilde{\lambda}_{n_{\nes q}}\rho\hh\eta^{\tz}_q - \nabla_{\!\bk_s}\!\sip \big(G\nabla_{\!\bk_s} \eta^{\tz}_q\big)~=~\tilde{f}_{\bk}-\sum_s \frac{\rho}{\rho^{\tz}_s}\langle \tilde{f}_{\bk}\rangle^{n_s} \hh\tilde{\phi}_{n_s}, \label{eta03}
\end{eqnarray} 
\noindent with $G(\nabla_{\!\bk_s}\hh\bchi^{\so}_q\!+\tilde{\phi}_{n_{\nes q}} \boldsymbol{I})$, $G(\nabla_{\!\bk_s}\hh\bchi^{\st}_q\!+\{\boldsymbol{I}\!\otimes\bchi^{\so}_q\}')$ and $G\nabla_{\!\bk_s}\eta^{\tz}_q$ each being subject to the flux boundary conditions~\eqref{gbcs}. Note also that Remark \ref{eta02_lin} still applies in this case. Our goal is then to find the coupled effective equations satisfied by $u_{0q}, u_{1q}$ and~$u_{2Ã¥q}$. To this end, we (i)~insert \eqref{utildeexpc} in \eqref{RPDE}; (ii)~integrate~$\langle \eqref{RPDE}\rangle_{\!\rho}^{n_p}$ by parts using the boundary conditions \eqref{BCPDE}, and (iii)~expand the result in powers of $\eps$ as
\begin{eqnarray}
O(\eps^{-1}): &&(\boldsymbol{A}^{\gamma}(\hat\bk) + \chsig \check{\omega}^2\boldsymbol{D}) \hh \bu_0  \;=\; \boldsymbol{0}, \label{leadeqc}\\
O(1): &&-(\boldsymbol{B}^{\tz}(\hat\bk)  +\hasig \hat{\omega}^2 \boldsymbol{D}) \hh \bu_0 - (\boldsymbol{A}^{\gamma}(\hat\bk)  +\chsig \check{\omega}^2 \boldsymbol{D}) \hh \bu_1 \;=\; \mathbf{f}_0, \label{seceqc} \\
O(\eps): && - \boldsymbol{A}^{\so}(\hat\bk) \hh \bu_0-(\boldsymbol{B}^{\tz}(\hat\bk)  +\hasig \hat{\omega}^2 \boldsymbol{D}) \hh \bu_1 - (\boldsymbol{A}^{\gamma}(\hat\bk)  +\chsig \check{\omega}^2 \boldsymbol{D}) \hh \bu_2 \;=\; \mathbf{f}_1, \label{thirdeqc}
\end{eqnarray}
where 
\begin{equation}\label{Agam}
 \boldsymbol{A}^{\gamma}(\hat\bk)  \;=\; \boldsymbol{A}^{\tz}(\hat\bk) \hh+\hh \boldsymbol{\Gamma}\boldsymbol{D}, \qquad \boldsymbol{\Gamma}_{pq} = \delta_{pq} \gamma_q 
 \end{equation}
accounts for the eigenvalue separations in~\eqref{EVPclus}, while $\boldsymbol{A}^{\tz}$, $\boldsymbol{B}^{\tz}$, $\boldsymbol{A}^{\so}$, $\boldsymbol{D}$, $\mathbf{f}_0$ and $\mathbf{f}_1$ are given by \eqref{ADU0}, \eqref{BDU1}, \eqref{theta1pq} and \eqref{A1DU2} as before. 

 \begin{remark}\label{Agamma0}
We observe a clear similarity between~\eqref{leadeqc}, \eqref{seceqc}, \eqref{thirdeqc} and their repeated-eigenvalue  predecessors~\eqref{solconmat1}, \eqref{solconmat2} and \eqref{solconmat3} respectively. In fact, the differences are in this case confined to the diagonal matrix $\boldsymbol{\Gamma}\boldsymbol{D} = \boldsymbol{A}^{\gamma} - \boldsymbol{A}^{\tz}$ that accounts for separations between the neighboring eigenvalues according to~\eqref{eig_scal}. Further, since~$\boldsymbol{A}^{\tz}$ is Hermitian, so is~$\boldsymbol{A}^{\gamma}$.
\end{remark}
%--------------------------------------------------------------------------------------------------------------%
\subsection{Eigenfunction basis}\label{eig_basis}
%--------------------------------------------------------------------------------------------------------------%
\noindent Let $\boldsymbol{P}$ denote the matrix of eigenvectors associated with the generalized eigenvalue problem
\[
\boldsymbol{A}^{\gamma}(\hat\bk) \hh \boldsymbol{v} \;=\; \tau \boldsymbol{D}\boldsymbol{v}. 
\]
In order to diagonalize~$\boldsymbol{A}^{\gamma}$, we express $\bu_0,\bu_1,\bu_2, \mathbf{f}_0$ and~$\mathbf{f}_1$ in terms of~$\boldsymbol{P}$ as
\begin{equation}\label{change_basis}
\bu_0 = \boldsymbol{P}\hh  \bu'_0, \quad \bu_1 = \boldsymbol{P}\hh  \bu'_1, \quad \bu_2 = \boldsymbol{P}\hh  \bu'_2,\quad \mathbf{f}_0 = \boldsymbol{P}\hh  \mathbf{f}'_0, ,\quad \mathbf{f}_1 = \boldsymbol{P}\hh  \mathbf{f}'_1
\end{equation}
and we premultiply \eqref{leadeqc}--\eqref{thirdeqc} by $\overline{\boldsymbol{P}}^{T}$. For simplicity, we then drop the prime symbol from the ``rotated'' vectors~$\bu_0',\bu_1',\bu_2', \mathbf{f}_0'$ and~$\mathbf{f}_1'$, and we keep the original notation of the transformed matrices in~\eqref{leadeqc}--\eqref{thirdeqc}. In this setting, we have 
\begin{equation} \label{Adiaggg}
\boldsymbol{A}^{\gamma}(\hat\bk) \;=\; \text{diag}(\tau_1,\tau_2,\dots, \tau_{Q}), 
\end{equation}
noting for future reference that~$\tau_q=0$ ($q=\overline{1,N_0}$) when $\text{rank}(\boldsymbol{A}^{\gamma})=Q-N_0$.

%--------------------------------------------------------------------------------------------------------------%
\subsection{Leading-order approximation}\label{eig_basis}
%--------------------------------------------------------------------------------------------------------------%

\noindent Thanks to the presence of the ``penalty'' term~$\boldsymbol{\Gamma}\boldsymbol{D}$ in~\eqref{Agam}, $\boldsymbol{A}^{\gamma}$ is of at least partial rank when~$\boldsymbol{A}^{\tz}\!=\boldsymbol{0}$. As a result, in the sequel we present the effective models for full- and partial-rank~$\boldsymbol{A}^{\gamma}$ only.

%--------------------------------------------------------------------------------------------------------------%
\subsubsection{Effective solution for full-rank $\boldsymbol{A}^{\gamma}$ }\label{full_rankc}
%--------------------------------------------------------------------------------------------------------------%

\noindent When~$\text{rank}(\boldsymbol{A}^{\gamma})=Q$, $\tilde{f}_\bk\neq 0$, and $\omega^2\!-\omega_n^2=O(\eps)$ i.e.~$|\hat\sigma|\!=\!1$, we must have $\bu_0=\boldsymbol{0}$ due to~\eqref{leadeqc}. As a result, \eqref{seceqc} yields the leading-order effective equation
\begin{equation}\label{leadsolfuran}
- (\boldsymbol{A}^{\gamma} +\chsig \check{\omega}^2 \boldsymbol{D})\hh \bu_1 \;=\; \mathbf{f}_0.
\end{equation}
When~$\tilde{f}_\bk=0$, eigenvalues of the GEP stemming from~\eqref{leadeqc} (or equivalently~\eqref{leadsolfuran}) define the leading-order asymptotic approximation of the dispersion relationships in direction $\hat{\bk}/\|\hat\bk\|$ as
\begin{equation}\label{leadasympdb}
\omega_{n_{\nes q}}^2 \;=\; \omega_{n_{0}}^2 - \eps \frac{\tau_q}{\rho^{\tz}_{q}}, 
\end{equation}
where~$\tau_q$ denotes the~$q$th eigenvalue of~$\boldsymbol{A}^{\gamma}(\hat\bk)$. 

\begin{remark}
When $\boldsymbol{A}^{\gamma}$ is of full rank, \eqref{leadsolfuran} and~\eqref{leadasympdb} provide a general framework to handle the clusters of nearby dispersion branches, regardless of the fact whether they intersect or ``almost touch'' for example at $\bk=\bk_s$.
\end{remark}
%--------------------------------------------------------------------------------------------------------------%
\subsubsection{Effective solution for partial-rank $\boldsymbol{A}^{\gamma}$ }\label{par_rankc}
%--------------------------------------------------------------------------------------------------------------%

\noindent When~$\text{rank}(\boldsymbol{A}^{\gamma})=Q-N_0$ for some~$N_0\!>\!0$ and  $\tilde{f}_\bk\neq 0$, we first consider the situation where $\omega^2-\lambda_{n_{0}}=O(\eps)$ i.e.~$|\check\sigma|=1$. In this case the leading-order effective equation is again given by~\eqref{leadsolfuran}, while the last $Q-N_0$ dispersion branches are approximated by~\eqref{leadasympdb} for $q\in\overline{N_0\!+\!1,Q}$. 

On the other hand, when $\omega^2-\omega{^2_{n_{0}}}=O(\eps^2)$ i.e.~$|\hat\sigma|=1$, the leading-order effective model~$\bu_0$ is given by 
\begin{eqnarray}\label{secasympdb}
-\sum_{q=1}^{N_0}(B^{\tz}_{pq}(\hat\bk) +\hasig \hat{\omega}^2 D_{pq})\hh u_{0q} &\!\!\!=\!\!\!& \text{f}_{0p}, \quad p\in\overline{1,N_0}, \label{leadsolparran}  \\*[-1mm]
u_{0p} &\!\!\!=\!\!\!& ~\:0, \quad p\in\overline{N_0+1,Q}. \label{secasympdbx}
\end{eqnarray}
When~$\tilde{f}_\bk=0$, the leading-order approximation of the first~$N_0$ dispersion branches is obtained  by solving the GEP affiliated with~\eqref{secasympdb}.

%--------------------------------------------------------------------------------------------------------------%
\subsection{First-order correctors}\label{eig_basis}
%--------------------------------------------------------------------------------------------------------------%
 %--------------------------------------------------------------------------------------------------------------%
\subsubsection{Effective solution for full-rank $\boldsymbol{A}^{\gamma}$ }\label{Afullrankx}
%--------------------------------------------------------------------------------------------------------------%
\noindent When $\boldsymbol{A}^{\gamma}$ is of full rank, $\tilde{f}_\bk\neq 0$, and $\omega^2\!-\omega_n^2=O(\eps)$ i.e.~$|\check\sigma|=1$, one can show that the first-order corrector~$\bu_2$ solves 
\begin{equation}\label{secordcor12}  
-\boldsymbol{B}^{\tz}(\hat\bk)\hh \bu_1 \hh-\hh (\boldsymbol{A}^{\gamma}(\hat\bk) +\chsig \check{\omega}^2 \boldsymbol{D}) \bu_2 \;=\; \mathbf{f}_1,
\end{equation}
with~$\bu_1$ being given by~\eqref{leadsolfuran}. Thanks to~\eqref{solvec1}, $\eqref{leadsolfuran}+\eps \eqref{secordcor12}$ yields the first-order effective model 
\begin{equation}\label{bojan3}
- (\boldsymbol{B}^{\tz}(\eps\hat\bk) + \boldsymbol{A}^{\gamma}(\eps\hat\bk) + \eps\chsig\check{\omega}^2 \boldsymbol{D})\hh \langle\tilde{\bu}\rangle_{\!\rho} \;\stackrel{\eps^2}{=}\; \mathbf{f}_0 + \eps\hh\mathbf{f}_1.
\end{equation}

%--------------------------------------------------------------------------------------------------------------%
\subsubsection{Effective solution for partial rank $\boldsymbol{A}^{\gamma}$ }\label{Apartialrankx}
%--------------------------------------------------------------------------------------------------------------%

\noindent When $\text{rank}(\boldsymbol{A}^{\gamma})<Q$, $\tilde{f}_\bk\neq 0$ and $\omega^2-\lambda_{n_0}=O(\eps)$, the first-order corrector~$\bu_2$ is given by~\eqref{secordcor12}. In contrast, when $\omega^2-\lambda_{n_0}=O(\eps^2)$, the first-order corrector is given by~$\bu_1$ whose components can be shown to satisfy 
\begin{eqnarray}\label{secorccor3}
-\sum_{q=1}^{N}(B^{\tz}_{pq}(\hat\bk) +\hasig \hat{\omega}^2 D_{pq})\hh u_{1q}~= &\text{f}_{1p}+ \sum_{q=1}^{N_0} A^{\so}_{pq}(\hat\bk) \hh u_{0q}, & p\in\overline{1,N_0},  \\
-A^{\gamma}_{pp}(\hat\bk) \hh u_{1p}~= & \text{f}_{0p}+ \sum_{q=1}^{N_0}B^{\tz}_{pq}(\hat\bk) \hh u_{0q}, & p\in\overline{N_0\!+1,Q}, \label{secorccor3x}
\end{eqnarray}
with~$u_{0q}$ ($q=\overline{N_0\!+1,Q}$) being subject to~\eqref{leadsolparran}. Accordingly, we obtain the first-order model as~$\langle\tilde{\bu}\rangle_{\!\rho} \stackrel{1}{=} \eps^{-2}\bu_0+\eps^{-1}\bu_1$, with the scaled summands being directly computable from~\eqref{secasympdb}--\eqref{secasympdbx} and~\eqref{secorccor3}--\eqref{secorccor3x} on replacing~$\hat\bk$ and~$\hat\sigma\hat\omega$ respectively by~$\eps\hh\hat\bk$ and~$ \eps^2\hat\sigma\hat\omega$.

%--------------------------------------------------------------------------------------------------------------%
\section{Discussion} \label{Num_examp}
%--------------------------------------------------------------------------------------------------------------%

\noindent In this section, we share new insights stemming from the general analysis, and we discuss several special cases in support of the numerical simulations (Section~\ref{Numerical examples}). 

\begin{remark}
A common thread of our developments is that we approximate the Bloch wave in terms of its projection to the nearest~$Q$ branches, $Q\geqslant 1$. In this vein, term~$\tilde{u}(\bx)$ on the left-hand side of ansatz~\eqref{sol_scal} and its descendants such as~\eqref{utildeexpc} should be interpreted in the sense of restriction of~\eqref{uts} to the nearest~$Q$ dispersion branches, namely 
\begin{equation}\label{utsQ}
-\sum_{n=1}^{\infty} \frac{(\tilde{f}_{\bk},\tilde{\phi}_n)~\tilde{\phi}_n(\bx)}{ (\rho \tilde{\phi}_n,\tilde{\phi}_n) (\omega^2-\tilde{\lambda}_n) } \quad \Longrightarrow \quad
-\sum_{q=1}^{Q} \frac{(\tilde{f}_{\bk},\tilde{\phi}_{n_q})~\tilde{\phi}_{n_q}(\bx)} {(\rho \tilde{\phi}_{n_q},\tilde{\phi}_{n_q}) (\omega^2-\tilde{\lambda}_{n_q})}.
\end{equation}
This leaves an open question regarding the contribution of ``remote'' branches ($n\neq n_q, q=\overline{1,Q}$) that is beyond the scope of this study, see for instance the recent discussion in~\cite{MOG20}. 
\end{remark}
%--------------------------------------------------------------------------------------------------------------%
\subsection{Energy considerations}\label{synth_1}
%--------------------------------------------------------------------------------------------------------------%

\noindent On the basis of the results in Section~\ref{cluster_eig} which covers the instances of simple and repeated eigenvalues as degenerate cases, we find from~\eqref{utildeexpc} that the instantaneous power density $(\tilde{f}_{\bk},i\omega\tilde{u})=-i\omega (\tilde{f}_{\bk},\tilde{u})$ generated by the source term~$\tilde{f}_\bk$ can be approximated as 
\begin{eqnarray}
\text{Leading order:} &  -i\omega(\tilde{f}_{\bk},\tilde{u})  &   \stackrel{\eps^{-1}}{=}~  -i\omega\sum_q \langle \tilde{f}_{\bk} \rangle^{n_{\nes q}}  \overline{\langle \tilde{u} \rangle^{n_{\nes q}}_{\!\rho}}, \notag \\
\text{First order:} & -i\omega (\tilde{f}_{\bk},\tilde{u})  &  \stackrel{1}{=}~ -i\omega\sum_q \big( \langle \tilde{f}_{\bk} \rangle^{n_{\nes q}} - (\tilde{f}_{\bk},\bchi^{\so}_q)\cdot(i\eps\hat{\bk}) \big) \overline{\langle \tilde{u} \rangle^{n_{\nes q}}_{\!\rho}}, \notag \\*[-2mm]
\text{Second order:} & -i\omega (\tilde{f}_{\bk},\tilde{u})  &  \stackrel{\eps}{=}~ -i\omega\sum_q \big( \langle \tilde{f}_{\bk} \rangle^{n_{\nes q}} - (\tilde{f}_{\bk},\bchi^{\so}_q)\cdot(i\eps\hat{\bk})+ (\tilde{f}_{\bk},\bchi^{\st}_q):(i\eps\hat{\bk})^2 \big) \overline{\langle \tilde{u} \rangle^{n_{\nes q}}_{\!\rho}}  -i\omega\sum_{q=1}^{\bar{Q}}  (\tilde{f}_{\bk},\eta_q^{\tz}). \notag
\end{eqnarray}
The above result in particular demonstrates that the instantaneous power density and threfore the work, averaged in space over~$Y$, equal -- up to the first order -- that exerted by the effective source term on the averaged displacement. This result is in line with the well known Hill-Mandel condition~\cite{Hill63}, requiring that the volume average of the increment of work performed on the representative volume element be equal to the increment of local work performed by the macroscopic i.e. averaged quantities.

%--------------------------------------------------------------------------------------------------------------%
\subsection{Asymptotic solution in physical space near the edge of a band gap}\label{synth_1}
%--------------------------------------------------------------------------------------------------------------%

\noindent With reference to the class~\eqref{forceexp} of source distributions, one immediate application of the foregoing analysis is the case where: (i) $\bk_s = \frac{1}{2}(\sum_j n_j\hh\be^j)$, $n_j\in\{-1,0,1\}$; (ii) the driving frequency is within a band gap near simple eigenfrequency $\omega_n(\bk_s)$, and (iii) the source function $\tilde{f}_\bk$ is given by
\begin{eqnarray}
\tilde{f}_{\bk}(\bx) &\!\!\!=\!\!\!& F(\bk)\phi(\bx), \qquad \phi\in L^2_p(Y),  \quad \text{supp}(F)=\mathcal{C}\subset\mathcal{B}. 
\end{eqnarray}
On recalling BWE~\eqref{BWE} and ansatz~\eqref{sol_scal}, we conveniently introduce the~$p$th-order asymptotic solution in the physical space as 
\begin{eqnarray}
u^{\bpp}(\bx) \;:=\; \sum_{m=0}^p \eps^{m-2}\hh u_m(\boldsymbol{x}), \quad \bx\in S,  \qquad \text{where}\quad 
u_m(\boldsymbol{x}) \;:=\; |\mathcal{C}|^{-1}\int_{\mathcal{C}} \tilde{u}_m(\bx) \, e^{i(\bk_s\!+\eps\hat\bk)\cdot\bx}  \dd(\eps\hat\bk). \label{leadt} 
\end{eqnarray}
From~\eqref{ut0exp}, \eqref{ut1exp}, \eqref{solvacon0}, \eqref{ut2exp}, \eqref{solvacon1}, \eqref{solvacon2}, Remark~\ref{eta0_lin} and Claim~\ref{eta1_id}, we specifically find that
\begin{eqnarray}
\tilde{u}_0(\boldsymbol{x}) &\!\!\!=\!\!\!& -\eps^{2} \,  \frac{\tilde{\phi}_n(\bx) \,  \langle \tilde{f}_{\bk} \rangle}{\rho^{\tz}(\omega^2\!-\omega_n^2) + \bmu^{\tz}\!:\!(i\eps\hat\bk)^2}, \label{leadt0} \\
\tilde{u}_1(\boldsymbol{x}) &\!\!\!=\!\!\!& \eps \, \frac{\big[\tilde{\phi}_n(\bx)\hh (\tilde{f}_{\bk},\bchi^{\so})  -  \langle\tilde{f}_{\bk}\rangle~\bchi^{\so}(\bx)\big] \sip (i\eps\hat\bk) } {\rho^{\tz}(\omega^2\!-\omega_n^2) + \bmu^{\tz}\!:\!(i\eps\hat\bk)^2}, \label{leadt1} \\
\tilde{u}_2(\boldsymbol{x})  &\!\!\!=\!\!\!& -\,  \frac{\big[\tilde{\phi}_n(\bx)\hh \big((\tilde{f}_{\bk},\bchi^{\st}) + \frac{\langle\tilde{f}_{\bk}\rangle}{\rho^{\tz}}\{(\rho\bchi^{\so}\!\otimes \overline{\bchi^{\so}},1)\} \big) - \{(\tilde{f}_{\bk},\bchi^{\so})\!\otimes\bchi^{\so}(\bx)\} + \langle\tilde{f}_{\bk}\rangle \hh\bchi^{\st}(\bx) \big] \!:\! (i\eps\hat\bk)^2} {\rho^{\tz}(\omega^2\!-\omega_n^2) + \bmu^{\tz}\!:\!(i\eps\hat\bk)^2}  \notag \\
  &  & - \,  \frac{\eps^{-2} \tilde{u}_0(\boldsymbol{x}) \; \bmu^{\st}\!:(i\eps\hat\bk)^4} {\rho^{\tz}(\omega^2\!-\omega_n^2) + \bmu^{\tz}\!:\!(i \eps\hat\bk)^2 } \,+\, \zeta^{\tz}(\bx) \, F(\bk_s\!+ \eps\hat\bk). \label{leadt2} 
\end{eqnarray}

In terms of the \emph{effective solution}, by~\eqref{effsoluds} we can similarly introduce the $p$th-order mean motion in the physical space as
\begin{eqnarray}
\langle u\rangle_{\!\rho}^{\bpp}(\bx) \;:=\; \sum_{m=0}^p \eps^{m-2}\hh |\mathcal{C}|^{-1}\int_{\mathcal{C}} u_m(\eps\hat\bk) \, e^{i(\bk_s\!+\eps\hat\bk)\cdot\bx}  \dd(\eps\hat\bk),  \quad \bx\in \mathbb{R}^d, \label{leade}
\end{eqnarray}

via superposition of the averaged Bloch-wave solutions, $u_m(\eps\hat\bk)=\langle \tilde{u}_m\rangle_{\!\rho}$. From~\eqref{leadt0}--\eqref{leadt2}, we clearly have
\begin{eqnarray}
u_0(\eps\hat\bk) &\!\!\!=\!\!\!& -\eps^{2} \,  \frac{\langle \tilde{f}_{\bk} \rangle}{\rho^{\tz}(\omega^2\!-\omega_n^2) + \bmu^{\tz}\!:\!(i\eps\hat\bk)^2}, \label{leade0} \\
u_1(\eps\hat\bk) &\!\!\!=\!\!\!& \eps \, \frac{(\tilde{f}_{\bk},\bchi^{\so}) \sip (i\eps\hat\bk) } {\rho^{\tz}(\omega^2\!-\omega_n^2) + \bmu^{\tz}\!:\!(i\eps\hat\bk)^2}, \label{leade1} \\
u_2(\eps\hat\bk)  &\!\!\!=\!\!\!& -\,  \frac{\big((\tilde{f}_{\bk},\bchi^{\st}) + \frac{\langle\tilde{f}_{\bk}\rangle}{\rho^{\tz}}\{(\rho\bchi^{\so}\!\otimes \overline{\bchi^{\so}},1)\} \big) \!:\! (i\eps\hat\bk)^2 + \eps^{-2} u_0(\eps\hat\bk) \; \bmu^{\st}\!:\!(i\eps\hat\bk)^4} {\rho^{\tz}(\omega^2\!-\omega_n^2) + \bmu^{\tz}\!:\!(i\eps\hat\bk)^2}. \label{leade2} 
\end{eqnarray}
In Section~\ref{Numerical examples}, we shall make use of the above results toward approximating the ``full'' and mean wave motion near the edge of a band gap. 
%--------------------------------------------------------------------------------------------------------------%
\subsection{Dirac behavior in~$\mathbb{R}^2$ for~$Q=2$}\label{Dirac cones}
%--------------------------------------------------------------------------------------------------------------%

\noindent Consider a two-dimensional periodic medium, $S\!\subset\!\mathbb{R}^2$, whose spectral neighborhood~\eqref{scal} features two nearby eigenfrequencies~$\omega_{n_1}=\omega_{n_1}(\bk_s)$ and~$\omega_{n_2}=\omega_{n_2}(\bk_s)$ ($Q=2$). In this case, matrix $\boldsymbol{A}^{\gamma}$ reads 
\begin{equation}\label{Agamma2}
\boldsymbol{A}^{\gamma} = 
\begin{pmatrix}
\bth^{\tz}_{11}\cdot i \hat{\bk}  &\bth^{\tz}_{12}\cdot i \hat{\bk} \\
-\overline{\bth^{\tz}_{12}}\cdot i \hat{\bk} & \bth^{\tz}_{22}\cdot i \hat{\bk}+\gamma\rho^{\tz}_2
\end{pmatrix}, \qquad \gamma= \eps^{-1}(\omega_{n1}^2-\omega_{n2}^2).
\end{equation}
By way of~\eqref{leadeqc}, the two dispersion relationships are accordingly given by
\begin{equation}\label{eig_3}
\omega_{n_{1/2}}^2(\bk) \;=\; \omega_{n_1}^2 - \frac{\eps}{2} \Big( \gamma+\frac{\bth^{\tz}_{11}\cdot i \hat{\bk}}{\rho^{\tz}_1} +\frac{\bth^{\tz}_{22}\cdot i \hat{\bk}}{\rho^{\tz}_2} \pm \sqrt{\Big( \gamma - \frac{\bth^{\tz}_{11}\cdot i \hat{\bk}}{\rho^{\tz}_1} +\frac{\bth^{\tz}_{22}\cdot i \hat{\bk}}{\rho^{\tz}_2} \Big)^2 + \frac{4\{\bth^{\tz}_{12}\!\otimes\overline{\bth^{\tz}_{12}}\}: (\hat{\bk})^2}{\rho^{\tz}_1\rho^{\tz}_2}}\Big),
\end{equation}
where the matrix $\{\bth^{\tz}_{12}\!\otimes \overline{\bth^{\tz}_{12}}\}\in\mathbb{R}^{2\times2}$ is in general positive semi-definite, and specifically positive definite when 
\begin{equation}\label{condef}
\bth^{\tz}_{12}\!\cdot\bi_1\neq0, \quad \bth^{\tz}_{12}\!\cdot\bi_2\neq0, \quad \text{and} \quad \arg(\bth^{\tz}_{12}\!\cdot\bi_1)-\arg(\bth^{\tz}_{12}\!\cdot\bi_2) \neq n \pi, \quad n\in\mathbb{Z}. 
\end{equation}
Equations \eqref{eig_3} describe ``almost touching'' (resp.~crossing) branches when $\gamma\neq0$ (resp.~$\gamma=0$) featuring the middle  plane
\begin{equation}\label{planeD}
 (\mathcal{P}): ~~ \omega_{\mathcal{P}}^2(\bk) \;=\; \omega_{n_1}^2 - \frac{\eps}{2} \Big( \gamma+\frac{\bth^{\tz}_{11}\!\cdot i \hat{\bk}}{\rho^{\tz}_1} +\frac{\bth^{\tz}_{22}\!\cdot i \hat{\bk}}{\rho^{\tz}_2} \Big).
\end{equation}
When $(\mathcal{P})$ is horizontal, we further have 
\begin{equation}\label{horisym}
\frac{1}{\rho^{\tz}_1} \bth^{\tz}_{11}+\frac{1}{\rho^{\tz}_2} \bth^{\tz}_{22} \;=\; \boldsymbol{0},
\end{equation}
which holds true for any $\rho$-orthogonal eigenfunction basis, by the conservation of the trace of $\boldsymbol{A}^{\gamma}$. In this case, dispersion relationship~\eqref{eig_3} simplifies to
\begin{equation}\label{eig_3s}
\omega_{n_{1/2}}^2(\bk) = \omega_{n_1}^2 - \frac{\eps}{2} \Big(\gamma \pm \sqrt{ \Big( \gamma -\frac{2\bth^{\tz}_{11}\!\cdot i \hat{\bk}}{\rho^{\tz}_1}\Big)^2+ \frac{4\{\bth^{\tz}_{12}\!\otimes \overline{\bth^{\tz}_{12}}\}\!: \! (\hat{\bk})^2}{\rho^{\tz}_1\rho^{\tz}_2}}\Big).
\end{equation}

%--------------------------------------------------------------------------------------------------------------%
\subsubsection{Dirac cones}\label{DC}
%--------------------------------------------------------------------------------------------------------------%

\begin{figure}[h!] 
\centering{\includegraphics[width= 140 mm]{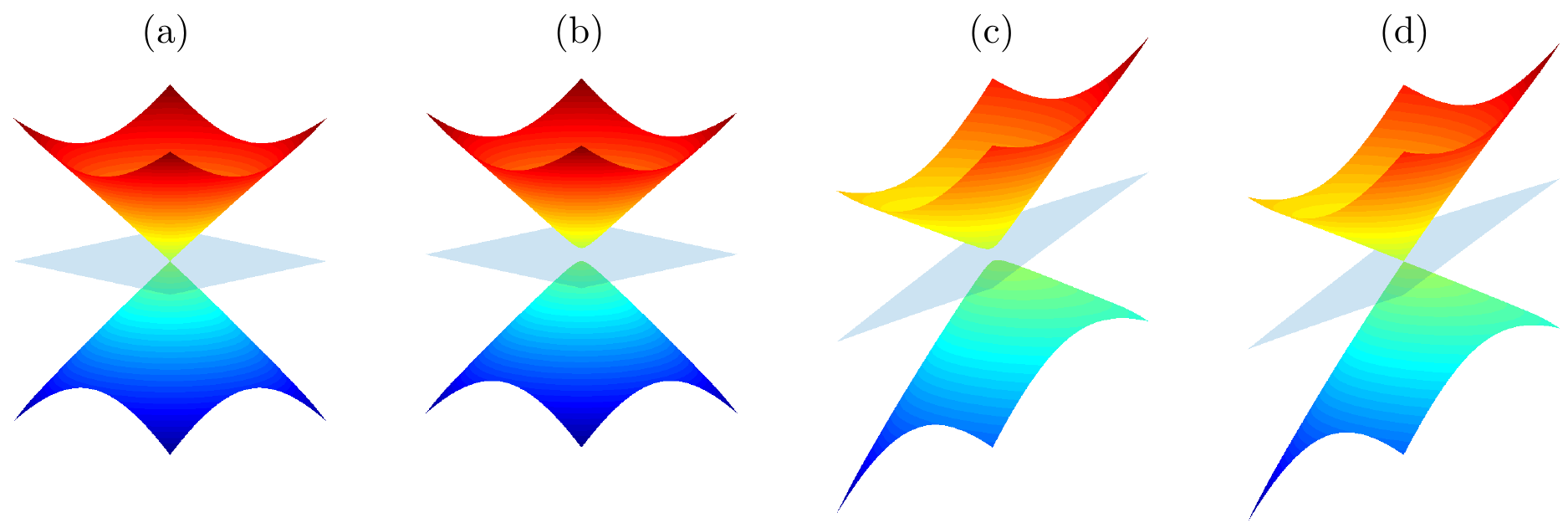}}\vspace*{-2mm}
\caption{Examples of: (a) Dirac cones, (c) Blunted Dirac cones, (c) tilted-blunted Dirac cones, and (d) tilted Dirac cones. All panels include the ``middle plane'' $(\mathcal{P})$ according to~\eqref{planeD}. } \label{fig3} \vspace*{-3mm}
\end{figure}

\noindent With reference to \eqref{eig_3s}, when $\gamma=0$ i.e.~$\omega_{n_1}=\omega_{n_2}$, we choose the eigenfunction basis that diagonalizes $\boldsymbol{A}^{\gamma}$ in given direction $\hat{\bk}/\|\hat{\bk}\|$. In this case all vectors $\bth^{\tz}_{pq}$ are imaginary-valued up to a complex multiplier, since $\bth^{\tz}_{11}, \bth^{\tz}_{22}\in i\mathbb{R}^2$ by~\eqref{teta&rhome} and $\bth^{\tz}_{12}\cdot i \hat{\bk}=\boldsymbol{0}$.  As a result, the dispersion relationships \eqref{eig_3s} are characterized by (i) elliptical isocontours when the vectors $(\rho^{\tz}_1)^{-1}\bth^{\tz}_{11}$ and~$(\rho^{\tz}_1\rho^{\tz}_2)^{-1/2}\bth^{\tz}_{12}$ are linearly independent, and (ii) circular isocontours when they are orthogonal with equal norms. In the latter case, \eqref{eig_3s} reduces to  
\begin{equation}\label{eig_3sss}
\omega_{n_{1/2}}^2(\bk) = \omega_{n_1}^2 \mp \frac{\|\bth^{\tz}_{11}\|}{\rho^{\tz}_1}\|\eps\hat{\bk}\|,
\end{equation}
which describe axisymmetric Dirac cones as depicted in Fig.~\ref{fig3}(a). 

%--------------------------------------------------------------------------------------------------------------%
\subsubsection{Blunted Dirac cones}\label{BDC}
%--------------------------------------------------------------------------------------------------------------%

\noindent When~$\gamma\neq0$ in~\eqref{eig_3s}, assuming $\text{d}\omega_{n_{1/2}}/\text{d}\bk|_{\bk=\bk_s}=\boldsymbol{0}$ implies that~$\bth^{\tz}_{11}=\boldsymbol{0}$ and thus~$\bth^{\tz}_{22}=\boldsymbol{0}$ by~\eqref{horisym}. Further, if the matrix $\{\bth^{\tz}_{12}\!\otimes \overline{\bth^{\tz}_{12}}\}$ is positive-definite due to~\eqref{condef}, the dispersion relationships in~\eqref{eig_3s} are characterized by elliptic isocontours and thus exhibit cone-like geometry. As a special case, the isotropy of $\{\bth^{\tz}_{12}\!\otimes \overline{\bth^{\tz}_{12}}\}$ is attained when~$\bth^{\tz}_{12} \parallel \bi_1\pm i \bi_2$, in which situation~\eqref{eig_3s} becomes
\begin{equation}\label{eig_3ss}
\omega_{n_{1/2}}^2(\bk) = \omega_{n_1}^2 - \frac{\eps}{2} \Big(\gamma \pm \sqrt{\gamma^2 + \frac{2\|\bth^{\tz}_{12}\|^2}{\rho^{\tz}_1\rho^{\tz}_2} \|\hat{\bk}\|^2}\Big).
\end{equation}
Geometrically, \eqref{eig_3ss} describe an axisymmetric variant of the \emph{blunted} Dirac cones shown in Fig.~\ref{fig3}(b). 

\begin{remark}\label{spcaseapex}
At ``apexes'' of the Brillouin zone $\bk_s=\frac{1}{2}(\sum_jn_j\be^j)$, $n_j\in\{-1,0,1\}$, we have $\bth^{\tz}_{11}=\boldsymbol{0}$ and $\bth^{\tz}_{12}\in\mathbb{R}^d$ by Claim~\ref{Aimag}. As a result, we cannot have axisymmetric Dirac cones~\eqref{eig_3sss} nor axisymmetric blunted Dirac cones~\eqref{eig_3ss} there. This claim also applies to situations without axial symmetry.
\end{remark}
%--------------------------------------------------------------------------------------------------------------%
\subsubsection{Tilted- and tilted-blunted Dirac cones}\label{TBDC}
%--------------------------------------------------------------------------------------------------------------%

\noindent Let us next forgo condition~\eqref{horisym} which ensures that $(\mathcal{P})$ is horizontal. When~$\gamma\neq 0$, \eqref{eig_3} describes a pair of \emph{tilted-blunted} Dirac cones only if: (i) the first term under the square-root sign reduces to~$\gamma^2$, namely  
\begin{equation}\label{tilted1}
\frac{1}{\rho^{\tz}_1} \bth^{\tz}_{11} = \frac{1}{\rho^{\tz}_2} \bth^{\tz}_{22}\neq\boldsymbol{0}, 
\end{equation}
and (ii)~$\{\bth^{\tz}_{12}\!\otimes \overline{\bth^{\tz}_{12}}\}$ is positive-definite according to~\eqref{condef}. If further~$\bth^{\tz}_{12} \parallel \bi_1\pm i \bi_2$, the two cones become ``symmetric'' in that~$\omega^2_{1/2}(\bk)-\omega_{\mathcal{P}}^2(\bk)$ is axisymmetric in terms of~$\hat\bk$.

When~$\gamma=0$, and~\eqref{horisym} is violated, on the other hand, \eqref{eig_3} describe \emph{tilted} Dirac cones (see Fig.~\ref{fig3}(c)) provided that the vectors~$(\rho^{\tz}_1)^{-1} \bth^{\tz}_{11}-(\rho^{\tz}_2)^{-1} \bth^{\tz}_{22}$ and~$(\rho^{\tz}_1\rho^{\tz}_2)^{-1/2}\bth^{\tz}_{12}$ are linearly independent. In this case, the ``cone symmetry'' is relaxed and applies as long as the two vectors are mutually orthogonal with the same norm.

%--------------------------------------------------------------------------------------------------------------%
\subsection{{Dirac behavior in~$\mathbb{R}^2$ for~$Q=3$}}\label{Dirac-like cones}
%--------------------------------------------------------------------------------------------------------------%

\noindent Consider a cluster of three eigenfrequencies, $\omega_{n_1}\!=\omega_{n_2}$ and~$\omega_{n_3}$, at $\bk_s=\frac{1}{2}(\sum_j n_j\be^j)$, $n_j\in\{-1,0,1\}$. Thanks to Claim~\ref{Aimag}, we find that
\begin{equation}\label{Agamma3}
\boldsymbol{A}^{\gamma} = 
\begin{pmatrix}
0  &\bth^{\tz}_{12}\cdot i \hat{\bk} &\bth^{\tz}_{13}\cdot i \hat{\bk} \\
-\bth^{\tz}_{12}\cdot i \hat{\bk} & 0 & \bth^{\tz}_{23}\cdot i \hat{\bk} \\
-\bth^{\tz}_{13}\cdot i \hat{\bk} & -\bth^{\tz}_{23}\cdot i \hat{\bk} &\gamma\rho^{\tz}_3
\end{pmatrix}, \qquad \gamma = \eps^{-1}(\omega_{n_1}^2- \omega_{n_3}^2), \quad \bth^{\tz}_{pq}\in\mathbb{R}^2.
\end{equation}

When $\gamma\neq0$ and~$\boldsymbol{A}^{\gamma}$ has partial rank in all directions $\hat{\bk}/\|\hat\bk\|$, condition $\det\boldsymbol{A}^{\gamma}=0$, implies that $\bth^{\tz}_{12}=\boldsymbol{0}$. In this case, we have $\omega_{n1}^2(\hat\bk)=\omega_{n1}^2$ (a horizontal plane), and 
\begin{equation}\label{eig_4}
\omega_{n_{2/3}}^2(\bk) \;=\; \omega_{n_1}^2 - \frac{\eps}{2} \Big(\gamma \pm \sqrt{ \gamma^2+ \frac{4(\bth^{\tz}_{13})^2: (\hat{\bk})^2}{\rho^{\tz}_1\rho^{\tz}_3}+ \frac{4(\bth^{\tz}_{23})^2: (\hat{\bk})^2}{\rho^{\tz}_2\rho^{\tz}_3}} \Big),
\end{equation}
which describe a pair of Dirac cones (with elliptic isocontours) as long as the vectors $(\rho^{\tz}_1\rho^{\tz}_3)^{-\frac{1}{2}}\bth^{\tz}_{13}$ and $(\rho^{\tz}_2\rho^{\tz}_3)^{-\frac{1}{2}}\bth^{\tz}_{23}$ are linearly independent. As before, the axial symmetry of \eqref{eig_4} is attained when the two  vectors are orthogonal and have equal norms.

Assuming~$\gamma=0$, on the other hand, $\boldsymbol{A}^{\gamma}$ becomes anti-symmetric and thus necessarily rank deficient. In this case, the counterpart of~\eqref{eig_4} reads 
\begin{equation}\label{eig_5}
\omega_{n_{2/3}}^2(\bk) = \omega_{n_1}^2 \mp \eps \, \sqrt{ \frac{(\bth^{\tz}_{12})^2: (\hat{\bk})^2}{\rho^{\tz}_1\rho^{\tz}_2}+ \frac{(\bth^{\tz}_{13})^2: (\hat{\bk})^2}{\rho^{\tz}_1\rho^{\tz}_3}+ \frac{(\bth^{\tz}_{23})^2: (\hat{\bk})^2}{\rho^{\tz}_2\rho^{\tz}_3}}, 
\end{equation}
which describe a pair of Dirac cones provided that the sum inside the square root is axially-symmetric in terms of~$\hat{\bk}$. We illustrate this case by letting $\bth^{\tz}_{12}=\boldsymbol{0}$ and assuming that $(\rho^{\tz}_1\rho^{\tz}_3)^{-\frac{1}{2}}\bth^{\tz}_{13}$ and $(\rho^{\tz}_2\rho^{\tz}_3)^{-\frac{1}{2}}\bth^{\tz}_{23}$ are orthogonal with equal norms. In such instance, \eqref{eig_5} reduces to
\begin{equation}\label{eig_5s}
\omega_{n_{2/3}}^2(\bk) = \omega_{n_1}^2 \mp \eps  \frac{\|\bth^{\tz}_{13}\|}{(\rho^{\tz}_1\rho^{\tz}_3)^{1/2}} \|\hat{\bk}\|,
\end{equation}
which yield the respective group velocities as
\begin{equation}\label{gpv}
\boldsymbol{c}_{g_{2/3}}(\bk) =  \mp \frac{\|\bth^{\tz}_{13}\|}{2\omega_{n_1}(\rho^{\tz}_1\rho^{\tz}_3)^{1/2}} \frac{\hat\bk}{\|\hat{\bk}\|}.
\end{equation}

\begin{remark}
Equations \eqref{eig_5s} and \eqref{gpv} describe the behavior of the so-called Zero Index Metamaterials (ZIM) \cite{AF15,HCWPK18,HLHZC11} where Dirac-like dispersion relationship occurs (for some  $\omega_{n_1}\!=\omega_{n_2}\!=\!\omega_{n_3}$) at the origin of the Brillouin zone $\bk_s\!=\!\boldsymbol{0}$. In this neighborhood, the phase velocity of branches~$n_2$ and~$n_3$ approaches infinity, while the affiliated group velocities are non-trivial -- and in fact constant in any given direction~$\hat\bk/\|\hat\bk\|$. This allows for the propagation of energy with near-zero phase delay across finite distances which has applications to e.g.~cloaking, wave tunneling, and directive emission \cite{AF15}.
\end{remark} 

%--------------------------------------------------------------------------------------------------------------%
\subsection{Dirac behavior in~$\mathbb{R}^3$ for~$Q=2$}\label{Dirac cones3D}
%--------------------------------------------------------------------------------------------------------------%

\begin{figure}[h!] 
\centering{\includegraphics[height= 45 mm]{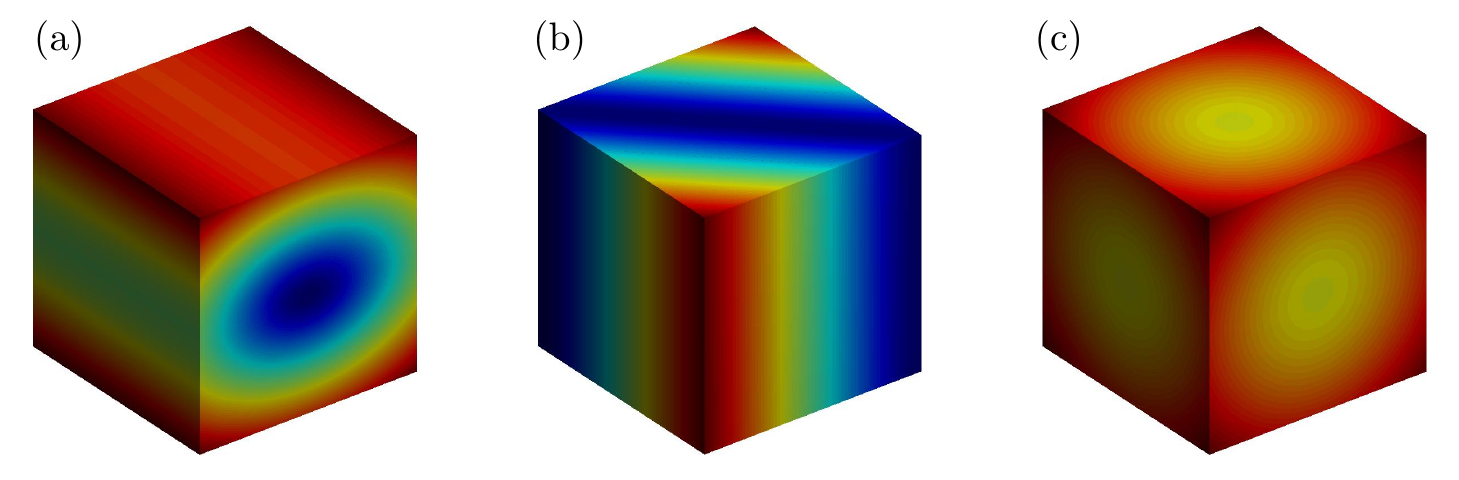}}\vspace*{-2mm}
\caption{Cross-sections of the dispersion relationship $\omega_{n_2}(\bk)$ for~$d=3$ according to: (a)  equation~\eqref{eig_6} where $\text{Re}(\bth^{\tz}_{12})$ and~$\text{Im}(\bth^{\tz}_{12})$ are linearly independent; (b) equation~\eqref{eig_6} where $\text{Re}(\bth^{\tz}_{12})$ and~$\text{Im}(\bth^{\tz}_{12})$ are parallel, and (c) equation~\eqref{eig_3sss} depicting a single hyper-cone.} \label{fig4} \vspace*{-3mm}
\end{figure}

\noindent Consider a periodic medium $S\subset\mathbb{R}^3$ that presents two nearby (or repeated) eigenfrequencies, $\omega_{n_1}$ and~$\omega_{n_2}$.  The two dispersion relationships are in this case also given by~\eqref{eig_3}, where the matrix $\{\bth^{\tz}_{12}\!\otimes \overline{\bth^{\tz}_{12}}\}\in\mathbb{R}^{3\times3}$ is positive semi-definite and necessarily rank deficient. When $\gamma\neq0$ and~$\text{d}\omega_{n_{1/2}}/\text{d}\bk|_{\bk=\bk_s}=\boldsymbol{0}$ , 
we find that $\bth^{\tz}_{11}=\bth^{\tz}_{22}=\boldsymbol{0}$ which reduces~\eqref{eig_3} to 
\begin{equation}\label{eig_6}
\omega_{n_{1/2}}^2(\bk) = \omega_{n_1}^2 - \frac{\eps}{2} \Big(\gamma \pm \sqrt{ \gamma^2 +  \frac{4\{\bth^{\tz}_{12}\!\otimes \overline{\bth^{\tz}_{12}}\}: (\hat{\bk})^2}{\rho^{\tz}_1\rho^{\tz}_2}}\Big).
\end{equation}
Expressions~\eqref{eig_6} describe anisotropic dispersion relationships that are: (a)~invariant in the direction $\text{Re}(\bth^{\tz}_{12})\times\text{Im}(\bth^{\tz}_{12})$ of $\hat{\bk}$ when $\text{Re}(\bth^{\tz}_{12})$ and $\text{Im}(\bth^{\tz}_{12})$ are linearly independent (see Fig.~\ref{fig4}(a)), and (b)~invariant within the planes orthogonal to $\text{Re}(\bth^{\tz}_{12})$ when the latter two vectors are parallel (se Fig.~\ref{fig4}(b)). 

On the other hand, when $\gamma=0$, equations \eqref{eig_3} describe \emph{hyper-cones}  provided that (i)~condition \eqref{horisym} holds, and (ii)~the vectors $(\rho^{\tz}_1)^{-1}\bth^{\tz}_{11}$, $\;(\rho^{\tz}_1\rho^{\tz}_2)^{-\frac{1}{2}}\text{Re}(\bth^{\tz}_{12})$ and~$(\rho^{\tz}_1\rho^{\tz}_2)^{-\frac{1}{2}}\text{Im}(\bth^{\tz}_{12})$ are mutually orthogonal with equal norms. In such case, the dispersion relationships are given by~\eqref{eig_3sss}, a scenario that is illustrated in Fig.~\ref{fig4}(c).

\begin{remark}\label{spcaseapex2}
  With reference to Claim~\ref{Aimag}, the real-valuedness of $\bth^{\tz}_{pq}$ prevents the existence of hyper-conical dispersion relationships at ``apexes'' $\bk_s=\frac{1}{2}(\sum_j n_j\be^j)$, $n_j\in\{-1,0,1\}$ of the first Brillouin zone. 
\end{remark}

%--------------------------------------------------------------------------------------------------------------%
\section{Numerical examples }\label{Numerical examples}
%--------------------------------------------------------------------------------------------------------------%

\noindent In this section, we seek to illustrate the utility of the proposed homogenization framework by considering both free ($\tilde{f}_\bk=0$) and forced ($\tilde{f}_\bk\neq 0$) wave motion problems.  

%--------------------------------------------------------------------------------------------------------------%
\subsection{Dispersion relationships}\label{disper_approx}
%--------------------------------------------------------------------------------------------------------------%

\noindent Let us first examine the performance of the asymptotic solution in terms of local approximation of the dispersion relationships. To this end, we consider a non-orthogonal lattice with Neumann exclusions, and an orthogonal lattice with Dirichlet exclusions. 

%--------------------------------------------------------------------------------------------------------------%
\subsubsection{Kagome lattice}\label{kagome_exp}
\begin{figure}[h!] 
\centering{\includegraphics[width=\textwidth]{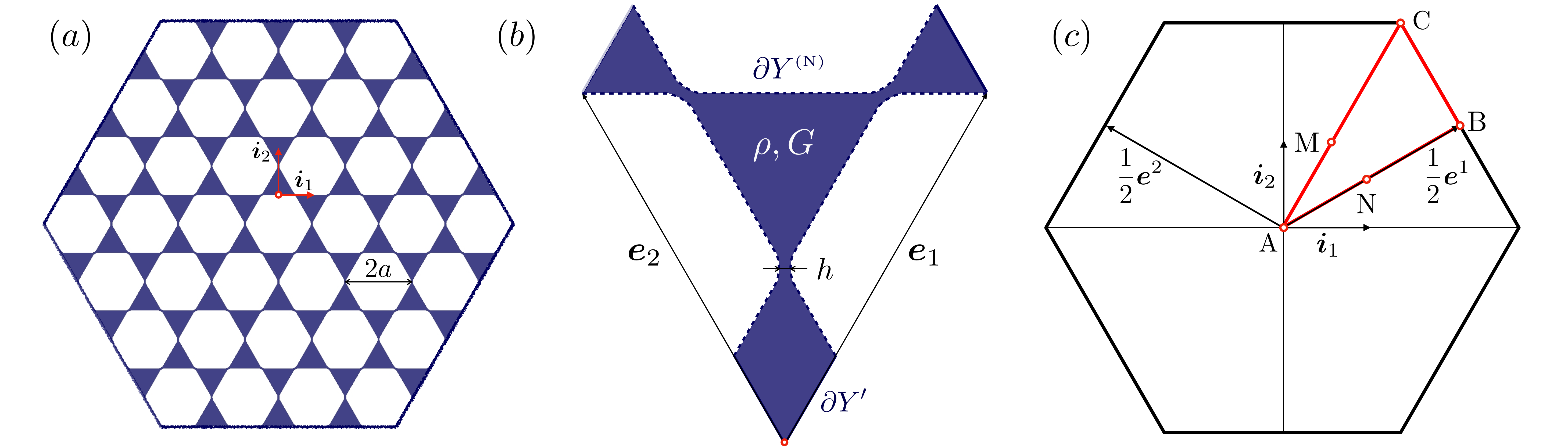}}\vspace*{-2mm}
\caption{Schematics of (a) Kagome lattice $\SSS\subset\mathbb{R}^2$ including the origin of the canonical basis; (b) unit cell of periodicity $Y$, with $\delta Y'$ (solid lines) and $\delta Y^{\tiny{\mbox{(N)}}}$ (dashed lines) indicating respectively the support of periodic and homogenous Neumann boundary conditions; and (c) the first Brillouin zone featuring ``test'' points~A, B, C, M and N.} \label{fig5} \vspace*{0mm}
\end{figure}

%--------------------------------------------------------------------------------------------------------------%
\noindent As a first example, we consider the anti-plane shear wave motion in a Kagome lattice $S\subset\mathbb{R}^2$. This configuration is motivated by a recent experimental study~\cite{MZSMG18} of the wave transport in symmetric and asymmetric Kagome lattices, that revealed frequency-dependent directive behavior in the bulk and the existence of (evanescent) edge modes. With reference to Fig.~\ref{fig5}(a), our lattice is characterized by a trihexagonal tiling geometry where the equilateral triangles of side $a\!=\!1$ are linked by hinges of thickness $h\!=\!0.04a$, yielding the porosity of $v\!=\!0.75$. For completeness, Fig.~\ref{fig5}(b) shows the unit cell of periodicity including the lattice basis vectors $\be_1$ and $\be_2$ are, while Fig.~\ref{fig5}(c) displays the first Brillouin zone including the reciprocal basis vectors $\be^1$ and $\be^2$. The motion in the medium is governed by the wave equation~\eqref{PDE} with $\rho(\bx)=1$ and $G(\bx)=1$, subject to the traction-free boundary condition along the perimeter of hexagonal voids. In this case, the lattice basis vectors and the reciprocal basis vectors are given by
 \begin{equation} \notag
 \be_1 = a \big(\bi_1+\sqrt{3}\hh \bi_2 \big),\quad   \be_2= a \big (-\bi_1+\sqrt{3}\hh \bi_2 \big),\quad  \be^1 = \tfrac{\pi}{a} \big (\bi_1+\tfrac{1}{\sqrt{3}}\hh \bi_2 \big),\quad \be^2 = \tfrac{\pi}{a} \big (-\bi_1+\tfrac{1}{\sqrt{3}} \hh \bi_2 \big).
 \end{equation}
In the absence of the source term ($\tilde{f}_{\bk}=0$), the foregoing homogenization framework enables local approximation the dispersion relationship in the vicinity of an \emph{arbitrary pair} $(\bk_s,\omega_n(\bk_s))$, $\bk_s\in\overline{\mathcal{B}}$, which is a way to access the effective properties of the medium. With reference to  Fig.~\ref{fig5}(c), we illustrate this by taking~$\bk_s$ as the origin of the Brillouin zone (point A), apex points~B and~C, and internal points~M and~N given respectively by    
 \begin{equation}\label{ptsm1} \notag
\overrightarrow{\text{AB}} = \tfrac{1}{2}\hh \be^1,\quad~  \overrightarrow{\text{AC}} = \tfrac{\pi}{a} \big(\tfrac{1}{3}\hh \bi_1+\tfrac{1}{\sqrt{3}}\hh \bi_2\big), \quad~ \overrightarrow{\text{AM}} = 0.4125\, \overrightarrow{\text{AC}}, \quad~~ \overrightarrow{\text{AN}} = 0.4761\, \overrightarrow{\text{AB}}.
 \end{equation}
 
\begin{figure}[h!] 
\centering{\includegraphics[height= 100mm]{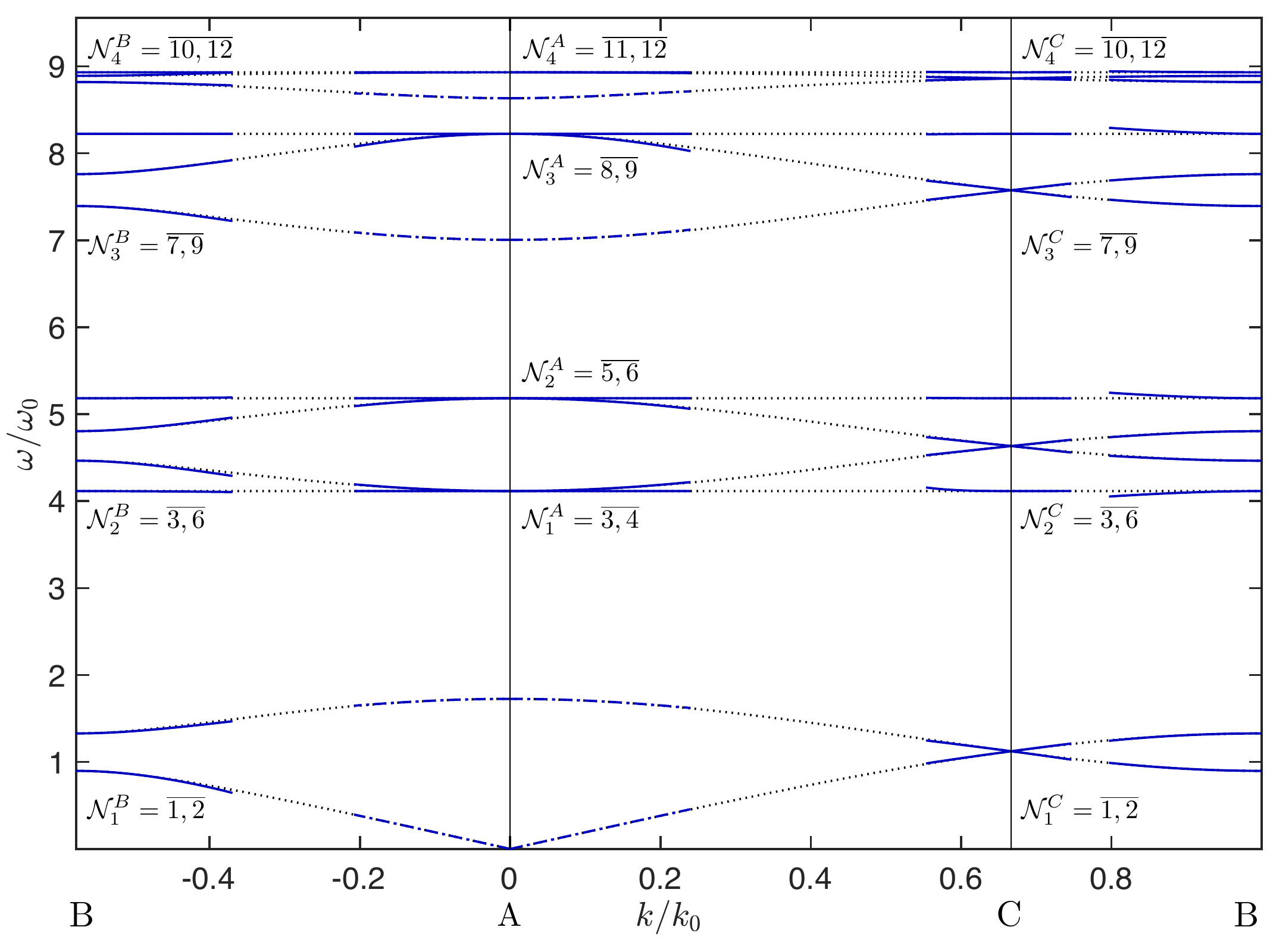}}\vspace*{-3mm}
\caption{Approximation of the first twelve dispersion branches for the Kagome lattice near points A, B and C in Fig~\ref{fig5}(c). In the display, dotted lines track the reference numerical results; solid lines signify the leading-order approximation of the clusters of nearby branches ($Q\!>\!1$), and dash-dotted lines plot the second-order approximation of isolated dispersion branches ($Q\!=\!1$). The normalization parameters are defined as $k_0=\pi/a$ and $\omega_0=\sqrt{G/(\rho a^2)}$.} \label{fig6} \vspace*{0mm}
\end{figure}

\begin{figure}[h!] 
\centering{\includegraphics[height= 100 mm]{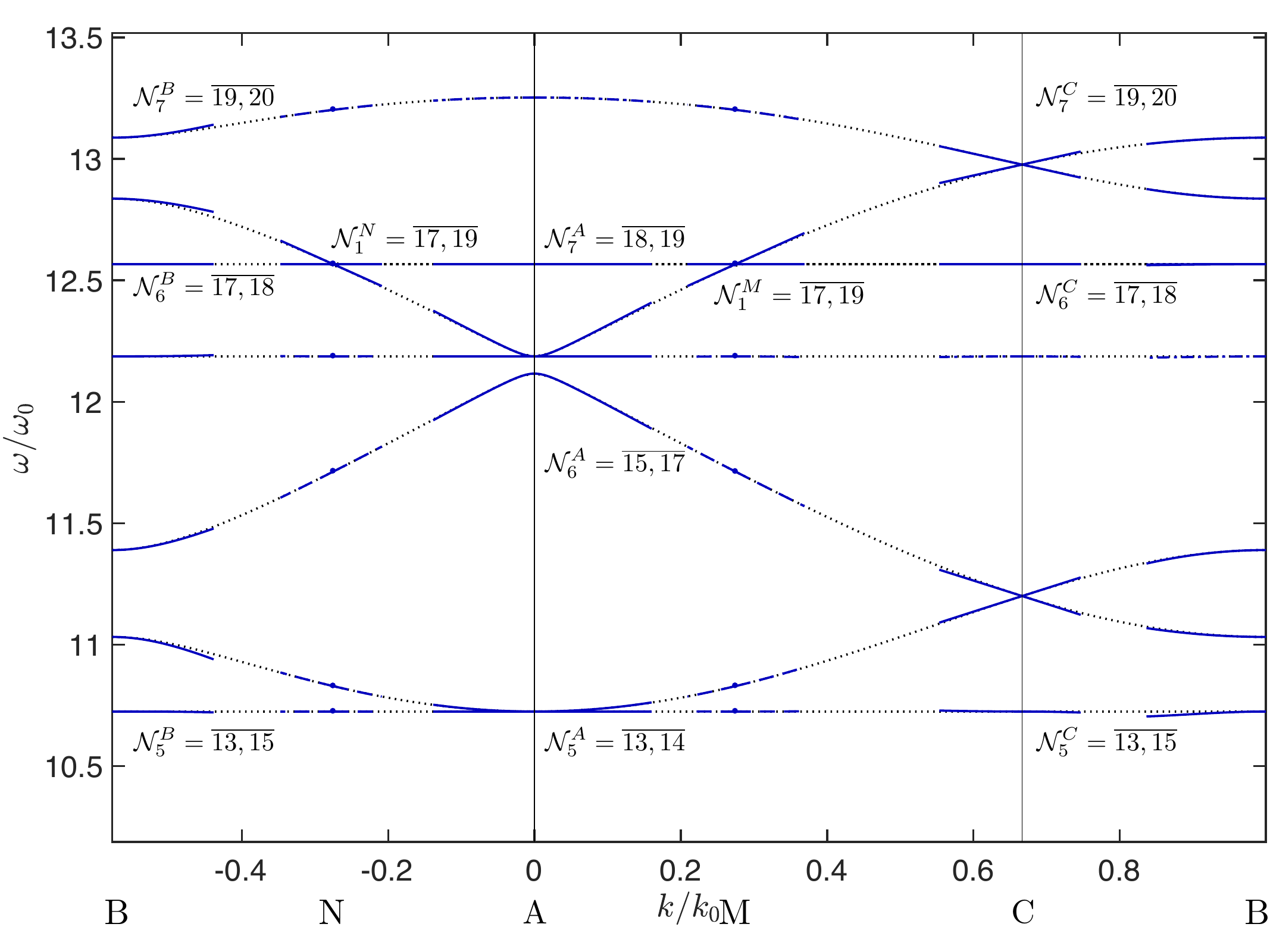}}\vspace*{-3mm}
\caption{Approximation of the dispersion branches 13--20 for the Kagome lattice near points A, B, C, M and N in Fig~\ref{fig5}(c). In the display, dotted lines track the reference numerical results; solid lines signify the leading-order approximation of the clusters of nearby branches ($Q\!>\!1$); dashed lines indicate first-order approximation of isolated dispersion branches ($Q\!=\!1$) at M and N, and dash-dotted lines plot the second-order approximation of isolated dispersion branches at A. The normalization parameters are defined as $k_0=\pi/a$ and $\omega_0=\sqrt{G/(\rho a^2)}$.} \label{fig7} \vspace*{0mm}
\end{figure}

The reference dispersion relationship along path BACB, as well as the cell functions at each wavenumber-eigenfrequency pair $(\bk_s,\omega_n(\bk_s))$ required to evaluate the asymptotic approximation, are computed numerically via the finite element platform NGSolve \cite{NGSolve} by discretizing the unit cell with triangular elements of order 5 and maximum size $h_{max}=0.02a$. Fig.~\ref{fig6} compares the first 12 dispersion branches with their respective approximations in the neighborhood of points A, B and C, while Fig.~\ref{fig7} focuses on branches 13--20 and the neighborhood of points A, B, C, M and N. In each case, we specify the extent of repeated- or cluster-eigenvalue asymptotic approximation (as applicable) by the set
\[
\mathcal{N}_{\ell}^{\star} = \{n_1,n_2,\ldots n_Q\}, \qquad \star\in\{A,B,C,M,N\}
\]
where, for given index~$\star$, $\ell$ locates the cluster in the order of increasing frequency.

In Figs.~\ref{fig6} and~\ref{fig7}, clusters $\mathcal{N}^A_1$ through~$\mathcal{N}^A_5$ and $\mathcal{N}^A_{7}$ each feature a repeated eigenfrequency of multiplicity $Q=2$, where $\boldsymbol{A}(\hat\bk)=\boldsymbol{0}$ in all perturbation directions;  hence the dispersion relationship in those neighborhoods is uniformly described by~\eqref{drasympta2}. In the cluster $\mathcal{N}^A_{6}$ with~$Q=3$, matrix $\boldsymbol{A}^{\gamma}(\hat\bk)$ obtained after expanding about the 15th branch ($n_0\!=\!15$) is found to be of full rank in all perturbation directions; as a result, the local description of the dispersion relationship is in this case provided by~\eqref{leadasympdb}. Alternatively, if the same cluster were expanded about the 16th branch ($n_0\!=\!16$) instead, one would find the local asymptotic description to be given by~$\omega_{n1}^2(\bk)=\omega_{n1}^2$ and~\eqref{eig_4} with~$(n_1,n_2,n_3)=(16,17,15)$.

Clusters $\mathcal{N}^A_{4}$ in Fig.~\ref{fig6} and $\mathcal{N}^A_{7}$, $\mathcal{N}^B_{6}$ and $\mathcal{N}^C_{6}$ in Fig.~\ref{fig7}, on the other hand, each feature two distinct eigenvalues ($Q\!=\!2, \gamma\neq 0$) and trivial effective matrices, $\boldsymbol{A}^{\tz} (\hat\bk)\!=\!\boldsymbol{0}$ and $\boldsymbol{B}^{\tz} (\hat\bk)\!=\!\boldsymbol{0}$, in all perturbation directions. As a result, matrix~$\boldsymbol{A}^{\gamma}(\hat\bk)$ has partial rank, whereby the local dispersion relationships are described by~\eqref{leadasympdb} with $q=2$ and \eqref{secasympdb} with $N_0=1$. Concerning the cluster $\mathcal{N}^B_{1}$ in Fig.~\ref{fig6} where similarly~$Q=2$, $\gamma\!\neq\!0$, it is worth noting that the effective matrix $\boldsymbol{B}^{\tz}$ is in this case \emph{sign-indefinite} as seen from the ``inverted'' curvatures in the directions $BA$ and $BC$.

In Fig.~\ref{fig6} and Fig.~\ref{fig7}, clusters $\mathcal{N}^C_{1}$ and $\mathcal{N}^C_{7}$ feature pure Dirac behavior ($Q\!=\!2$) in directions~CB and~CA, as described by~\eqref{eig_3s} with $\gamma=0$ and $\bth^{\tz}_{11}\nparallel\bth^{\tz}_{12}$. Clusters $\mathcal{N}^C_2$ with~$Q\!=\!4$, $\mathcal{N}^C_3$ with~$Q\!=\!3$, $\mathcal{N}^C_4$  with~$Q\!=\!3$ and $\mathcal{N}^C_5$  with~$Q\!=\!3$ feature similar Dirac behavior; however these clusters also account for the interaction with nearby (non-Dirac) branches for a better description of the local dispersion behavior. From Fig.~\ref{fig7}, we also note that the local approximation of wave dispersion at \emph{internal points} M and N describes with high fidelity the respective numerical results. This holds true for both isolated frequencies and clusters $\mathcal{N}^M_1$ and $\mathcal{N}^N_1$ of size $Q=3$.

%--------------------------------------------------------------------------------------------------------------%
\subsubsection{Pinned square lattice}\label{}
%--------------------------------------------------------------------------------------------------------------%

\begin{figure}[h!] 
\centering{\includegraphics[width=\textwidth]{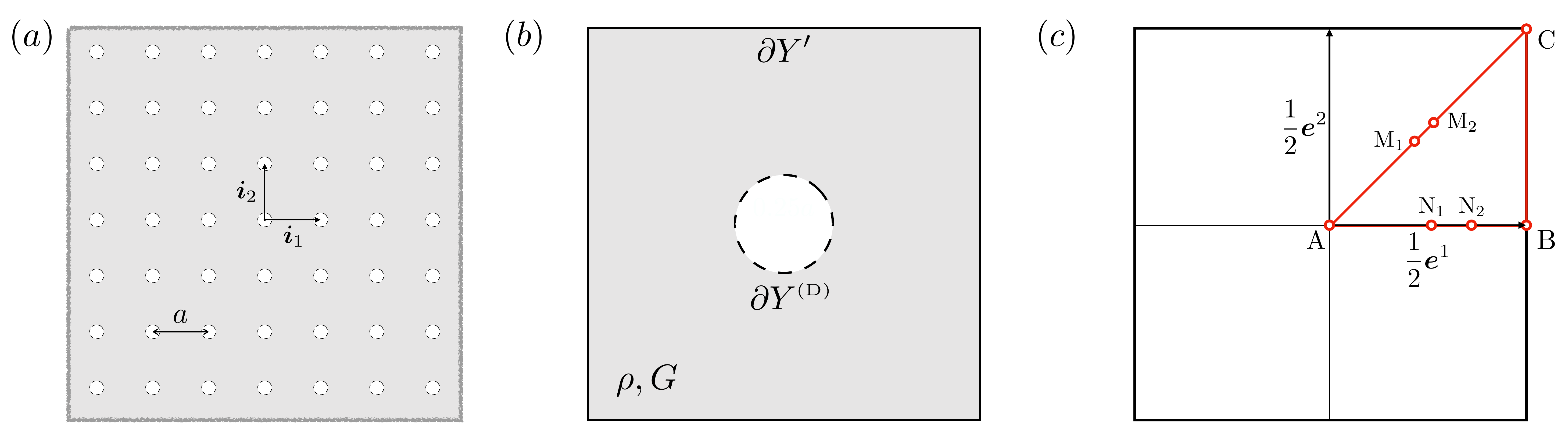}}\vspace*{-2mm}
\caption{Schematics of (a) pinned square lattice $\SSS\!\subset\!\mathbb{R}^2$; (b) unit cell of periodicity $Y$, with $\partial Y'$ (solid lines) and $\partial Y^{\tiny{\mbox{(D)}}}$ (dashed lines) indicating respectively the support of periodic and homogenous Dirichlet boundary conditions; and (c) first Brillouin zone featuring ``test'' points~A, B, C, M$_1$,  M$_2$,  N$_1$, and  N$_2$.} \label{fig8} \vspace*{0mm}
\end{figure}

\noindent As a second example, we take~$\SSS\!\subset\!\mathbb{R}^2$ as a homogenous medium ($G\!=\!1, \,  \rho\!=\!1$) endowed with a square array of circular exclusions, referred to as ``pins'', where~$u(\bx)=0$. Referring to  Fig.~\ref{fig8}(a), the array of pins is characterized by the spacing~$a=1$ and diameter~$0.25a$, resulting in the lattice porosity of $v=0.05$. In this case, the lattice basis vectors and the reciprocal basis vectors are given by  
 \begin{equation} \notag
 \be_1 = a\hh \bi_1,\quad~   \be_2= a\hh \bi_2,\quad~  \be^1 = \tfrac{2\pi}{a}\hh \bi_1,\quad~ \be^2 = \tfrac{2\pi}{a}\hh \bi_2.
 \end{equation}
For completeness, Fig.~\ref{fig8}(b) details the unit cell of periodicity~$Y$, and Fig.~\ref{fig8}(c) illustrates the first Brillouin zone including the ``test'' points~A, B, C, M$_1$,  M$_2$,  N$_1$, and  N$_2$ given by  
 \begin{equation}\label{ptsm2} \notag
\overrightarrow{\text{AB}} = \tfrac{1}{2}\hh \be^1,\quad~  \overrightarrow{\text{AC}} = \tfrac{1}{2}(\be_1+ \be_2\big), \quad~  \overrightarrow{\text{AM$_1$}}= 0.4250 \overrightarrow{\text{AC}} \quad~ \overrightarrow{\text{AM$_2$}} = 0.5250\, \overrightarrow{\text{AC}}, \quad~ \overrightarrow{\text{AN$_1$}} = 0.5125\, \overrightarrow{\text{AB}} \quad~ \overrightarrow{\text{AN$_2$}} = 0.7125\, \overrightarrow{\text{AB}}.
 \end{equation}

\begin{figure}[h!] 
\centering{\includegraphics[height= 100mm]{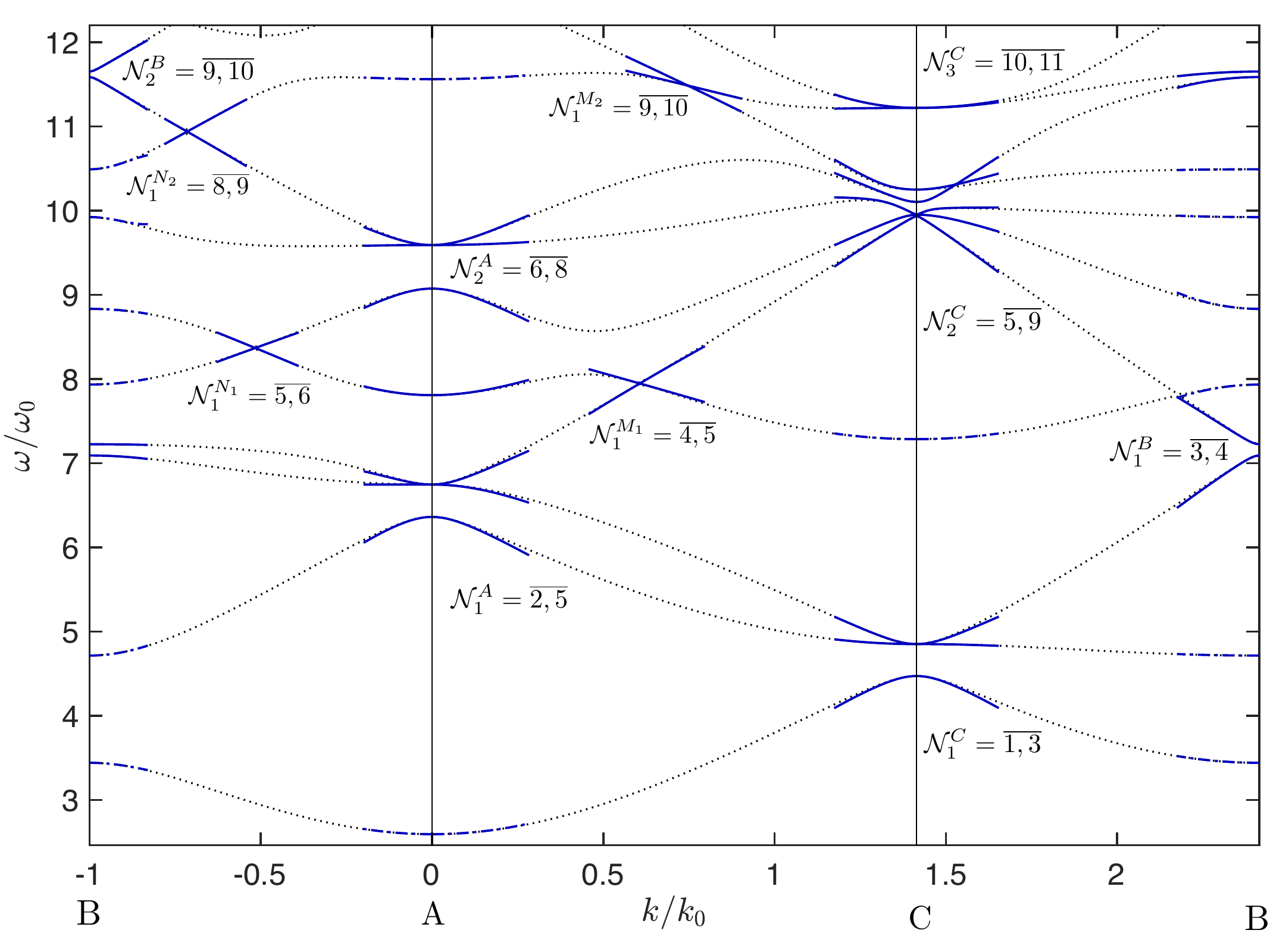}}\vspace*{-3mm}
\caption{Approximation of the dispersion branches 1--11 for the pinned square lattice near points~A, B, C, M$_1$, M$_2$, N$_1$ and N$_2$ in Fig~\ref{fig8}(c). In the diagram, dotted lines track the reference numerical results; solid lines signify the leading-order approximation of the clusters of nearby branches ($Q\!>\!1$), and dash-dotted lines plot the second-order approximation of isolated dispersion branches ($Q\!=\!1$). The normalization parameters are defined as $k_0=\pi/a$ and $\omega_0=\sqrt{G/(\rho a^2)}$.} \label{fig9} \vspace*{0mm}
\end{figure}

Fig.~\ref{fig9} examines the performance of the asymptotic models in terms of the first eleven branches of the dispersion relationship.  The reference numerical values along path BACB, as well as the cell functions at each wavenumber-eigenfrequency pair $(\bk_s,\omega_n(\bk_s))$ required to evaluate the asymptotic approximation, are computed via NGSolve by discretizing the unit cell with triangular elements of order 5 and maximum size $h_{max}=0.0175a$. As indicated earlier, the comparison is made in a neighborhood of the origin~A, apex points~B and~C, and internal points~M$_1$, M$_2$, N$_1$ and N$_2$, of the first Brillouin zone.

From Fig.~\ref{fig9}, one first observes that the introduction of ``pins'' (where~$u\!=\!0$) in an otherwise homogeneous medium results in both zero-frequency band gap, and another complete band gap just above the first dispersion branch. In terms of the asymptotic approximation, we note that clusters $\mathcal{N}^B_{1}$ and $\mathcal{N}^B_{2}$ of size $Q\!=\!2$ exhibit \emph{direction-dependent} behavior. Concerning $\mathcal{N}^B_{1}$ for example, we specifically find that the cluster's behavior in direction BC (where  $\boldsymbol{A}^{\gamma}(\hat\bk)$ is of full rank) is approximated by~\eqref{leadasympdb}, whereas in direction BA  (where  $\text{rank}(\boldsymbol{A}^{\gamma}(\hat\bk)\!=\!1$) \eqref{secasympdb} applies. Analogous comment applies to $\mathcal{N}^B_{2}$, with the roles of directions BC and BA reversed. 

Further, cluster $\mathcal{N}^C_{3}$ in Fig.~\ref{fig9} describes a repeated eigenvalue of multiplicity $Q=2$;  in this case $\boldsymbol{A}^{\tz}(\hat\bk)\!=\!\boldsymbol{0}$, and the dispersion relationships are approximated {by~\eqref{drasympta2}}. In contrast, clusters $\mathcal{N}^{N_1}_{1}$, $\mathcal{N}^{N_2}_{1}$, $\mathcal{N}^{M_1}_{1}$, and~$\mathcal{N}^{M_2}_{1}$, which carry the same size, are commonly characterized by the full-rank~$\boldsymbol{A}^{\tz}(\hat\bk)$ in the respective directions of the band diagram. Accordingly, these clusters are described (to the leading order) {by~\eqref{drasympta1}}.

Moving our attention to larger clusters, we note that $\mathcal{N}^A_{2}$ and $\mathcal{N}^C_{1}$ (with $Q\!=\!3$) each describe the case where {$\gamma\neq 0$ and} $\text{rank}(\boldsymbol{A}^{\gamma})\!=\!2$, examined in {sections~\S\ref{Apartialrankx} and~\S\ref{Dirac-like cones}}. With reference to Section~\ref{Dirac-like cones}, we specifically let $(n_1,n_2,n_3)_{\mathcal{N}^A_{2}}\!=\!(6,7,5)$ with~$n_0\!=\!6$, and $(n_1,n_2,n_3)_{\mathcal{N}^C_{1}}\!=\!(2,3,1)$ with~$n_0\!=\!2$. Finally, we note that the ``superclusters'' $\mathcal{N}^A_{1}$ (with $Q\!=\!4$) and $\mathcal{N}^C_{2}$ (with $Q\!=\!5$) are approximated very well by the leading-order model~\eqref{secasympdb}, with $\boldsymbol{A}^{\gamma}(\hat{\bk})$ being of full rank {in each relevant direction of the band diagram}.

%--------------------------------------------------------------------------------------------------------------%
\subsection{Forced medium motion ($\tilde{f}_\bk\neq 0$)}\label{disper_approx}
%--------------------------------------------------------------------------------------------------------------%

\noindent To illustrate the asymptotic approximation of the forced motion problem, we next examine the response of a two-dimensional Kagome lattice detailed in Section \ref{kagome_exp}, at frequency $\omega^2 = \omega_2^2(\boldsymbol{0})+\eps^2$ within the first band gap (see Fig.~\ref{fig6}). On recalling Claim \ref{BPWESR} and Remark \ref{commut}, for the source term~$f(\bx)$ according to~\eqref{forceexp} we assume Gaussian source distribution with~$\bk_s=\boldsymbol{0}$ (i.e. $\bk=\eps\hat\bk$) given by 
\begin{eqnarray} \label{FFF}
\tilde{f}_{\bk}(\bx) &\!\!\!=\!\!\!& F(\bk) \hh \phi(\bx), \qquad 
F(\bk) \:=\: |\mathcal{C}|\, \frac{a^2}{\pi}e^{-a^2\|\hat\bk\|^2} e^{- i\eps\hat\bk\cdot\bx_\circ}, \qquad \phi\in L^2_p(Y), 
\end{eqnarray}
and
\begin{equation}\label{dpoleM}
\phi(\boldsymbol{x}) = \sum_{k=0}^{M} \frac{4}{3\pi(2k\!+\!1)} \big[\sin\big(2\pi(2k\!+\!1) (\tfrac{1}{4}-x^1)\big)  + \sin\big(2\pi(2k\!+\!1) (\tfrac{1}{4}-x^2)\big) +\sin\big(2\pi(2k\!+\!1)(x^1+x^2)\big)\big],
\end{equation}
with $\bx_\circ=-0.25(\be_1+\be_2)$. Note that~\eqref{dpoleM} approximates a dipole on $Y$ as $M\to\infty$.

{We consider the response of the Kagome lattice for $\eps\in\{0.25,0.5\}$ in terms of the full field~$u^{\bpp}$ and the mean field $\langle u\rangle_{\!\rho}^{\bpp}$ ($p=0,1,2)$ according to~\eqref{leadt} and~\eqref{leade}. Thanks fo the exponential decay with~$\bk$ of~$\tilde{f}_\bk$ in~\eqref{leadt0}--\eqref{leadt2} and~\eqref{leade0}--\eqref{leade2}, in computations we conveniently supersede the domain of integration~$\mathcal{C}$ by~$[-\pi/a,\pi/a]\times [-\pi/a,\pi/a]$. To compute the featured eigenfunction, the cell functions, and the effective coefficients in~\eqref{leadt0}--\eqref{leadt2} and~\eqref{leade0}--\eqref{leade2}, we solve the the respective unit cell problems via NGSolve by discretizing the unit cell in Fig.~\ref{fig5}(b) with triangular elements of order 5 and maximum size $h_{max}=0.05a$.}

{ For the purposes of numerical verification, the reference global solution is computed via NGSolve over a parallelepiped $\Omega=\{\bx\!\in\!S,\: |x^j|\!<\!23,~j\: \in\!\overline{1,2}\}$ (subject to homogeneous Dirichlet boundary conditions) that is discretized with triangular elements of order 3 and maximum size $h_{max}= 0.05a$. Since the finite element simulations require the source distribution in space~$f(\bx)$, we similarly extend the domain of integration in~\eqref{forceexp} to~$\mathbb{R}^2$ in order to facilitate analytical evaluation. Note that the use of Dirichlet (in lieu of radiation) conditions on~$\partial\Omega$ is permitted by the exponential decay of~$u(\bx)$ at driving frequencies inside the band gap.}  

\begin{figure}[h!] 
\centering{\includegraphics[width= 160mm]{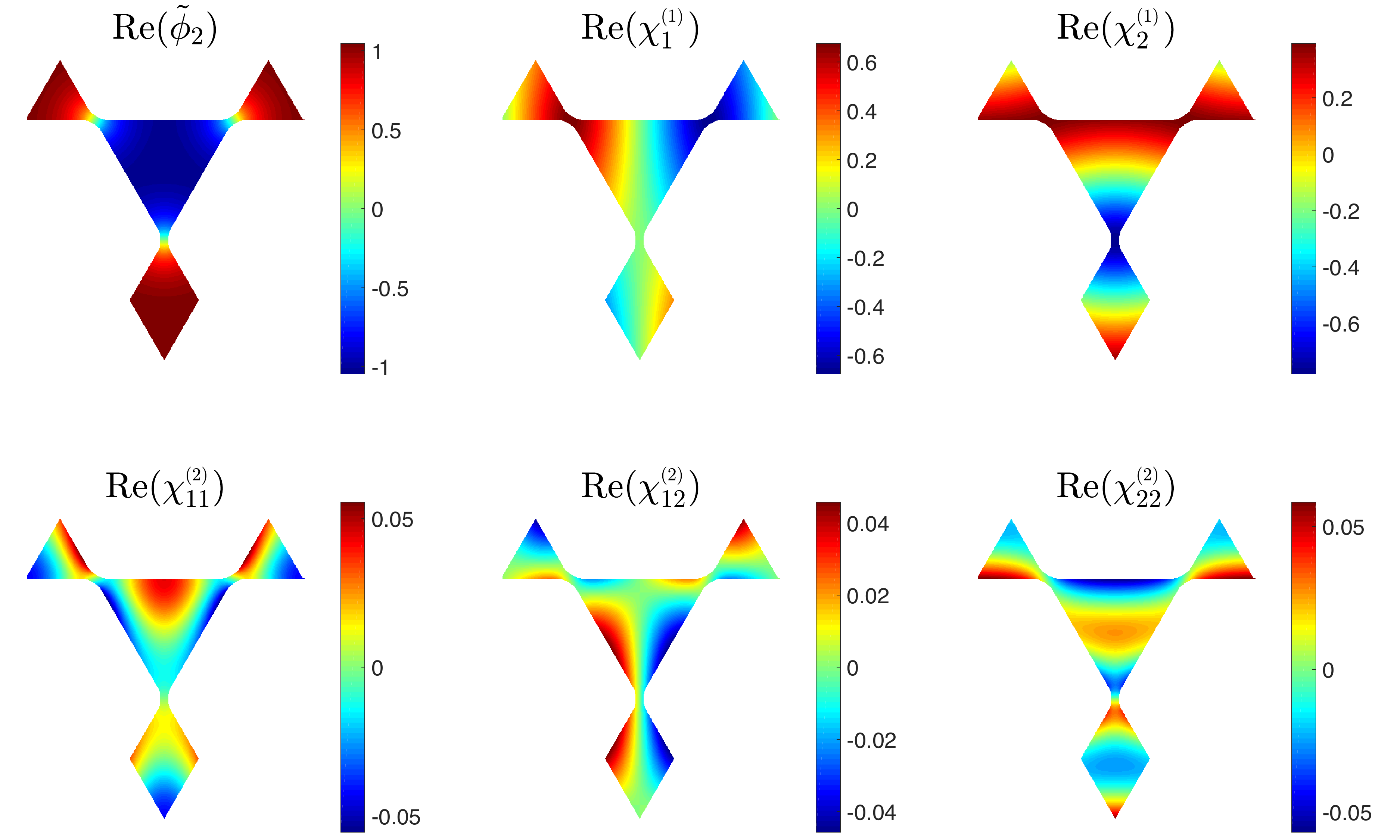}}\vspace*{-0mm}
\caption{Eigenfunction $\tilde{\phi}_2(\bx)$ and components of the cell functions $\bchi^{\so}(\bx)$ and $\bchi^{\st}(\bx)$ ($\bk_s=\boldsymbol{0}$, real parts only).} \label{fig10} \vspace*{0mm}
\end{figure}

\begin{figure}[h!] 
\centering{\includegraphics[width= 130mm]{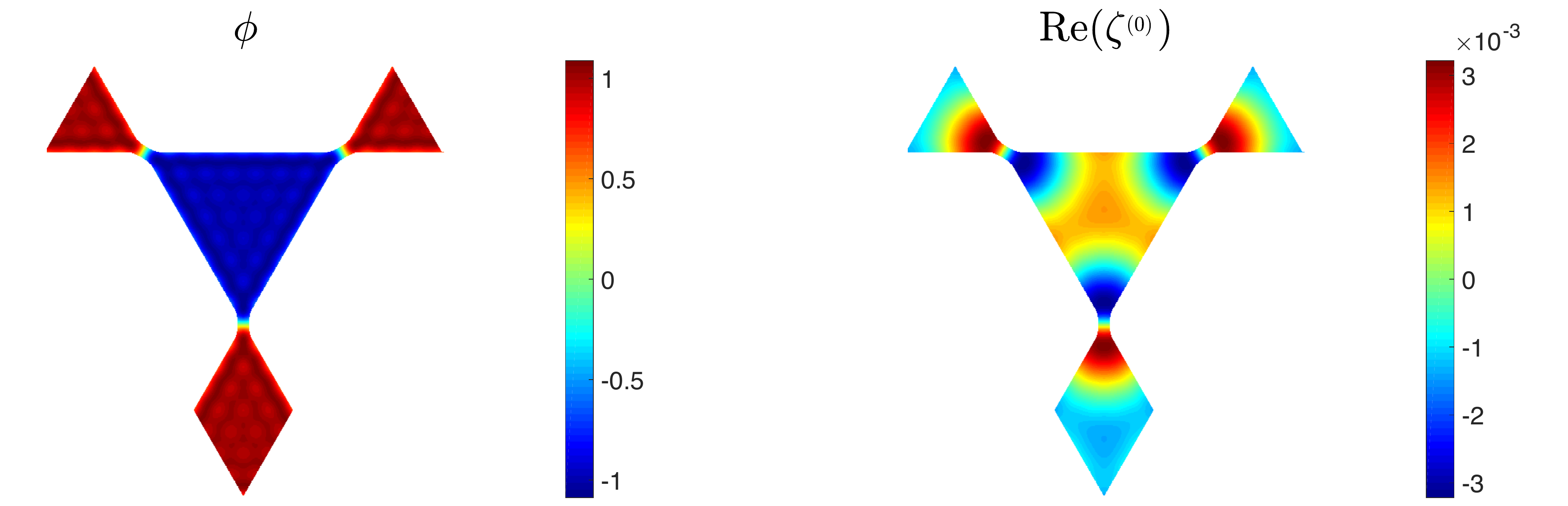}}\vspace*{-0mm}
\caption{ Periodic function $\phi(\bx)$ and affiliated cell function $\zeta^{\tz}(\bx)$ (real part only). } \label{fig11} \vspace*{0mm}
\end{figure}

\begin{figure}[h!] 
\centering{\includegraphics[width= 165mm]{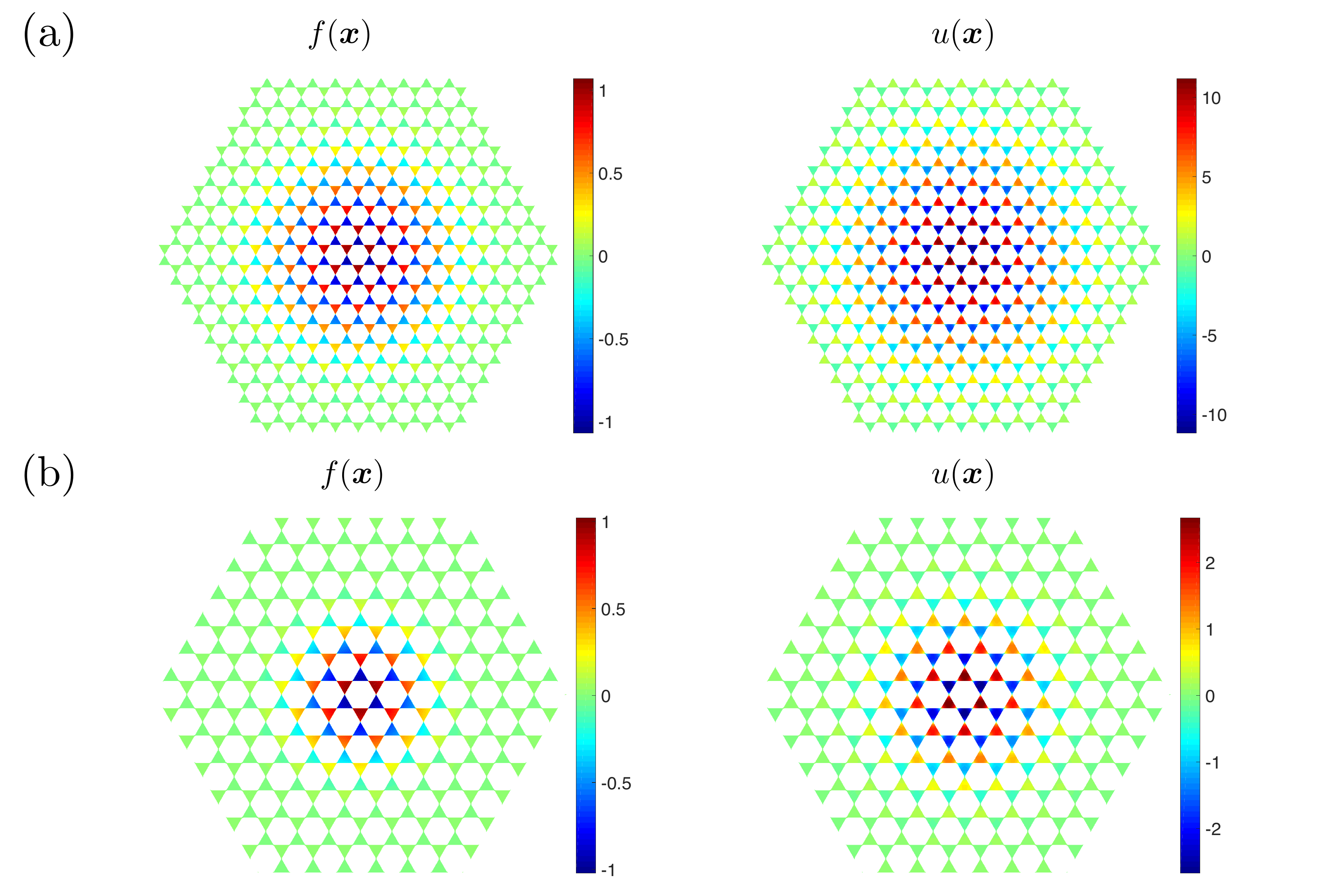}}\vspace*{-0mm}
\caption{ Source distributions~$f(\bx$) (left panels) and respective responses~$u(\bx)$ of the Kagome lattice (right panels) for: (a) $\eps=0.25$, and (b) $\eps=0.5$. The plots are made over a hexagonal subdomain of~$\Omega$ centered at~$\bx_\circ$ } \label{fig12} \vspace*{0mm}
\end{figure}

\begin{figure}[h!] 
\centering{\includegraphics[width= 175mm]{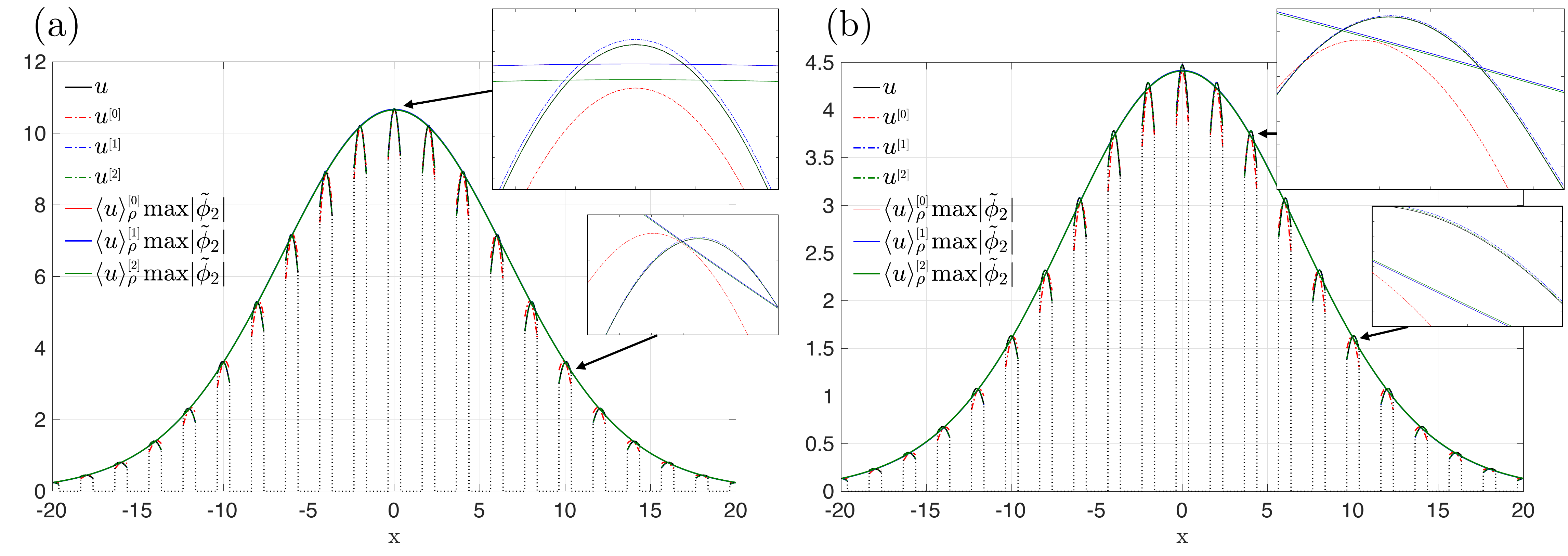}}\vspace*{-0mm}
\caption{Response of the Kagome lattice for $\eps=0.25$: numerical values versus leading-, first-, and second-order (full and effective) asymptotic approximations for (a) $\bx\cdot\boldsymbol{j}=1.5$, and (b) $\bx\cdot\boldsymbol{j}=8.43$. The leading- and first- order effective approximations overlap due to symmetry of the problem.} \label{fig13} \vspace*{0mm}
\end{figure}

\begin{figure}[h!] 
\centering{\includegraphics[width= 175mm]{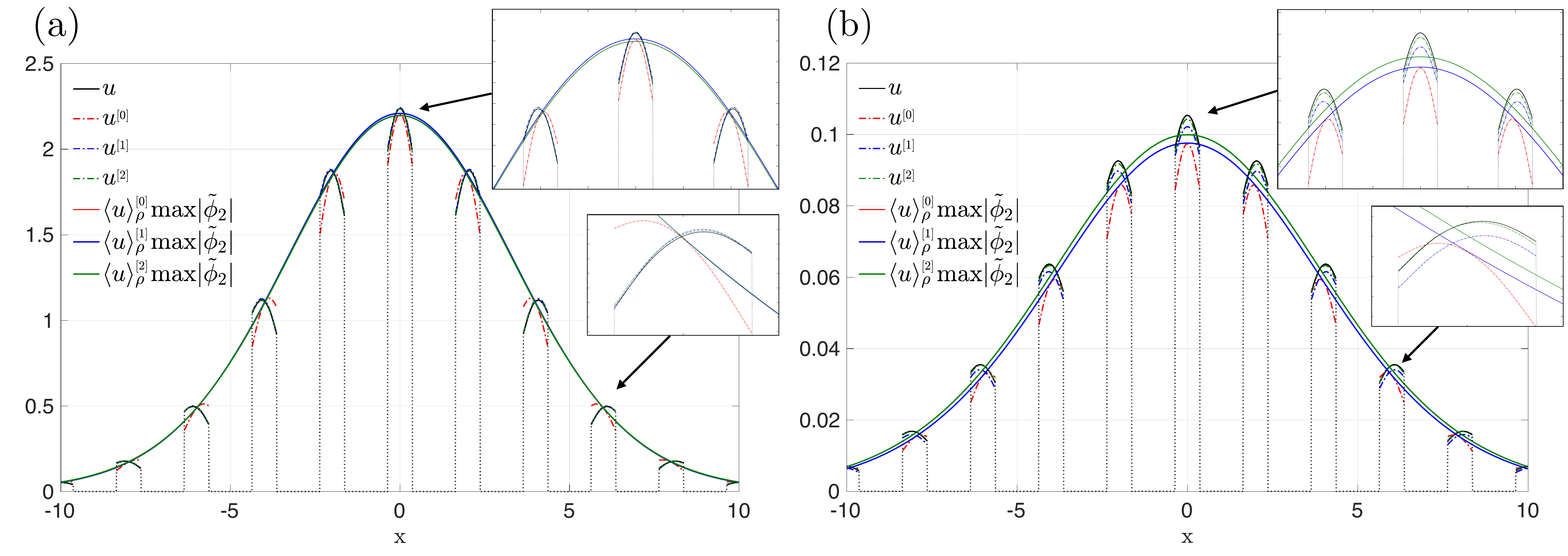}}\vspace*{-0mm}
\caption{Response of the Kagome lattice for $\eps=0.5$: numerical values versus leading-, first-, and second-order (full and effective) asymptotic approximations for (a) $\bx\cdot\boldsymbol{j}=1.5$, and (b) $\bx\cdot\boldsymbol{j}=8.43$. The leading- and first- order effective approximations overlap due to symmetry of the problem.} \label{fig14} \vspace*{0mm}
\end{figure}

{As an illustration, Fig.~\ref{fig10} plots the real parts of the eigenfunction $\tilde{\phi}_2$ and the components of $\bchi^{\so}$ and $\bchi^{\st}$. With reference to~\eqref{FFF}--\eqref{dpoleM}, Fig.~\ref{fig11} shows the periodic function $\phi$ for $M=8$ and the associated cell function $\zeta^{\tz}$ featured in the expression $\eta^{\tz}(\bk,\bx)=F(\bk)\zeta^{\tz}(\bx)$, see Remark \ref{eta0_lin}. Fig.~\ref{fig12} plots both the computed source distribution, $f(\bx)$, in the physical space and the NGSolve simulations of~$u(\bx)$ for $\eps\in\{0.25,0.5\}$ over a hexagonal truncation of~$\Omega$ centered at $\bx_\circ$. From the display, we observe that the medium response (i) conforms with the symmetries of the source and the Kagome lattice, and (ii) decays fast with~$\|\bx-\bx_\circ\|$ as expected for solutions inside a band gap. Fig.~\ref{fig13} (resp. Fig.~\ref{fig14}) compares, for  $\eps=0.25$ (resp.  $\eps=0.5$), the finite element response with the asymptotic approximations~$u^{\bpp}$ and~$\langle u\rangle_{\!\rho}^{\bpp}$ ($p=0,1,2$) across example cross-sections of the lattice. From the panels, one may observe both a clear increase in the fidelity of asymptotic approximation with $p$, and the ability of the effective solution (supported in~$\mathbb{R}^2$ instead of~$S$) to describe the essential response of the Kagome lattice.}

%%---------------------------------------------------------------------------------------------------------------------------------------------------------------------------------%
\section{Summary} \label{Summary}
%--------------------------------------------------------------------------------------------------------------%
\noindent In this work, we establish a rational framework for finite-wavenumber, finite-frequency (FW-FF) homogenization of the scalar wave equation in (generally non-orthogonal) periodic media with Dirichlet and Neumann exclusions. The proposed asymptotic ansatz applies to spectral neighborhoods of (i) an arbitrary wavenumber within the first Brillouin zone, and (ii) either simple, repeated, or nearby eigenfrequencies. With the aid of the Bloch-wave expansion that provides us with control over the wavenumber, we also account for, and homogenize, the source term featured in the scalar wave equation. A systematic asymptotic analysis of tightly-spaced eigenfrequency clusters (covering the most general asymptotic configuration) reveals an effective system of field equations featuring a ``matrix'' operator weaving Dirac- and wave-like behaviors, and a vector source term built from the projections of the (scalar) source term onto participating phonons, i.e. Bloch eigenfunctions. As numerical examples, we provide asymptotic description of the dispersion relationship for a Kagome lattice with Neumann exclusions, and a pinned square lattice. We also examine the performance of the effective model with a source term, up to the second order of asymptotic correction, by considering the response of the Kagome lattice to a dipole-like source acting near the edge of an internal band gap.

\section{Acknowledgments} \label{Ack}
%--------------------------------------------------------------------------------------------------------------%
\noindent BG kindly acknowledges the support provided by the endowed Shimizu Professorship during the course of this investigation.

%--------------------------------------------------------------------------------------------------------------%
\appendix\section{Supporting results} \label{Appendix}
%--------------------------------------------------------------------------------------------------------------%
%--------------------------------------------------------------------------------------------------------------%
\subsection{Bloch wave expansion (BWE)} \label{BWEappa}
%--------------------------------------------------------------------------------------------------------------%

\noindent Assume $f\in L^2(\SSS)$. Leting $\bx\in\SSS$, we define the function $\tilde{f}_{\bk}\in L^2_p(Y)$ as
\begin{equation}\notag 
\tilde{f}_{\bk}(\bx)= \sum_{\br\in\boldsymbol{R}}f(\bx+\br)e^{- i \bk\cdot(\bx+\br)}.
\end{equation}
\noindent  For $\bk_s\in\mathbb{R}^d$, we have
\begin{eqnarray}
\int_{\bk_s\!+\mathcal{B}} \frac{1}{|\mathcal{B}|} \tilde{f}_{\bk}(\bx)e^{ i \bk\cdot\bx}d\bk  ~=\: \int_{\bk_s\!+\mathcal{B}} \frac{1}{|\mathcal{B}|} \sum_{\br\in\boldsymbol{R}} f(\bx+\br)e^{- i \bk\cdot\br} d\bk ~=\:
\sum_{\br\in\boldsymbol{R}} f(\bx+\br)\, \frac{e^{- i \bk_s\cdot\br}}{|\mathcal{B}|} \int_{\mathcal{B}} e^{- i \bk\cdot\br}d\bk ~=~ f(\bx), 
\end{eqnarray}
since the last integral vanishes for all~$\br\neq\boldsymbol{0}$ due to~\eqref{Brillouin}.

%--------------------------------------------------------------------------------------------------------------%
\subsection{Relationship between the plane wave expansion (PWE) and BWE} \label{RB2PWEappa}
%--------------------------------------------------------------------------------------------------------------%

\noindent Let $\psi(\bk)$ be a $Y_0^*$-periodic function defined on the reciprocal space $\mathbb{R}^d$ as 
\begin{equation} \notag
\psi(\bk)= \sum_{\br\in\boldsymbol{R}}~e^{i\br\cdot\bk}. 
\end{equation}
For $\br^*\!\in\boldsymbol{R}^*$, it is clear that $\psi(\bk+\br^*)= \psi(\bk)$ since~$\br\cdot\br^*\in 2\pi \mathbb{Z}$ by~\eqref{Brav} and~\eqref{recBrav}. Next, let $\phi(\bk)$ be $Y_0^*$-periodic and square-integrable over $Y_0^*$ and let $\Phi(\br)$, with $\br\in\boldsymbol{R}$, denote its Fourier series coefficient. For $\bk_{\circ}\!\in\mathbb{R}^d$, we then have 
\begin{eqnarray} \notag
\int_{Y_0^*} \psi(\bk-\bk_{\circ}) \phi(\bk) d\bk &\!\!=\!\!\!& \int_{Y_0^*}~\Big(\sum_{\br\in\boldsymbol{R}}~e^{i\br\cdot(\bk-\bk_{\circ})}\Big)  \phi(\bk)d\bk ~=\: 
\sum_{\br\in\boldsymbol{R}}~e^{-i\br\cdot\bk_{\circ}}~\int_{Y_0^*}~e^{i\br\cdot\bk}\phi(\bk)d\bk\\
&=\!& \sum_{\br\in\boldsymbol{R}} e^{-i\br\cdot\bk_{\circ}} \: |Y_0^*|\: \Phi(-\br) \notag  ~=\: |Y_0^*|\sum_{\br\in\boldsymbol{R}} e^{i\br\cdot\bk_{\circ}} \: \Phi(\br) \notag \\*[1mm]
&=\!& |Y_0^*| \; \phi(\boldsymbol{\bk_{\circ}}) ~=~  |\mathcal{B}| \; \phi(\boldsymbol{\bk_{\circ}}\!+\bc^*), \quad \bc^*\in\bR^*. \label{psi}  \label{aux1}
\end{eqnarray}

Next, recall the PWE of $f\!\in\!L_2(\SSS)$ given by~\eqref{PWE}. For given $\bk_{\circ}\!\in\bk_s\!+\mathcal{B}$, we first define a unique translation vector~$\bc^*(\bk_\circ)\in\bR^*$ such that 
\[
\bk_{\circ}\!\in\bk_s\!+\mathcal{B} \quad \rightarrow \quad \bk_\circ +\bc^*(\bk_{\circ}) \in Y_0^*. 
\]
With reference to definition~\eqref{BWE2} of the Bloch transform and~\eqref{PWE}, one accordingly has
\begin{eqnarray}\label{RPBWEd}
\tilde{f}_{\bk_{\circ}}(\bx) &= & \sum_{\br\in\boldsymbol{R}}f(\bx+\br)~e^{ -i \bk_{\circ}\cdot(\bx+\br)} ~=~ \sum_{\br\in\boldsymbol{R}}~\int_{\mathbb{R}^d} \mathfrak{F}(\bk) \hh e^{i (\bk-\bk_{\circ})\cdot(\bx+\br)} \dd \bk \notag \\
&= & \int_{\mathbb{R}^d}  \mathfrak{F}(\bk) \hh e^{i (\bk-\bk_{\circ})\cdot\bx}~\Big( \sum_{\br\in\boldsymbol{R}}~e^{i (\bk-\bk_{\circ})\cdot\br} \Big) \dd \bk ~=~
\int_{\mathbb{R}^d}  \mathfrak{F}(\bk) \hh e^{i (\bk-\bk_{\circ})\cdot\bx}~ \psi(\bk-\bk_{\circ}) \dd \bk \notag \\
&= & \sum_{\br^*\in\boldsymbol{R}^*} \int_{\br^*+Y_0^*}  \mathfrak{F}(\bk) \hh e^{i (\bk-\bk_{\circ})\cdot\bx} \: \psi(\bk-\bk_{\circ}) \dd \bk. \notag \\
&= & \sum_{\br^*\in\boldsymbol{R}^*}~\int_{Y_0^*} \mathfrak{F}(\br^*\!+\bk) \hh e^{i (\br^*\!+\bk-\bk_{\circ})\cdot\bx} \: \psi(\bk-\bk_{\circ}) \dd \bk \label{proof2x} \\
&= &  |\mathcal{B}| \sum_{\br^*\in\boldsymbol{R}^*} \mathfrak{F}(\br^*\!+\bk_{\circ}+\bc^*(\bk_{\circ})) \hh e^{i(\br^*+\bc^*(\bk_{\circ}))\!\cdot\bx} ~=~ |\mathcal{B}| \sum_{\br^*\in\boldsymbol{R}^*} \mathfrak{F}(\br^*\!+\bk_{\circ}) \hh e^{i\br^*\!\cdot\bx}. \label{proof2}
\end{eqnarray}
Note that~\eqref{proof2} is obtained from~\eqref{proof2x} by (i) setting, for any given~$\br^*$, $\phi(\bk):=\mathfrak{F}(\br^*\!+\bk) \hh e^{i (\br^*\!+\bk-\bk_{\circ})\cdot\bx}$ over~$Y_0^*$; (ii) extending the support of such defined~$\phi(\bk)$ to~$\mathbb{R}^d$ by the application of~$Y_0^*$-periodicity, and (iii) applying~\eqref{aux1} with~$\bc^*=\bc^*(\bk_\circ)$ to \emph{each term} in the sum. This establishes claim~\eqref{RPBWE}. Assuming further that~$\mathfrak{F}(\bk)$ is compactly supported inside $\bk_s\!+\mathcal{B}$, we obtain~\eqref{RPBWES} thanks to the fact that~$\mathfrak{F}(\br^*+\bk)=0$ for all~$\br^*\neq \boldsymbol{0}$.

%--------------------------------------------------------------------------------------------------------------%
\subsection{Effective coefficients} \label{Effcoef}
% -----------------------------------------------------------------------------------------------------------------------------------------------------------------------------------%

\noindent \textbf{Proof of Claim \ref{mu0real}.}~~We integrate $\{ \overline{(\eqref{chi1}\!\otimes  \overline{\bchi^{\so}},1)}\}$ using the divergence theorem and the flux boundary conditions \eqref{gbcs} written in terms of $G(\nabla_{\!\bk_s}\hh\bchi^{\so}\!+\tilde{\phi}_n \boldsymbol{I})$ to obtain
\begin {eqnarray}\label{mu0ae}
\bmu^{\tz} = \langle G \tilde{\phi_n} \boldsymbol{I} \rangle + \lambda_n \{(\rho \bchi^{\so}\!\otimes  \overline{\bchi^{\so}},1)\} - \{(G( \overline{\nabla_{\!\bk_s}\hh\bchi^{\so}})\cdot (\nabla_{\!\bk_s}\hh\bchi^{\so}),1)\} ~\in \mathbb{R}^{d\times d}.
\end{eqnarray}
\begin{remark}
The ``dot'' operator in \eqref{mu0ae} and thereafter assumes inner tensor contraction with respect to the first index. 
\end{remark}

\noindent \textbf{Proof of Claim \ref{theta1r}.}~~We evaluate integrals $\{(\eqref{chi2}\!\otimes  \overline{\bchi^{\so}},1)\}$ and $\{ \overline{(\eqref{chi1}\!\otimes  \overline{\bchi^{\st}},1)}\}$ by applying the divergence theorem and exploiting the boundary conditions \eqref{gbcs} satisfied by $G(\nabla_{\!\bk_s}\hh\bchi^{\so}\!+\tilde{\phi}_n \boldsymbol{I})$ and $G(\nabla_{\!\bk_s}\hh\bchi^{\st}+ \{\boldsymbol{I}\!\otimes\bchi^{\so}\}')$. In this way, we find that
\begin {eqnarray}
\bth^{\so} &\!\!\!=\!\!\!& \{(G\hh\bchi^{\so}\!\otimes  \overline{\nabla_{\!\bk_s}\bchi^{\so}},1)\} -\{(G\hh\overline{\bchi^{\so}}\!\otimes \nabla_{\!\bk_s}\bchi^{\so},1)\} \,+\, \{(G\hh\boldsymbol{I}\!\otimes\bchi^{\so},\tilde{\phi}_n)\} -\{ \overline{(G\boldsymbol{I}\!\otimes\bchi^{\so},\tilde{\phi}_n)}\} \notag \\
&& ~+\, \{ \bth^{\tz} \!\otimes \frac{1}{\rho^{\tz}} \{(\rho\bchi^{\so}\!\otimes \overline{\bchi^{\so}},1 )\}\} ~\in i\hh \mathbb{R}^{d\times d \times d}.\label{bth1exp}
\end{eqnarray}

\noindent \textbf{Proof of Claim \ref{mu2r}.}~~We evaluate  integrals $\{(\eqref{chi3}\!\otimes  \overline{\bchi^{\so}},1)\}$, $\{ \overline{(\eqref{chi1}\!\otimes  \overline{\bchi^{\sh}},1)}\}$ and $\{ \overline{(\eqref{chi2}\!\otimes  \overline{\bchi^{\st}},1)}\}$ by applying the divergence theorem and exploiting the boundary conditions \eqref{gbcs} satisfied by $G(\nabla_{\!\bk_s}\hh\bchi^{\so}\!+\tilde{\phi}_n \boldsymbol{I})$, $G(\nabla_{\!\bk_s}\hh\bchi^{\st}+ \{\boldsymbol{I}\!\otimes\bchi^{\so}\}')$ and $G(\nabla_{\!\bk_s}\hh\bchi^{\sh}+ \{\boldsymbol{I}\!\otimes\bchi^{\st}\}')$. We then obtain 
\begin {eqnarray}
\bmu^{\st} &\!\!\!=\!\!\!& -\lambda_n \{(\rho \bchi^{\st}\!\otimes  \overline{\bchi^{\st}},1)\} + \{(G( \overline{\nabla_{\!\bk_s}\hh\bchi^{\st}})\cdot (\nabla_{\!\bk_s}\hh\bchi^{\st}),1)\}
- \{(G\boldsymbol{I}\!\otimes\bchi^{\so}\!\otimes \overline{\bchi^{\so}},1)\} \notag \\
&& ~+\, \{\bmu^{\tz}\!\otimes \frac{1}{\rho^{\tz}} \{(\rho\bchi^{\so}\!\otimes \overline{\bchi^{\so}},1 )\}\} ~\in \mathbb{R}^{d\times d \times d \times d}.
\end{eqnarray}

\noindent \textbf{Proof of Claim \ref{Bherm}.}~~We integrate $\{(\eqref{chi1q}\!\otimes  \overline{\bchi^{\so}_p},1)\}$ by applying the divergence theorem and making use of the boundary conditions \eqref{gbcs} satisfied by $G(\nabla_{\!\bk_s}\hh\bchi^{\so}_q+\tilde{\phi}_{n_p} \boldsymbol{I})$ to obtain
\begin {eqnarray}
\bmu^{\tz}_{qp} &\!\!\!=\!\!\!& \langle G \tilde{\phi}_{n_p} \boldsymbol{I} \rangle^{n_{\nes q}} + \lambda_n \{(\rho  \overline{\bchi_q^{\so}}\!\otimes \bchi_p^{\so},1)\} - \{(G(\nabla_{\!\bk_s}\hh\bchi_p^{\so})\cdot ( \overline{\nabla_{\!\bk_s}\hh\bchi_q^{\so}}),1)\}. \notag
\end{eqnarray}
As a result, for $p\in\overline{1,Q}$ we have
\begin {eqnarray}
\bmu^{\tz}_{qp}:(i\hat{\bk})^2 &\!\!\!=\!\!\!&\Big( \langle G \tilde{\phi}_{n_p} \boldsymbol{I} \rangle^{n_{\nes q}} + \lambda_n \{(\rho  \overline{\bchi_q^{\so}}\!\otimes \bchi_p^{\so},1)\} - \{(G(\nabla_{\!\bk_s}\hh\bchi_p^{\so})\cdot ( \overline{\nabla_{\!\bk_s}\hh\bchi_q^{\so}}),1)\} \Big):(i\hat{\bk})^2 ~=~  \overline{\bmu^{\tz}_{pq}:(i\hat{\bk})^2}, \notag 
\end{eqnarray}
which demonstrates that the matrix $\boldsymbol{B}^{\tz}$ is Hermitian.\\

\noindent \textbf{Proof of Claim \ref{psiortho}.}~~By construction, family $\tilde{\psi}_q$, $q\in\overline{N_0,N}$ is $\rho$-orthogonal. Thus for $p,q\in\overline{N_0,N}$ we have 
\begin {eqnarray} \notag 
(\rho\tilde{\psi}'_p,\tilde{\psi}'_q) &\!\!\!=\!\!\!& (\rho\tilde{\psi}_p,\tilde{\psi}_q) ~=~ \delta_{pq}  (\rho\tilde{\psi}_p,\tilde{\psi}_p).
\end{eqnarray}
\noindent For $p,q\in\overline{1,N_0}$, we have 
\begin {eqnarray}
(\rho\tilde{\psi}'_p,\tilde{\psi}'_q) &=& \sum_{s=1}^{N_0}\sum_{r=1}^{N_0} \bar{P}_{sp} \overline{\bar{P}_{rq}}(\rho\tilde{\psi}_s,\tilde{\psi}_r) ~=~ \sum_{s=1}^{N_0} \bar{P}_{sp} \overline{\bar{P}_{sq}}(\rho\tilde{\psi}_s,\tilde{\psi}_s) ~=~ \sum_{s=1}^{N_0} \bar{P}_{sp} \overline{\bar{P}_{sq}}\rho^{\tz}_{s} ~=~ \sum_{s=1}^{N_0}  \overline{\bar{P}_{sq}}D_{ss}\bar{P}_{sp}\notag \\
&=& \delta_{pq} \Big(\sum_{s=1}^{N_0}  \overline{\bar{P}_{sp}}D_{ss}\bar{P}_{sp} \Big). \notag
\end{eqnarray}
Further, for $p\in\overline{1,N_0}$ and $q\in\overline{N_0,N}$, we have 
\begin {eqnarray}
(\rho\tilde{\psi}'_p,\tilde{\psi}'_q) ~=~ \sum_{s=1}^{N_0}\bar{P}_{sp} (\rho\tilde{\psi}_s,\tilde{\psi}_q) ~=~ 0, \notag
\end{eqnarray}
whereby family $\tilde{\psi}'_p$ is $\rho$-orthogonal as well.
%--------------------------------------------------------------------------------------------------------------%
\subsubsection{Cell functions identities } \label{identities}
%--------------------------------------------------------------------------------------------------------------%
\noindent \textbf{Proof of Claim \ref{eta_chi1_id}.}~~Identity \eqref{eta0_chi1_id} is obtained by evaluating integrals $ \overline{(\eqref{chi1}\hh\overline{\eta^{\tz}},1)}$ and $(\eqref{eta0}\hh\overline{\bchi^{\so}},1)$ via the divergence theorem and use of the boundary conditions \eqref{gbcs} satisfied by $G\nabla_{\!\bk_s}\eta^{\tz}$ and  $G(\nabla_{\!\bk_s}\hh\bchi^{\so}\!+\tilde{\phi}_n \boldsymbol{I})$.\\

\noindent \textbf{Proof of Claim \ref{eta2_idc}.}~~Identity \eqref{eta2_ide} is obtained by evaluating integrals $ \overline{(\eqref{chi1}\hh\overline{\eta^{\st}},1)}$ and $(\eqref{eta_2}\hh\overline{\bchi^{\so}},1)$ via the divergence theorem and use of the boundary conditions \eqref{gbcs} satisfied by $G\nabla_{\!\bk_s}\eta^{\st}$ and  $G(\nabla_{\!\bk_s}\hh\bchi^{\so}\!+\tilde{\phi}_n \boldsymbol{I})$.\\

\noindent \textbf{Proof of Claim \ref{eta1_id}.}~~Identity \eqref{eta1_id} is obtained by evaluating integrals $ \overline{(\eqref{chi1}\hh \overline{\boldsymbol{\eta}^{\so}},1)}$, $(\eqref{eta_1}\hh \overline{\bchi^{\so}},1)$ and $ \overline{(\eqref{chi2}\hh \overline{\eta^{\tz}},1)}$ via the divergence theorem and use of the boundary conditions \eqref{gbcs} satisfied by $G\nabla_{\!\bk_s}\eta^{\tz}$,  $G(\nabla_{\!\bk_s}\hh\boldsymbol{\eta}^{\so}+\eta^{\tz} \boldsymbol{I})$ and  $G(\nabla_{\!\bk_s}\hh\bchi^{\so}\!+\tilde{\phi}_n \boldsymbol{I})$.
%--------------------------------------------------------------------------------------------------------------%
\subsubsection{Cell functions and effective coefficients at special wavenumbers  } \label{specases}
%--------------------------------------------------------------------------------------------------------------%
\noindent \textbf{Proof of Claim \ref{specase1}.}~~At a wavenumber-eigenfreqency pair $(\bk_s,\omega_n)$, $n\in\mathbb{Z}^+$, where $\bk_s = \frac{1}{2}(\sum_j n_j\hh\be^j)$, $n_j\in\{-1,0,1\}$, and $\omega_n$ is simple, the Bloch function $\phi_n(\bx)=\tilde{\phi}_n(\bx) e^{ i \bk_s\cdot\bx}$ with $\tilde{\phi}_n\!\in\!H^1_{p0}(Y)$ uniquely solves the field equation
\begin{eqnarray}\label{neigp}
-\tilde{\lambda}_n\rho(\bx)\hh {\phi}_n - \nabla\sip \big(G(\bx)\nabla{\phi}_n\big) ~=~0 \quad\text{in~~}Y,
\end{eqnarray}
subject to the boundary conditions  
\begin{eqnarray}
{\phi}_n|_{\partial{Y^{\tpp}_{j0}}} &\!\!\!=\!\!\!& e^{-i\bk_s\cdot\be_j}{\phi}_n|_{\partial{Y^{\tpp}_{j1}}}, \notag\\
\bnu\cdot G\nabla{\phi}_n|_{\partial{Y^{\tpp}_{j0}}} & =&- e^{-i\bk_s\cdot\be_j}\bnu\cdot G\nabla{\phi}_n|_{\partial{Y^{\tpp}_{j1}}} ,  \notag \\%\quad j\in\overline{1,d},  \\
\bnu\cdot G\nabla{\phi}_n|_{\partial Y^{\tnn}} &\!\!\!=\!\!\!& 0,\notag \\
{\phi}_n|_{\partial Y^{\tdd}}  &\!\!\!=\!\!\!& 0,  \label{nbcs}
\end{eqnarray}
where $e^{-i\bk_s\cdot\be_j}=(-1)^{n_j}$, i.e. $e^{-i\bk_s\cdot\be_j}\in\{-1,1\}$. Hence, by the uniqueness of the solution, we find that the real and imaginary parts of $\phi_n(\bx)$ are proportional, which factorizes $\phi_n(\bx)$ into a (real-valued function) $\times$ $e^{i\varphi_0}$, and we let~$\varphi_0=0$ without affecting the normalization \eqref{evpo}. We then obtain 
\begin{eqnarray}
\bth^{\tz}  &\!\!\!=\!\!\!&  \langle G\nabla_{\!\bk_{s}} \tilde{\phi}_n \rangle -  \overline{\langle G\nabla_{\!\bk_{s}}\tilde{\phi}_n\rangle} \notag \\
&\!\!\!=\!\!\!& 2i~\text{Im}( \langle G\nabla_{\!\bk_{s}} \tilde{\phi}_n \rangle) \,=\, 2i~\text{Im}( (G\nabla{\phi}_n, \phi_n )) \:=\: \boldsymbol{0}. \notag 
\end{eqnarray}
Similarly, we define $\bX^{\so}(\bx)=\bchi^{\so}(\bx) e^{ i \bk_s\cdot\bx}$ with~$\bchi^{\so}\in(\bar{H}^{1}_{p0}(Y))^d$ that uniquely solves the field equation
\begin{equation}
\tilde{\lambda}_n\hh\rho\hh\bX^{\so} + \nabla_{\!\bk_s}\!\sip \big(G(\nabla\hh\bX^{\so}+{\phi}_n \boldsymbol{I})\big) + G\nabla{\phi}_n \;=\; \boldsymbol{0},
\end{equation} 
subject to the boundary conditions
\begin{eqnarray}
\bX^{\so}|_{\partial{Y^{\tpp}_{j0}}} &\!\!\!=\!\!\!& e^{-i\bk_s\cdot\be_j}\bX^{\so}|_{\partial{Y^{\tpp}_{j1}}}, \notag\\
\bnu\cdot G(\nabla\hh\bX^{\so}+{\phi}_n \boldsymbol{I})|_{\partial{Y^{\tpp}_{j0}}} & =&- e^{-i\bk_s\cdot\be_j}\bnu\cdot G(\nabla\hh\bX^{\so}+{\phi}_n \boldsymbol{I})|_{\partial{Y^{\tpp}_{j1}}} ,  \notag \\%\quad j\in\overline{1,d},  \\
\bnu\cdot G(\nabla\hh\bX^{\so}+{\phi}_n \boldsymbol{I})|_{\partial Y^{\tnn}} & = &  \boldsymbol{0}, \notag \\
\bX^{\so}|_{\partial Y^{\tdd}}  &\!\!\!=\!\!\!& \boldsymbol{0}, \notag
\end{eqnarray}
\noindent from which we conclude that $\bX^{\so}(\bx)$ is real-valued up to a constant factor $e^{i\varphi_0}$. \\

\noindent \textbf{Proof of Claim \ref{specase2}.}~~On recalling Claim~\ref{specase1} and~\eqref{bth1exp}, we have
\begin {eqnarray}
\bth^{\so} &\!\!\!=\!\!\!& \{(G\bX^{\so}\!\otimes \overline{\nabla\bX^{\so}},1)\} -\{(G\hh\overline{\bX^{\so}}\!\otimes \nabla\bX^{\so},1)\} \notag \\
&& +\; \{(G\boldsymbol{I}\!\otimes\bX^{\so},{\phi}_n)\} -\{ \overline{(G\boldsymbol{I}\!\otimes\bX^{\so},{\phi}_n)}\} + \{ \bth^{\tz} \!\otimes \frac{1}{\rho^{\tz}} \{(\rho\bX^{\so}\!\otimes \overline{\bX^{\so}},1 )\}\} ~=~ \boldsymbol{0}. \label{bth1val}
\end{eqnarray}
We next define the Bloch functions $\bX^{\st}(\bx)=\bchi^{\st}(\bx) e^{ i \bk_s\cdot\bx}$ and $\bX^{\sh}(\bx)=\bchi^{\sh}(\bx) e^{ i \bk_s\cdot\bx}$,  with~$\bchi^{\st}\!\in\big(\bar{H}^{1}_{p0}(Y)\big)^{d\times d}$ and~$\bchi^{\sh}\!\in\big(\bar{H}^{1}_{p0}(Y)\big)^{d\times d \times d}$, that uniquely solve the respective boundary value problems 
\begin{eqnarray}
&\tilde{\lambda}_n\hh\rho\hh\bX^{\st} + \nabla_{\!\bk_s}\!\sip \big(G(\nabla\hh\bX^{\st}+\{\boldsymbol{I}\!\otimes\bX^{\so}\}')\big) + G(\nabla\hh\bX^{\so}+{\phi}_n \boldsymbol{I})= \frac{\rho}{\rho^{\tz}}{\phi}_n\hh\bmu^{\tz},
\end{eqnarray}
\begin{eqnarray}
\bX^{\st}|_{\partial{Y^{\tpp}_{j0}}} &\!\!\!=\!\!\!& e^{-i\bk_s\cdot\be_j}\bX^{\st}|_{\partial{Y^{\tpp}_{j1}}}, \notag\\
\bnu\cdot G(\nabla\hh\bX^{\st}+\{\boldsymbol{I}\!\otimes\bX^{\so}\}')|_{\partial{Y^{\tpp}_{j0}}}  & =&- e^{-i\bk_s\cdot\be_j}\bnu\cdot G(\nabla\hh\bX^{\st}+\{\boldsymbol{I}\!\otimes\bX^{\so}\}')|_{\partial{Y^{\tpp}_{j1}}} ,  \notag \\%\quad j\in\overline{1,d},  \\
\bnu\cdot G(\nabla\hh\bX^{\st}+\{\boldsymbol{I}\!\otimes\bX^{\so}\}')|_{\partial Y^{\tnn}} & = &  \boldsymbol{0}, \notag \\
\bX^{\st}|_{\partial Y^{\tdd}}  &\!\!\!=\!\!\!& \boldsymbol{0}, \notag
\end{eqnarray}
where $\bth^{\tz}=\boldsymbol{0}$ due to Claim~\ref{specase1}, and
\begin{eqnarray}
&\tilde{\lambda}_n\hh\rho\hh\bX^{\sh} + \nabla\sip \big(G(\nabla\hh\bX^{\sh}+\{\boldsymbol{I}\!\otimes\bX^{\st}\}' )\big) + G\{ \nabla\hh\bX^{\st}+\boldsymbol{I}\!\otimes\bX^{\so}\} = \frac{\rho}{\rho^{\tz}}\hh\big \{\bmu^{\tz}\!\otimes \bX^{\so} \big \},
\end{eqnarray}
\begin{eqnarray}
\bX^{\sh}|_{\partial{Y^{\tpp}_{j0}}} &\!\!\!=\!\!\!& e^{-i\bk_s\cdot\be_j}\bX^{\sh}|_{\partial{Y^{\tpp}_{j1}}}, \notag\\
\bnu\cdot G(\nabla\hh\bX^{\sh}+\{\boldsymbol{I}\!\otimes\bX^{\st}\}')|_{\partial{Y^{\tpp}_{j0}}}  & =&- e^{-i\bk_s\cdot\be_j}\bnu\cdot G(\nabla\hh\bX^{\sh}+\{\boldsymbol{I}\!\otimes\bX^{\st}\}')|_{\partial{Y^{\tpp}_{j1}}}  \notag \\%\quad j\in\overline{1,d},  \\
\bnu\cdot G(\nabla\hh\bX^{\sh}+\{\boldsymbol{I}\!\otimes\bX^{\st}\}')|_{\partial Y^{\tnn}} & = &  \boldsymbol{0}, \notag \\
\bX^{\sh}|_{\partial Y^{\tdd}}  &\!\!\!=\!\!\!& \boldsymbol{0}, \notag
\end{eqnarray}
where $\bth^{\tz}=\boldsymbol{0}$ due to Claim~\ref{specase1} and $\bth^{\so}=\boldsymbol{0}$ by \eqref{bth1val}. Recalling the properties of $\phi_n(\bx)$ and $\bX^{\so}(\bx)$, we then conclude that $\bX^{\st}(\bx)$ and $\bX^{\sh}(\bx)$ are each real-valued up to a constant factor $e^{i\varphi_0}$. \\

\noindent \textbf{Proof of Claim \ref{Aimag}.} At wavenumber-eigenfrequency pair $(\bk_s,\omega_n)$ where $\bk_s = \frac{1}{2}(\sum_j n_j\hh\be^j)$ $n_j\in\{-1,0,1\}$, and $\omega_n$ is of multiplicity~$Q\!>\!1$, Bloch functions $\phi_{n_p}(\bx)=\tilde{\phi}_{n_p}(\bx) e^{ i \bk_s\cdot\bx}$ ($p\!=\!\overline{1,Q}$) independently solve the boundary value problem~\eqref{neigp}--\eqref{nbcs}. In a way that is similar to the simple eigenfrequency case, one finds that $ \overline{\phi_{n_p}}$ are also solutions of~\eqref{neigp}--\eqref{nbcs}, and so are~$\text{Re}(\phi_{n_p})$ and~$\text{Im}(\phi_{n_p})$. Thus, by proceeding with the Gram-Schmidt $\rho$-orthogonalization of the $2Q$ \emph{real-valued} solutions, we can obtain $Q$ real-valued solutions that are $\rho$-orthogonal and normalized in the sense of~\eqref{evpo}. On relabeling the new (real-valued) basis using original notation, we find  
\begin{eqnarray}
\bth_{pq}^{\tz} &\!\!\!=\!\!\!&  \langle G\nabla_{\!\bk_{s}} \tilde{\phi}_{n_{\nes q}} \rangle^{n_p} -  \overline{ \langle G\nabla_{\!\bk_{s}}\tilde{\phi}_{n_p}\rangle^{n_{\nes q}}} \notag \\
&\!\!\!=\!\!\!&  (G\nabla\phi_{n_{\nes q}}, \phi_{n_p})-  \overline{(G\nabla\phi_{n_p}, \phi_{n_{\nes q}})} \;=\; (G\nabla\phi_{n_{\nes q}}, \phi_{n_p})- (G\nabla\phi_{n_p}, \phi_{n_{\nes q}}) ~\in\mathbb{R}^d, 
\end{eqnarray}
which vanishes identically for $p=q$. We also note that in such case $\bth_{qp}^{\tz}=-\bth_{pq}^{\tz}$, which makes $i\boldsymbol{A}^{\tz}$ skew-symmetric and~$\boldsymbol{A}^{\tz}$ itself Hermitian.

To complete the proof of the claim, we note that $\boldsymbol{A}^{\tz}$ admits real eigenvalues thanks to its Hermitian nature. Let $\tau$ be a non-zero eigenvalue, and let $\boldsymbol{v}$ denote the affiliated eigenvector so that $\boldsymbol{A}^{\tz}\boldsymbol{v}=\tau\boldsymbol{v}$. By conjugating this relationship, we find that $ \overline{\boldsymbol{A}^{\tz}} \overline{\boldsymbol{v}}=\tau \overline{\boldsymbol{v}}$, which yields $\boldsymbol{A}^{\tz} \overline{\boldsymbol{v}}=-\tau \overline{\boldsymbol{v}}$ thanks to the fact that $ \overline{\boldsymbol{A}^{\tz}}=-\boldsymbol{A}^{\tz}$. Hence, $-\tau$ is also an eigenvalue of $\boldsymbol{A}^{\tz}$, and $ \overline{\boldsymbol{v}}$ is the affiliated eigenvector.

%%--------------------------------------------------------------------------------------------------------------%
%\section*{References} \label{refs}
%%--------------------------------------------------------------------------------------------------------------%

\end{document}